\documentclass{article} 
   \usepackage{amssymb}  
\newcommand{\g}{\mbox{$\bf g$}}
\newcommand{\gal}{\mbox{$ \bf \mb{\bf g}_\alpha $}}

\newcommand{\h}{\mbox{\textbf{h}}}
\newcommand{\n}{\mbox{\textbf{n}}}
\newcommand{\np}{\mbox{$\textbf{n}^+$}}
\newcommand{\nm}{\mbox{$\textbf{n}^-$}}
\newcommand{\al}{\alpha}
\newcommand{\eps}{\epsilon}
\newcommand{\la}{\lambda}
\newcommand{\La}{\Lambda}
\newcommand{\de}{\delta}
\newcommand{\Th}{\Theta}
\newcommand{\gt}{\theta}
\newcommand{\mb}{\mbox}

\newcommand{\Mklz}[2]{\left\{\left.\;#1\;\right|\; #2\;\right\}}
\newcommand{\W}{\mbox{$\Delta$}}

\newcommand{\pW}{\mbox{$\Delta^+$}}
\newcommand{\rW}{\mbox{$\Delta_{re}$}}
\newcommand{\nrW}{\mbox{$\Delta_{re}^-$}}
\newcommand{\prW}{\mbox{$\Delta_{re}^+$}}
\newcommand{\iW}{\mbox{$\Delta_{im}$}}

\newcommand{\HI}{\mbox{$H_I$}}
\newcommand{\HR}{\mbox{$H_{Rest}$}}

\newcommand{\F}{\mathbb{F}}
\newcommand{\K}{\mathbb{K}}
\newcommand{\C}{\mathbb{C}}
\newcommand{\N}{\mathbb{N}}
\newcommand{\Nn}{\mathbb{N}_0}
\newcommand{\Q}{\mathbb{Q}}
\newcommand{\Qp}{\mathbb{Q}^+}

\newcommand{\R}{\mathbb{R}}

\newcommand{\Z}{\mathbb{Z}}
\newcommand{\We}{\mbox{$\cal W$}}

\newcommand{\Cq}{\mbox{$\overline{C}$}}
\newcommand{\FTq}{\mbox{$\overline{F}_\Theta$}}
\newcommand{\RkX}{\mbox{${\cal R}(X)$}}
\newcommand{\iB}[2]{\left(#1\mid#2\right)}
\newcommand{\kB}[2]{\left\langle \left\langle #1\mid #2
                       \right\rangle \right\rangle}
\newcommand{\kkB}[2]{\langle \langle #1\mid #2\rangle \rangle}                       
\newcommand{\kBl}{\left\langle \left\langle \;\ \mid \;\ \right\rangle 
                       \right\rangle}
\newcommand{\dcup}{\dot{\cup}}
\newcommand{\Vi}[1]{I\left(#1\right)}
\newcommand{\Vm}[1]{{\cal V}\left(#1\right)}
\newcommand{\grAdj}{\mbox{$gr$-$Adj$}\left(\,\bigoplus_{\La\in P^+}\,L(\La)\,
                    \right)}                   
              
\newcommand{\Endo}{\mbox{$End$}\left(\,\bigoplus_{\La\in P^+}\,L(\La)\,
                    \right)}
\newcommand{\EndoX}{\mbox{$End$}_X\left(\,\bigoplus_{\La\in P^+}\,L(\La)\,
                    \right)}
\newcommand{\grEndoX}{\mbox{$gr$-$End$}_X\left(\,\bigoplus_{\La\in P^+}\,L(\La)\,
                    \right)}
\newcommand{\Ze}[2]{\mb{Z}_{#1}\left(#2\right)}
\newcommand{\Mb}[2]{\mb{M}_{#1}^{\,\subset}\left(#2\right)}
\newcommand{\Mp}[2]{\mb{M}_{#1}^{\,\supset}\left(#2\right)}
\newcommand{\No}[2]{\mb{N}_{#1}\left(#2\right)}
\newcommand{\Gq}{\mbox{$ \overline{G} $}}
\newcommand{\Nq}{\mbox{$ \overline{N} $}}   
\newcommand{\Tq}{\mbox{$ \overline{T} $}}
\newcommand{\TD}{\mbox{$\hat{T}$}}
\newcommand{\ND}{\mbox{$\hat{N}$}}
\newcommand{\GD}{\mbox{$\hat{G}$}}
\newcommand{\WeD}{\mbox{$\hat{\We}$}}
\newcommand{\ve}[1]{\mbox{$\varepsilon\left(#1\right)$}} 
\newcommand{\FGq}{\mbox{$ \F\, [ \, \overline{\mbox{G}} \, ] $}}
\newcommand{\KK}[1]{\mbox{$\K\,[#1]$}}
\newcommand{\KXV}{\mbox{$\K_X[V]$}}
\newcommand{\FK}[1]{\mbox{$\F\,[#1]$}}
\newcommand{\CK}[1]{\mbox{${\mathbb C}\,[#1]$}}
\newcommand{\KA}{\mbox{$\K\,[\,A\,]$}}
\newcommand{\KB}{\mbox{$\K\,[\,B\,]$}}
\newcommand{\KC}{\mbox{$\K\,[\,C\,]$}}
\newcommand{\KD}{\mbox{$\K\,[\,D\,]$}}
\newcommand{\KCD}{\mbox{$\K\,[\,C\times D\,]$}}

\newcommand{\KAdjV}{\mbox{$\K\,[\,\mb{Adj}(V)\,]$}}

\newcommand{\FA}{\mbox{${\cal F}_A$}}
\newcommand{\FB}{\mbox{${\cal F}_B$}}
\newcommand{\FC}{\mbox{${\cal F}_C$}}
\newcommand{\FD}{\mbox{${\cal F}_D$}}

\newcommand{\FDg}{\mbox{${\cal F}_{D(g)}$}}
\newcommand{\FCD}{\mbox{$ {\cal F}_{C\times D}$}}
\newcommand{\FAdjV}{\mbox{${\cal F}_{Adj(V)}$}}

\newcommand{\FV}{\mbox{${\cal F}_V$}}
\newcommand{\Fa}[1]{ Fa(#1) }
\newcommand{\Proof}{\mbox{\bf Proof: }}
\newcommand{\End}{\mb{}\hfill\mb{$\square$}\\}
\newcommand{\Ende}{\mb{}\hfill\mb{$\square$}\vspace*{1ex}\\}
\newcommand{\cc}[1]{\mbox{cc}\left(#1\right)}
\newcommand{\ccp}[1]{\mbox{$\dot{\mbox{cc}}$}\left(#1\right)}
\newcommand{\Sp}{\mbox{Spec\,}}
\newcommand{\Spm}{\mbox{Specm\,}}
\newcommand{\FSpm}{\mbox{$F$-Specm\,}}
\newcommand{\Kern}{\mb{kernel\,}}

\newcommand{\ti}{\tilde}
\newcommand{\res}[1]{\!\mid_{#1}}
\newcommand{\CD}{C\times D}
\newcommand{\obot}{{\bigcirc\hspace*{-0.9em}\bot}}
\newcommand{\mapsot}{\begin{picture}(5,5)
\makebox(0,3){$\leftarrow$}\put(5,0){\line(0,1){5}}
\end{picture}}
\newcommand{\spann}[1]{\mbox{span}\left(#1\right)}
\begin{document}
\newtheorem{Theorem}{Theorem}[section]
\newtheorem{Def}[Theorem]{Definition}
\newtheorem{Prop}[Theorem]{Proposition}
\newtheorem{Cor}[Theorem]{Corollary}
\newtheorem{Lemma}[Theorem]{Lemma}
%
%
%
%
\title{An analogue of a reductive algebraic monoid, whose unit group is a Kac-Moody group}
\author{Claus Mokler\\\\ Universit\"at Wuppertal, Fachbereich Mathematik\\  Gau\ss stra\ss e 20\\ D-42097 Wuppertal, Germany\vspace*{1ex}\\ 
          mokler@math.uni-wuppertal.de}
\date{}
\maketitle
\begin{abstract}\noindent
By a generalized Tannaka-Krein reconstruction we associate to the admissible representions of the category ${\cal O}$ of a Kac-Moody algebra, and its 
category of admissible duals a monoid with a coordinate ring.\\
The Kac-Moody group is the Zariski open dense unit group of this monoid. The restriction of the coordinate ring to the Kac-Moody group is the 
algebra of strongly regular functions introduced by Kac and Peterson.\\
This monoid has similar structural properties as a reductive algebraic monoid. In particular it is unit regular, its idempotents related to the 
faces of the Tits cone. It has Bruhat and Birkhoff decompositions.\\
The Kac-Moody algebra is isomorphic to the Lie algebra of this monoid.
\end{abstract}
\vspace*{1ex}
\section*{Introduction}
%
%
The (minimal) Kac-Moody group $G$ constructed in \cite{KP1} by V. Kac and D. Peterson is a group analogue of a semisimple, simply connected algebraic 
group (without a coordinate ring, without a topology).\\
Kac and Peterson defined and investigated in \cite{KP2} the algebra of strongly regular functions $\CK{G}$ of a symmetrizable Kac-Moody group. 
This algebra has many properties in common with the coordinate ring of a semisimple, simply connected algebraic group. It is an integrally 
closed domain, even a unique fac\-tor\-iza\-tion domain. It admits a Peter and Weyl theorem, i.e.,
\begin{eqnarray*}
  \CK{G} &\cong & \bigoplus_{\La\in P^+} L^*(\La)\otimes L(\La)
\end{eqnarray*}
as $G\times G$-modules. But the following things, which hold in the non-classical case, are different:
\begin{itemize}
\item The group structure of $G$ does not induce a Hopf algebra structure on $\CK{G}$. There is no comultiplication due to the  
infinite-dimensionality of the highest and lowest weight representations. There is no antipode due to the inequivalence between highest 
and lowest weight representations.  
\item  Assigning to every element of $G$ its vanishing ideal, $G$ embeds in the set of 1-codimensional ideals of $\CK{G}$, which we 
denote by $\Spm\CK{G}$. But this map is not surjective.
\end{itemize}
Kac and Peterson posed the problem to determine the part 
\begin{eqnarray}\label{KPSpart} 
   \bigcup_{m\in\N}\;\;\bigcup_{\beta_1,\,\ldots,\,\beta_m\in\Delta_{re}}\overline{\,U_{\beta_1}\cdots U_{\beta_m}T\,}^{Specm}\;\,\subseteq \;\,
   \Spm\CK{G}\;\;.
\end{eqnarray}
(Here the elements of $G$ have been identified with the corresponding vanishing ideals. $T$ denotes a maximal torus given by the construction of the 
Kac-Moody group, and $U_{\beta_1},\,\ldots,\, U_{\beta_m}$ are root groups belonging to $T$.)\vspace*{1ex}\\
At the same time P. Slodowy wanted to connect the deformation theory of certain singularities with the Kac-Moody groups. Together with E. Looijenga 
he defined an adjoint quotient
\begin{eqnarray}\label{adjointquotient} 
  \chi\,:G &\to &\TD/\We
\end{eqnarray} 
as a map, \cite{Sl1}, \cite{Sl2}. In the non-classical case the target space $\TD/\We$ is a proper extension of $T/\We$, 
related to a certain partial compactification described in \cite{Loo}. (Such an extension is really necessary. This is indicated by the fact, that 
not every element of $G$ can be conjugated in a fixed Borel subgroup.)\\
Compared to its classical counterpart this adjoint quotient has bad properties. In the non-classical case there 
are 'missing' conjucagy classes in its fibres. Slodowy guessed that this quotient can be extended to a certain part $\GD\subseteq \Spm\CK{G}$.\\ 
There is an action of $G\times G$ on $\Spm\CK{G}$, extending the action on $G$. Slodowy conjectured
\begin{eqnarray}\label{conjectureSlodowy}
  \GD &=& G\,\cdot\,\TD\,\cdot\, G\;\;,
\end{eqnarray}
where $\TD$ should be realized as a torus embedding $T\subseteq \TD$ of nonfinite type, \cite{Sl2}. He also guessed the form of some elements 
$e(R(\Th))\in\TD$, $\Th$ special.\vspace*{1ex}\\
He discussed this with Peterson. Some time later Peterson announced a proof of this conjecture, $\GD$ the part (\ref{KPSpart}) of $\Spm\CK{G}$, 
together with some other structural properties of $\GD$, \cite{P}. The most important are:
\begin{itemize}
  \item There should be a monoid structure on $\GD$.
  \item Every idempotent of $\GD$ is conjugate to one of the elements $e(R(\Th))$ by an element of $G$. 
  \item $\GD = U_{\pm}\ND U_{\pm}$ where $\ND:=N\Mklz{e(R(\Th))}{\Th\mb{ special }}N$ is a submonoid.
\end{itemize}
But Peterson didn't give any proof.\vspace*{1ex}\\
About ten years later, as a student of P. Slodowy, the task of my PhD was to investigate all these conjectures and claims, \cite{M}.\vspace*{1ex}\\
Important for these investigations was the paper \cite{Kas} of M. Kashiwara, in which he had studied the whole spectrum $\Sp\CK{G}$ in connection with 
his infinite dimensional algebraic geometric approach to the flag varieties of Kac-Moody groups. For his construction 
of the flag variety he used an open subscheme of the full spectrum $\Sp\CK{G}$, which has a countable covering by big cells. He factored this subscheme 
by an action of the formal Borel subgroup of the formal Kac-Moody group $G_f$.\vspace*{1ex}\\
There had also been another development relevant for these investigations:\vspace*{0.5ex}\\
A reductive linear algebraic monoid $(M,\CK{M})$ is a connected linear algebraic monoid, whose principal open unit group $G$ is reductive. The linear 
algebraic monoid $M(n,\C)$ of $n\times n$-matrices is a typical example. Note that 
there are no nontrivial reductive linear algebraic monoids with semisimple unit group.\\
Reductive linear algebraic monoids have been studied carefully by M. Putcha and L. Renner. Many of their results are presented in the book \cite{Pu}.\\
Such a monoid $M$ is unit regular. Every idempotent of $M$ is $G$-conjugate to an idempotent in the closure $\overline{T}$ of a 
maximal torus $T$ of $G$. It has Bruhat and Birkhoff decompositions, the Weyl group $N/T$ replaced by the Renner monoid $\overline{N}/T$.\\
The closure $\Tq$ of the torus $T$ is an affine toric variety. The Weyl group $\We$ acts on $\Tq$, extending its action on $T$.  
In \cite{Re} Renner showed that a normal reductive linear algebraic monoid $M$ with zero and one-dim\-en\-si\-onal 
center is classified by its unit group $G$, and the toric variety $\overline{T}$, equipped with the Weyl group action.\\
In \cite{Vi} E. Vinberg remarked that the condition of a one-dim\-en\-si\-onal center is not really essential. He gave a representation theoretic 
approach to the classification problem.\\\\
To investigate part (\ref{KPSpart}) of $\Spm\CK{G}$, which we call Kac-Peterson-Slodowy part, to prove the conjectures of P. Slodowy, and 
the claims of D. Peterson we shall proceed as follows:\vspace*{1ex}\\
Since the Kac-Moody group $G$ acts faithfully on the sum $\bigoplus_{\La\in P^+} L(\La)$, it can be identified with with a subgroup of 
the monoid of linear endomorphisms $End(\bigoplus_{\La\in P^+} L(\La))$. In Section 1 we define a monoid $\GD$, which is generated by $G$ and certain 
idempotents $e(R(\Th))\in End(\bigoplus_{\La\in P^+} L(\La))$, $\Th$ special.\\
For this monoid we prove the conjecture (\ref{conjectureSlodowy}) of P. Slodowy, and all the claims of D. Peterson on the structure. But we obtain 
even more results. We show 
\begin{eqnarray*}
  \GD &=& \dot{\bigcup_{\hat{\sigma}\in\hat{\cal W}}} B^\pm \hat{\sigma} B^\pm\;\;,
\end{eqnarray*}
and investigate the monoid $\WeD=\ND/T$, which we call Weyl monoid. We prove a theorem, which allows to reduce any computation in $\GD$ to 
computations in the Kac-Moody group $G$, and in the face lattice $\RkX$ of the Tits cone $X$, equipped with its action of the Weyl group. 
We show how to construct the monoid $\GD$ from the system of root groups data $(U_\al)_{\al\in\Delta_{re}}$, $T$, and from the monoid 
$(\,\RkX\,,\,\cap\,)$.\vspace*{1ex}\\
It is not too difficult to see, that $\GD$ embeds in $\Spm\CK{G}$ in a natural way, its image contained in the Kac-Peterson-Slodowy part. 
But I could not find a direct proof, showing equality.\\
Also it is surprising, that there exists such a monoid. It would be desirable to gain a better understanding where this monoid comes from, and 
to make the analogy with a reductive algebraic monoid even more tight.\vspace*{1ex}\\ 
A semisimple, simply connected algebraic group can be obtained from its category of rational representations by Tannaka-Krein reconstruction.\\ 
It is also possible to use for this reconstruction the corresponding category of representations of its semisimple Lie algebra. Because we assumed 
'simply connected', this category coincides with the category of finite dimensional representations of the semisimple Lie algebra.\vspace*{1ex}\\
It is important to observe, that there are several natural generalizations of this category for a Kac-Moody algebra $\g$:
\begin{itemize}
\item The tensor category of admissible modules of $\g$. (These modules are used for the construction of the Kac-Moody group $G$.)
\item The tensor category ${\cal O}_{adm}$ of admissible modules of the category $\cal O$. This category is close to the classical case, because it 
keeps the complete reducibility theorem, i.e., every module of ${\cal O}_{adm}$ is isomorphic to a direct sum of irreducible highest weight modules 
$L(\La)$, $\La\in P^+$. 
\item A category, which is build in a similar way as ${\cal O}_{adm}$, but where we have a complete reducibility theorem using the lowest weight 
modules $L^*(\La)$, $\La\in P^+$.
\end{itemize}
In a way similar to the Tannaka-Krein reconstruction, a suitable category of representations of a Lie algebra, together with a suitable category of 
duals determine a monoid with coordinate ring. This monoid is the biggest monoid, acting reasonably on the representations, compatible with the duals. 
The coordinate ring is a coordinate ring of matrix coefficients. But the multiplication of the monoid does not have to induce a comultiplication of 
the coordinate ring.\vspace*{1ex}\\
In Section 4 we formulate this generalized Tannaka-Krein reconstruction for the category ${\cal O}_{adm}$ and its category of 
admissible duals.\\
In Section 5 we show that $\GD$ is the monoid determined by these categories, which is one of the main results of this paper. The Kac-Moody group $G$ is 
the Zarsiki open dense unit group of $\GD$. The algebra of strongly regular functions is isomorphic to the coordinate ring $\CK{\GD}$ by the restriction 
map. In this way, the algebra of strongly regular functions is really the coordinate ring of a monoid.\\
Note that this monoid is a purely infinite dimensional phenomenon. In the classical case, it reduces to a semisimple, simply connected algebraic group.
\vspace*{1ex}\\
It is quite a long way to prove this result. We use an easy infinite dimensional algebraic geometry, which is developed in Section 3. The 
monoids obtained by a generalized Tannaka-Krein reconstruction are weak algebraic monoids in this algebraic geometric setting.\vspace*{1ex}\\
As remarked above, it is not too difficult to see, that the part of $\Spm \CK{G}$ corresponding to $\GD$ is contained in the Kac-Peterson-Slodowy 
part. Using the characterization of $\GD$ by the generalized Tannka-Krein reconstruction, we show in Section 5, that the Kac-Peterson-Slodowy part 
is contained in the part of $\Spm \CK{G}$ corresponding to $\GD$. Therefore we have equality.\\ 
Furthermore we show in Section 5, that the Zariski closure $\Tq=\TD$ is an affine toric variety of nonfinite type. We show $\Nq=\ND$, which makes 
the analogy of the Weyl monoid $\WeD=\ND/T=\Nq/T$ with a Renner monoid even more tight.\\ 
Although the monoid $\GD$ does not act on the flag varieties, we show that it acts on the affine cones of the flag varieties.\vspace*{1ex}\\
In Section 6 we show that the Kac-Moody algebra is isomorphic to the Lie algebra of $\GD$.\vspace*{1ex}\\
In Section 1 we quickly review Kac-Moody algebras, Kac-Moody groups and the algebra of strongly regular functions, merely to introduce our 
notations. We also give an easy generalization of affine toric varieties.\\\\
This preprint is a revised version of the preprint 'A monoid completion of a Kac-Moody group', which had been submitted to a journal three years ago, but 
did not get valued.\\ 
The monoid $\GD$ can be characterized as a closure in a very big space. To investigate this closure has been motivated in the preprint 
'A monoid completion of a Kac-Moody group' by the approach of Vinberg to the classification problem of reductive algebraic monoids.\\
I found the interpretation by the generalized Tannaka-Krein reconstruction after my PhD. I included it here, because it is more beautiful.
\vspace*{1ex}\\ 
After this long time, also many other things in relation to this monoid have been investigated:\\
The $\C$-valued points of the algebra of strongly regular functions have been determined and investigated in \cite{M2}. The results on the 
$G_f\times G_f$-orbits of $\Spm\CK{G}$ (closures, big cells, stratified transversal slices) can be transfered to the $G\times G$-orbits of $\GD$, for 
which they are easier to prove. The closures of the Bruhat and Birkoff cells have been determined, and will be given in \cite{M3}, as well as an 
extension of the length function from the Weyl group to the Weyl monoid.\\
The generalized Tannaka-Krein reconstruction has been investigated in general. It was natural to determine the monoid, which corresponds to the 
category ${\cal O}_{adm}$ and its category of full duals. This is a monoid $\widehat{G_f}$, which is build in a similar way as the monoid $\GD$, but 
the minimal Kac-Moody group $G$ is replaced by the formal Kac-Moody group $G_f$, \cite{M4}.\\
Presumably the Kac-Moody group $G$ is the monoid associated to the category of admissible modules and its category of admissible duals. This is 
difficult to show, because there is no complete reducibility theorem.\vspace*{1ex}\\
The following investigations have already been started:\\
Extending the adjoint quotient (\ref{adjointquotient}) of Looijenga and Slodowy to an adjoint quotient $\chi:\GD\to \TD/\We$.\\
Extending the monoid $G^\geq\subseteq G$ of totally non-negative elements of Lusztig to a monoid $\hat{G}^\geq\subseteq \hat{G}$. 
Generalizing the characterizations of Fomin and Zelevinsky.\\
Defining and investigating the analogues of normal reductive algebraic monoids in the Kac-Moody case.
%
%
%
\newpage\tableofcontents
%
%
%
%
\newpage\section{Preliminaries}                                                                                                                   
%
%
%
We denote by $\N=\Z^+$, $\Q^+$, resp. $\R^+$ the sets of strictly positive numbers of $\Z$, $\Q$, resp. $\R\,$,
and the sets $\N_0=\Z^+_0$, $\Q^+_0$, $\R^+_0$ contain in addition the zero.\vspace*{0.5ex}\\
In the whole paper $\F$ is a field of characteristic 0, and $\K$ an arbitrary field with $|\K\,|=\infty\,$. Their
unit groups are denoted by $\F^\times$ and $\K^\times $.
%
%
%
\subsection{Kac-Moody algebras, Kac-Moody groups and the algebra of strongly regular functions\label{section2}}                           
%
%
%
In this subsection we collect some basic facts about Kac-Moody algebras, Kac-Moody groups and the algebra of strongly regular functions which 
are used later, and introduce our notation.\\
The Kac-Moody group given in \cite{KP1}, \cite{KP3} corresponds to the derived Kac-Moody algebra. We work 
with a slightly enlarged group corresponding to the full Kac-Moody algebra as in \cite{Ti1}, \cite{MoPi}.\\
All the material stated in this subsection about Kac-Moody algebras can be found in the books \cite{K2} (most results also valid for a field of 
characteristic zero, with the same proofs), \cite{MoPi}, about Kac-Moody groups in \cite{KP1}, \cite{KP3}, \cite{MoPi}, and about the algebra of 
strongly regular functions in \cite{KP2}.\vspace*{2ex}\\ 
%
%
%
{\bf Optimal realizations:}
%
%
%
Starting point for the construction of a Kac-Moody algebra and its associated simply connected Kac-Moody group is a {\it generalized Cartan matrix} 
$A=(a_{ij})\in M_{n}(\Z)$, which we also assume to be symmetrizable. Denote by $l$ the rank of $A$, and set $I:=\{1,2,\ldots n\}$.\\
A simply connected minimal free root base for $A$, which we call an {\it optimal realization} for short, consists of:\\ 
$\bullet$ Dual free  $\mathbb{Z}$-modules $H,P$ of rank $2n-l$.\\
$\bullet$ Linear independent sets $\Pi^\vee =\{h_1,\,\ldots,\, h_n\}\subseteq H $, 
$\Pi=\{\al_1,\,\ldots,\,\al_n\}\subseteq P$, such that $\al_i(h_j)=a_{ji}\,$, $i,j=1,\,\dots,\, n$. \\
$\bullet$ Furthermore $\HI:=Q^\vee:= \mathbb{Z}\mb{-span}\{h_1,\,\ldots,\, h_n\}$ is saturated in $H$, 
which means that for all $n\in\mathbb{N}\,$, $x\in H$ we have: $\quad nx\in H_I \;\Rightarrow\;  x\in H_I$\vspace*{0.5ex}\\
$P$ is called the {\it weight lattice}, $Q:=\Z\mb{-span}\Mklz{\al_i\,}{\,i\in I}$ the {\it root lattice}, and $Q^\vee$ 
the {\it coroot lattice}. Set $Q^\pm_0:=\Z^\pm_0\mb{-span}\Mklz{\al_i\,}{\,i\in I}$, and $Q^\pm:=Q^\pm_0\setminus\{0\}$.\vspace*{0.5ex}\\
We fix a complement $H_{rest}$ of $H_I$ in $H$. This complement determines a system of {\it fundamental dominant weights} $\La_1,\,\ldots,\,\La_n$ by
\begin{eqnarray*}
     \La_i(h_j) \;:=\;\de_{ij} \quad (j=1,\ldots n)  &\;\:,\:\;& \La_i(h) \;:=\; 0  \quad \;\;(h\in H_{rest})\;\;.
\end{eqnarray*} 
We extend $h_1,\,\ldots,\, h_n\in H_I$ with elements $h_{n+1},\,\ldots,\, h_{2n-l}\in H_{rest}$ to a base of $H$, and 
extend $\La_1,\,\ldots,\, \La_n$ to the corresponding dual base $\La_1,\,\ldots,\, \La_{2n-l}\,$.\vspace*{1ex}\\
%
%
%
%
%
%
%
{\bf The Weyl group and the Tits cone:} 
%
%
%
%
Identify $H$ and $P$ with the corresponding sublattices of the following vector spaces over $\F\,$:
\begin{eqnarray*}
   \h  \;\,:=\;\,  \h_\F    \;\,:=\;\,   H \otimes_{\mathbb Z} \F   &\quad,\quad &
   \h^* \;\,:=\;\, \h^*_\F  \;\,:=\;\,  P \otimes_{\mathbb Z} \F\;\;.
\end{eqnarray*}
$\h^*$ is interpreted as the dual of $\h$. Order the elements of $\h^*$ by $\la\leq\la'$ if and only if $\la'-\la\in Q_0^+$.\\ 
Choose a symmetric matrix $B\in M_n (\Q)$, and a diagonal matrix 
$D=\mb{diag}(\eps_1,\,\ldots,\, \eps_n)$, $\eps_1,\,\ldots,\,\eps_n\in\Qp$, such that $A=DB$.
Define a nondegenerate symmetric bilinear form on $\h$ by:
\begin{eqnarray*}
 \iB{h_i}{h} \;=\; \iB{h}{h_i} \;:=\; \al_i(h)\,\eps_i  &\qquad & i\in I\,,
 \quad h\in\h\;,\\
 \iB{h'}{h''} \;:=\; 0 \qquad\qquad &\qquad & h',\,h'' \in \h_{rest}:=H_{rest}\otimes\F \;.
\end{eqnarray*}
Denote by $\nu:\,\h\to\h^*$ the corresponding isomorphism. Denote by $\iB{\;}{\;}$ the in\-duced 
nondegenerate symmetric form on $\h^*$.\vspace*{1ex}\\
The {\it Weyl group} $\We=\We(A)$ is the Coxeter group with generators $\sigma_i$, $i\in I$, and
relations
\begin{eqnarray*}
       \sigma_i^2 \;=\; 1 \qquad (i\in I)    \;\:&\;,\;&\;\:
       {(\sigma_i\sigma_j)}^{m_{ij}} \;=\; 1 \qquad (i,j\in I,\,i\ne j)\;\;.
\end{eqnarray*}
The $m_{ij}$ are given by:
    $\quad \begin{tabular}{c|ccccc}
      $a_{ij}a_{ji}$ & 0 & 1 & 2 & 3 &  $\geq$ 4  \\[0.5ex] \hline 
         $m_{ij}$    & 2 & 3 & 4 & 6 &  no relation between $\sigma_i$ and $\sigma_j$ 
     \end{tabular}$\vspace*{1ex}\\
The Weyl group $\We$ acts faithfully and contragrediently on $\h$ and $\h^*$ by 
\begin{eqnarray*}
  \sigma_i h \;:=\; h - \al_i\,(h) h_i & \qquad i\in I,\quad h\in \h \;\;,\\
  \sigma_i \la \;:=\; \la - \la(h_i)\,\al_i & \qquad i\in I,\quad \la\in \h^*\;\;,
\end{eqnarray*}
leaving the lattices $Q^\vee$, $H$, $Q$, $P$, and the forms invariant.\\
The set $\Delta_{re}:=\We\Mklz{\al_i}{i\in I}\subseteq Q$ is called the set of {\it real roots}, and 
$\Delta_{re}^\vee:=\We\Mklz{h_i}{i\in I}\subseteq Q^\vee$ the set of {\it real coroots}. The map 
$\al_i\mapsto h_i\,$, $i\in I$, can be extended to a $\We$-equivariant bijection 
$\al \mapsto  h_\al=\frac{2\nu^{-1}(\al)}{\iB{\al}{\al}}\,$.\vspace*{1ex}\\
To illustrate the action of $\We$ on $\h^*_\R$ geometrically, for $J\subseteq I$ set
\begin{eqnarray*}
  F_J &:=& \Mklz{\la\in\h^*_\R}{\la(h_i)\,=\,0 \;\mb{ for }\; i\in J\,,\;\;\; \la(h_i)\,>\,0\; \mb{ for }\;i\in I\setminus J}\;\;,\\
  \overline{F}_J &:=& \Mklz{\la\in\h^*_\R}{\la(h_i)\,=\,0\; \mb{ for }\;i\in J \,,\;\;\;\la(h_i)\,\geq\,0\; \mb{ for }\;i\in I\setminus J}\;\;. 
\end{eqnarray*}  
$\overline{F}_J$ is a closed, finitely generated, convex cone with relative interior $F_J$, and
\begin{eqnarray*}
  \overline{F}_J  &=& \dot{\bigcup_{I\,\supseteq\, K\,\supseteq\, J}} F_K \:\;.
\end{eqnarray*}
Call $\overline{C} := \overline{F}_\emptyset = \Mklz{\la\in\h^*_\R}{\la(h_i)\,\geq\,0\; \mb{ for }\;i\in I}$ the {\it fundamental chamber}, and 
$C := F_\emptyset = \Mklz{\la\in\h^*_\R}{\la(h_i)\,>\,0\; \mb{ for }\;i\in I}$ the {\it open fundamental chamber}. 
Set $c := F_I =\overline{F}_I = \Mklz{\la\in\h^*_\R}{\la(h_i)\,=\,0 \,\mb{ for }\,i\in I}$.\\ 
The {\it Tits cone} $X:=\We\,\overline{C}\,$ is a convex $\We$-invariant cone 
with edge $c$, and fundamental region $\overline{C}=\bigcup_{J\subseteq I}F_J$. The parabolic subgroup $\We_J$ of 
$\We$ is the stabilizer of any element $\la\in F_J$, as well as the stabilizer of $F_J$ as a whole.\\ 
The set $\Mklz{\sigma F_J}{\sigma\in\We\,,\,J\subseteq I}$ gives a $\We$-invariant partition of $X$. We call $\sigma F_J$ a facet of 
type $J$.\vspace*{1ex}\\ 
Set $P^+:=P\cap\Cq$, and $P^{++}:=P\cap C$.\vspace*{1ex}\\
%
%
%
{\bf The Kac-Moody algebra:} 
%
%
%
%
The {\it Kac-Moody algebra} $\g=\g(A)$ is the Lie algebra over $\F$ generated by the abelian Lie algebra $\h$, and
$2n$ elements $e_i$, $f_i$, ($i\in I$), with the following relations, which 
hold for any $i,j \in I$, $h \in \h$: 
  \begin{eqnarray*}
    \left[ e_i,f_j \right] \,=\,   \delta_{ij} h_i   \;\;,\;\; 
    \left[ h,e_i \right]   \,=\,   \al_i(h) e_i    \;\;,\;\;  
    \left[ h,f_i \right]   \,=\,  -\al_i(h) f_i\;\;,     \\
    \left(ad\,e_i\right)^{1-a_{ij}}e_j  \,=\, \left(ad\,f_i\right)^{1-a_{ij}}f_j  \,=0  \qquad (i\neq j)\;\;.\qquad
  \end{eqnarray*}
The {\it Chevalley involution} $*$ of $\g$ is the involutive anti-automorphism determined by $e_i^*=f_i$, $f_i^*=e_i$, $h^*=h$, where $i\in I$ and 
$h\in\h$.\vspace*{1ex}\\
The space $\h$ and the elements $e_i$, $f_i$, $(i\in I)$, can be identified with their images in $\g$. The nondegenerate symmetric bilinear form 
( $|$ ) on $\h$ can be uniquely extended to a nondegenerate symmetric invariant bilinear form ( $|$ ) on $\g$.\\
We have the {\it root space decomposition}
\begin{eqnarray*}
  \g=\bigoplus_{\al \in \h^*}\g_{\al} \quad \mb{where} \quad 
  \g_\al := \Mklz{x\in \g}{[h,x]=\al(h)\,x\;\mb{ for all }\;h\in \h}\;\;.            
\end{eqnarray*}
In particular $\g_0 = \h$, $\g_{\al_i}=\F e_i$, and $\g_{-\al_i}=\F f_i$, $i\in I$.\\
The set of {\it roots} $\W:=\Mklz{\al\in\h^*\setminus\{0\}}{\g_\al\ne 0}$ is invariant under the Weyl group and spans the {\it root lattice} $Q$. 
We have $\rW\subseteq \W$, and $\iW:=\W\setminus\rW$ is called the set of {\it imaginary roots}.\\
$\W$, $\rW$, and $\iW$ decompose into the disjoint union of the sets of 
{\it positive} and {\it negative roots} $\W^\pm:=\W\cap Q^\pm$, $\rW^\pm:=\rW\cap Q^\pm$, $\iW^\pm:=\iW\cap Q^\pm$, and we have 
$\W^\pm = -\W^\mp$, $\rW^\pm = -\rW^\mp$, $\iW^\pm = -\iW^\mp$.\\ 
The decomposition $\g=\nm \oplus \h \oplus \np$, where $\n^\pm :=\bigoplus_{\al\in {\Delta}^\pm} \g_\al$, is called the {\it triangular decomposition}.
\vspace*{1ex}\\ 
For a real root $\al$ the subalgebra $\g_\al\oplus[\g_\al,\g_{-\al}]\oplus\g_{-\al}$ of $\g$, which is isomorphic to $sl(2,\F)$, is denoted by 
${\bf s}_\al$.\vspace*{1ex}\\
The derived Lie algebra $\g$ is given by $\g'=\bigoplus_{\al\in\Delta}\g_\al\oplus \h'$, where 
$\h'$ is spanned by the elements $h_i$, $i\in I$.\vspace*{1ex}\\
%
%
%
%
{\bf The category $\cal O$, and the irreducible highest weight representations:}
%
%
%
The category $\cal O$ is a tensor subcategory of the category of $\g$-modules. Its ob\-jects consist of the $\g$-modules $V$, 
with the following properties:\\
(1) $V$ is $\h$-diagonalizable with finite dimensional weight spaces.\\
(2) There exist finitely many elements $\la_1,\,\ldots,\,\la_m\in\h^*$, such that the set of weights of $V$, which we denote by $P(V)$, 
is contained in the union 
$\bigcup_{1=1}^m D(\la_i)$, where $D(\la_i):=\Mklz{\la\in\h^*}{\la\leq \la_i}$.\\
The morphisms of $\cal O$ are morphisms of $\g$-modules.\vspace*{1ex}\\ 
For every $\La\in\h^*$ there exists, unique up to isomorphism, an irreducible representation $(L(\La),\pi_\La)$
of $\g$, with a nonzero element $v_\La$ satisfying
\begin{eqnarray*}
   \pi_\La(\n^+)v_\La \;=\;0   &\;\:,\:\;&   \pi_\La(h)v_\La\;=\;\La(h)v_\La \quad (h\in\h)\;.
\end{eqnarray*}
It is is an object of $\cal O$, and is called an irreducible highest weight representation. We denote its set of weights by 
$P(\La)$.\vspace*{0.5ex}\\
Any such representation carries a nondegenerate symmetric bilinear form $\kBl:L(\La)\times L(\La)\to\F$,
which is contravariant, i.e. $\kB{v}{xw}=\kB{x^*v}{w}$ for all $v,w\in L(\La)$, $x\in\g$. This form is unique 
up to a nonzero multiplicative scalar.\vspace*{1ex}\\  
%
%
%
%
{\bf The Kac-Moody group:} 
%
%
%
%
To construct the Kac-Moody group call a representation $(V,\pi)$ of $\g$ {\it admissible} if:\\
$\quad$(1) $V$ is $\h$-diagonalizable with set of weights $P(V)\subseteq P$.\\ 
$\quad$(2) $\pi(x)$ is locally nilpotent for all $x\in\g_\al$, $\al\in\rW$.\\
(Note: If the generalized Cartan matrix is degenerate, then this definition is slightly different from 'integrable', because we demand 
$P(V)\subseteq P$. This guarantees, that we can integrate the full Cartan subalgebra $\h$ to an algebraic geometric torus.)\\ 
Examples are the irreducible highest weight representations ($L(\La)$, $\pi_\La$), $\La\in P^+$, 
and the adjoint representation ($\g\,$, $ad\,$).\vspace*{1ex}\\
We denote by ${\cal O}_{adm}$ the full subcategory of the category $\cal O$, whose objects are admissible modules. Every object of this category is 
isomorphic to a direct sum of irreducible highest weight modules $L(\La)$, $\La\in P^+$.\vspace*{1ex}\\
Let $\ti{G}$ be the free product of the additive groups $\g_\al\,$, $\al\in\rW$, and the torus Hom$((P,+),(\F^\times,\,\cdot\,))$. 
Due to (1), (2), and the coproduct property of $\ti{G}$, we get for any admissible representation $(V,\pi)$ a homomorphism $\ti{\pi}:\ti{G}\to GL(V)$, 
mapping $x_\al\in\g_\al$ to $\exp(\pi(x_\al))$, and $\eps\in\mb{Hom}(P,\F^\times)$ to an element 
$t_\pi(\eps)$, defined by
\begin{eqnarray*}
  t_\pi(\eps)v_\la &:=& \eps(\la)v_\la \quad,\quad v_\la\in V_\la\;,\;\la\in P(V)\;.
\end{eqnarray*}
Let $\ti{N}$ be the intersection of all kernels of homomorphisms $\ti{\pi}'$ corresponding to admissible 
representations $(V,\pi')$. The {\it Kac-Moody group} is defined as $G:=G(A):=\ti{G}/\ti{N}$, and due to its 
definition $\ti{\pi}:\ti{G}\to GL(V)$ factors to a representation $\Pi:G\to GL(V)$. 
Corresponding to ($\g\,$, $ad\,$) we get the adjoint representation $Ad: G\to Aut(\g)$. \vspace*{0.5ex}\\
By composing the injection of $\g_\al$ into $\ti{G}$ with the projection onto $G$ we get an injective 
homomorphism $\exp: \g_\al \to G$, its image $U_\al$ is called the {\it root group} belonging to $\al\in\rW$. Similar we get an 
injective homomorphism $t: \mb{Hom}(P,\F^\times)\to G$, its image denoted by $T$. Identifying $H\otimes_\Z\F^\times$ 
with $\mb{Hom}(P,\F^\times)$, the torus $T$ can be described by the isomorphism\vspace*{1ex}\\
\hspace*{8em}$\begin{array}{ccc}
 H\otimes_\Z\F^\times &\to & \quad\:T \\
 \sum_{\tau}h_\tau\otimes s_\tau &\mapsto & \prod_\tau t_{h_\tau}(s_\tau) 
\end{array}\;\;$,\vspace*{1ex}\\
where  $\,t_h(s)v_\la:=s^{\la(h)}v_\la$, $v_\la\in V_\la$, $\la\in P(V)$. Set $t_i(s):=t_{h_i}(s)$, 
$i=1,\,\ldots,\,2n-l$, $s\in\F^\times$.\vspace*{1ex}\\
The derived group $G'$, which is identical with the Kac-Moody group as defined in \cite{KP1}, is generated by 
the root groups $U_\al$, $\al\in \rW$. We have $G = G'\rtimes T_{rest}$, where $T_{rest}$ is the 
subtorus of $T$ generated by the elements $t_i(s)$, $i=n+1,\,\ldots,\, 2n-l$, $s\in \F^\times$.\vspace*{0.5ex}\\
The Kac-Moody group acts faithfully on $\bigoplus_{\La\in P^+}L(\La)$, and  
$\bigoplus_{i=1,\,\ldots,\, 2n-l}L(\La_i)$.\vspace*{0.5ex}\\
The Chevalley involution $*:\,G \to G$ is the involutive anti-isomorphism determined by
\begin{eqnarray*}
  \exp(x_\al)^*\;:=\;\exp(x_\al^*)  \quad  (x_\al\in\g_\al\,,\;\al\in\rW )  & \;\:,\:\; &  t^* \;:=\;t \quad (t\in T)\;.
\end{eqnarray*}
It is compatible with any nondegenerate symmetric contravariant form $\kBl$ on any of the modules $L(\La)$, $\La\in P^+$, i.e., 
$\kB{v}{gw}=\kB{g^*v}{w}$ for all $v,\,w\in L(\La)$, $g\in G$.\vspace*{1ex}\\
The Kac-Moody group has the following important structural properties:\vspace*{1ex}\\
1) Let $\al\in \prW$ and $x_\al\in\g_{\al}$, $x_{-\al}\in\g_{-\al}$, such that $[x_\al,x_{-\al}]=h_\al$. There exists an injective 
homomorphism of groups
\begin{eqnarray*}
  \phi_\al:\;\mb{SL}(2,\F) &\to & G
\end{eqnarray*}
with $
 \phi_\al\left(\begin{array}{cc}
  1 & s\\
  0 & 1
  \end{array}\right) \;:=\; \exp(s x_\al)  \;,\;
   \phi_\al\left(\begin{array}{cc}
  1 & 0\\
  s & 1
  \end{array}\right) \;:=\; \exp(sx_{-\al})$, $s\in \F^\times$, its image denoted by $G_\al$.\vspace*{1ex}\\ 
2) For $s\in \F^\times$ set $
  n_\al(s) = \phi_\al\left(\begin{array}{cc}
        0      & s  \\
  -\frac{1}{s} & 0
  \end{array}\right)$, and for short set $n_i(s):=n_{\al_i}(s)$, $i=1,\,\ldots,\, n$. Note that $t_{h_\al}(s) =  \phi_\al\left(\begin{array}{cc}
        s      &     0         \\
        0      & \frac{1}{s}
  \end{array}\right)= n_\al(s) n_\al(1)^{-1} $.\vspace*{1ex}\\
Denote by $N$ the subgroup generated by $T$ and $n_\al(1)$, $\al\in \Delta_{re}$. Let  
$B^\pm$ be the subgroups generated by $T$ and $U_\al$, $\al\in\Delta_{re}^\pm $, and let 
$U^\pm$ be the subgroups generated by $U_\al$, $\al\in\Delta_{re}^\pm $.\\ 
Then $(\,G\,,\,(U_\al)_{\al\in\Delta_{re}}\,,\,T\,)$ is a root groups data system, compare \cite{KP3}, Proposition 4.7,
and \cite{Ti}, leading to the twinned BN-pairs $B^\pm$, $N$, which have the property 
$B^+\cap B^- = B^+\cap N = B^-\cap N = T$. The Weyl group $\We$ can be identified with the common Coxeter group $N/T$ by the isomorphism 
$\kappa:\,\We \to  N/T$, which is given by $\kappa(\sigma_\al):=n_\al(1) T$, $\al\in \rW $.\\
In particular the twinned BN-pairs lead to the Bruhat and Birkhoff decompositions: 
\begin{eqnarray*}
  G &=& \dot{ \bigcup_{\sigma\in{\cal W}}} B^\epsilon\sigma B^\delta \qquad
  (\epsilon,\delta)\;=\;\underbrace{(+,+)\;,\;(-,-)}_{Bruhat}\;,\;\underbrace{(+,-)\;,\;(-,+)}_{Birkhoff}
\end{eqnarray*}
We denote an arbitrary element $n\in N$ with $\kappa^{-1}(nT)=\sigma\in\We$  by $n_\sigma$.
If $(V,\pi)$ is an admissible $\g$-module, its set of weights $P(V)$ is $\We$-invariant, 
and $ n_\sigma V_\la =V_{\sigma \la}$, $\la\in P(V)$.
We have $n_\sigma U_\al n_\sigma^{-1} = U_{\sigma\al}$, $\al\in\rW$. \vspace*{1ex}\\
3) The standard parabolic subgroups 
\begin{eqnarray*}
  P_J^\pm \;\,=\,\; B^\pm \We_J B^\pm &\quad,\quad &  J\subseteq I\;\;, 
\end{eqnarray*}
admit the Levi decompositions 
\begin{eqnarray*}
P_J^\pm \;\,=\,\;\ti{G}_J \ltimes (U^J)^\pm \;\;.
\end{eqnarray*} 
Here $\ti{G}_J$ is the group generated by the groups $U_{\pm\al_j}$, $j\in J$, and by $T$. This group coincides with the intersection 
$P_J^+ \cap P_J^-$. The groups $(U^J)^\pm$ are the normal subgroup of $U^\pm$, or of $(P_J)^\pm$, generated by $U_\al$, $\al\in\W^\pm_{re}\setminus
\sum_{j\in J}\Z\,\al_j$. They coincide with $\bigcap_{\sigma\in{\cal W}_J}\sigma U^\pm\sigma^{-1}$.\\
Usually the upper indices $+$ of all these groups are omitted.\vspace*{1ex}\\
We set $\ti{N}_J:=\ti{G_J}\cap N$. This group is the group generated by the elements $n_j(1)$, $j\in J$, and by $T$.\vspace*{1ex}\\
%
%
%
%
{\bf The algebra of strongly regular functions:}
%
%
%
%
%
The irreducible lowest weight module $L^*(\La)$ of lowest weight $-\La$ can be realized as the submodule 
$\bigoplus_{\la\in P(\La)}(L(\La)_\la)^*$ of the full dual module $L(\La)^*$ of $L(\La)$, $\La\in\h^*$.\\
For $\La\in P^+$, $v\in L(\La)$, and $\phi\in L^*(\La)$ call the function $\ti{f}_{\phi v}:\,G\to\F$ given by $\ti{f}_{\phi v}(g):=\phi(gw)$, 
$g\in G$, a {\it matrix coefficient} of $G$. The algebra $\FK{G}$ generated by all such matrix coefficients is called the {\it algebra of 
strongly regular functions} on $G$.\vspace*{0.5ex}\\
$\FK{G}$ is an integrally closed domain, and admits a Peter and Weyl theorem: Define an action of $G\times G$ 
on $\FK{G}$ by $\left((g,h)f\right)(x):=f(g^{-1}xh)$, $g,x,h\in G$, $f\in\FK{G}$. Then the map 
$\bigoplus_{\La\in P^+} L^*(\La)\otimes L(\La) \to  \FK{G}$ induced by $\phi\otimes v\mapsto \ti{f}_{\phi v}$, 
is an isomorphism of $G\times G$-modules.\vspace*{0.5ex}\\
Restricting the functions of $\FK{G}$ to $G'$ resp. $T_{rest}$ we get the algebras $\FK{G'}$ resp. $\FK{T_{rest}}$, the first
identical with the algebra of strongly regular functions as defined in \cite{KP2}, the second the classical 
coordinate ring of the torus $T_{rest}$. $\FK{G}$ is isomorphic to $\FK{G'}\otimes \FK{T_{rest}}$ by the comorphism of the 
multiplication map $G'\times T_{rest}\to G$.\vspace*{1ex}\\   
%
%
%
%
%
%
%
{\bf Substructures:}
%
%
%
For $\emptyset\neq J\subseteq I$ the submatrix $A_J:=(a_{ij})_{i,j\in J}$ of $A$ is a generalized Cartan matrix. There exist saturated sublattices 
$H(A_J)\subseteq H$, $P(A_J)\subseteq P$, with $(h_j)_{j\in J}\subseteq H(A_J)$, $(\al_j)_{j\in J}\subseteq P(A_J)$, giving an optimal realization 
for $A_J$.\\
We have $P=P(A_J)\oplus H(A_J)^\bot$, and the projections of $\La_j$, $j\in J$, onto $P(A_J)$ give a system of fundamental dominant weights of this 
optimal realization.\\
The corresponding Kac-Moody algebra $\g(A_J)$ embeds in $\g$. The root lattice $Q(A_J)$ identifies with the sublattice $Q_J:=\sum_{j\in J}\Z\al_j$ 
of $Q$. The Weyl group $\We(A_J)$ identifies with the parabolic subgroup $\We_J$. The set of roots $\W(A_J)$ identifies with $\W_J:=\W\cap Q_J$, 
and the set of real roots $\W(A_J)_{re}$ with $(\W_J)_{re}:= \rW\cap Q_J=\We_J\{\al_j\mid j\in J\}$. (To simplify the notation we sometimes identify 
the set of roots $\{\al_j\mid j\in J\}$ with $J$.)
The Kac-Moody group $G(A_J)$ embeds in $G$ in the obvious way.\vspace*{1ex}\\
The images of these embedding depend on the choice of the sublattice $H(A_J)$, only $H_J:=\Z\mb{-}span\Mklz{h_j}{j\in J}$ is 
uniquely determined by $A_J$.\\
The images of $\g(A_J)'$, $G(A_J)'$ are independent of this choice, and denoted by $\g_J$, $G_J$. It would be more consequent to write 
$\g_J'$ and $G_J'$ instead of $\g_J$ and $G_J$, but to simplify the notation we omit the prime.\\
The coordinate ring $\FK{G(A_J)'}$ identifies with the restrictions of $\FK{G}$ on $G_J$. (A similar statement for 
$\FK{G(A_J)}$ is not valid.)\vspace*{1ex}\\
Note that we have $\g_J=\bigoplus_{\al\in\Delta_J}\g_\al\oplus \h_J$, where $\h_J$ is spanned by the elements $h_j$, $j\in J$. 
$(G_J\,,\,(U_\al)_{\al\in(\Delta_J)_{re}}\,,\, T_J)$ is a root groups data system, where $T_J$ is the subtorus of $T$ generated by the elements 
$t_j(s)$, $j\in J$, $s\in\F^\times$.\vspace*{0.5ex}\\
We denote by $\n_J^\pm$ the subalgebras of $\g_J$ corresponding to $\n(A_J)^\pm$, and by $U_J^\pm$ the subgroups of
$G_J$ corresponding to $U(A_J)^\pm $, and set $B_J^\pm := T_J\ltimes U_J^\pm$.\\
We set $N_J:= G_J\cap N$. This group coincides with the group generated by the elements $n_j(1)$, $j\in J$, and by $T_J$.\vspace*{0.5ex}\\ 
To simplify the notation of many formulas, it is useful to set $\g_\emptyset:=\n_\emptyset:=\{0\}$, 
$G_\emptyset:= T_\emptyset:=U_\emptyset:=N_\emptyset:=\{1\}$, and $\FK{G_\emptyset}:=\F\,1$.\vspace*{1ex}\\
%
%
%
%
%
%
{\bf The case of a decomposable generalized Cartan matrix:}
%
%
%
If the generalized Cartan matrix is decomposable, $A=A_{I_1}\oplus A_{I_2}$, we can decompose the optimal realization and $H_{rest}$. Set
\begin{eqnarray*}
 H_{Rest\,1} &:=& \Mklz{x\in\HR}{\alpha_i (x)=0\;\mb{ for all }\;i\in I_2}\;\;,\\
 H_{Rest\,2} &:=& \Mklz{x\in\HR}{\alpha_i (x)=0\;\mb{ for all }\;i\in I_1}\;\;.
\end{eqnarray*}
We get optimal realizations of $A_1$, $A_2$ by:
\begin{eqnarray*}\begin{array}{ccccc}
  H_1 &:=& H_{I_1}\oplus H_{rest\,1} &\;,\; & \Pi^\vee_1\;:=\;\Mklz{h_i}{i\in I_1}\;,\\
  P_1 &:=& \Mklz{\la \in P}{\la(h)=0 \;\mb{ for }\;h\in H_2} &\;,\; & \Pi_1\;:=\;\Mklz{\alpha_i}{i\in I_1}\;.
\end{array}\\
\begin{array}{ccccc}
  H_2 &:=& H_{I_2}\oplus H_{rest\,2} &\;,\; & \Pi^\vee_2\;:=\;\Mklz{h_i}{i\in I_2}\;,\\
  P_2 &:=& \Mklz{\la \in P}{\la(h)=0 \;\mb{ for }\;h\in H_1} &\;,\; & \Pi_2\;:=\;\Mklz{\alpha_i}{i\in I_2}\;.
\end{array}\end{eqnarray*}
These satisfy:
\begin{eqnarray*}\begin{array}{ccc}
H\;=\;H_1\oplus H_2 &\quad,\quad & \Pi^\vee\;=\;\Pi^\vee_1\,\dcup\, \Pi^\vee_2 \;\;.\\
P\;=\;P_1\oplus P_2 &\quad,\quad & \Pi\;=\;\Pi_1\,\dcup\,\Pi_2 \;\;.  
\end{array}\end{eqnarray*}
The Weyl group $\We(A_i)$ can be identified with the parabolic subgroup $\We_{I_i}$ of $\We$, $i=1,2$, and we have $\We=\We_{I_1}\times\We_{I_2}$. 
If $X_{A_i}$ denotes the Tits cone of the optimal realization of $A_i$, $i=1,2$, then $X=X_{A_1}\oplus X_{A_2}$.\\
The Lie algebra $\g_i := \bigoplus_{\al\,\in\,\Delta_{I_i}}\g_\al\,\oplus\,\spann{H_i}$, 
the subgroup $G_i$ generated by $U_\al$, $\al\in (\Delta_{I_i})_{re}$ and by the elements $t_h$, $h\in H_i$, 
and the coordinate ring $\FK{G_i}$ obtained by restricting the functions of $\FK{G}$ onto $G_i$
are in the obvious way isomorphic to  $g(A_{I_i})$, $G(A_{I_i})$, and $\FK{G(A_{I_i})}$, $i=1,2$. 
Furthermore we have $\g =\g_1\oplus\g_2$, $G=G_1\times G_2$, and $\FK{G} = \FK{G_1}\otimes \FK{G_2}$.\vspace*{0.5ex}\\
%
%
%
%
%
\subsection{A generalization of affine toric varieties}                                                                                         
%
%
%
An affine toric variety is a normal affine variety containing a torus $T$ as a dense open subset, together with an action of $T$ on the variety, 
that extends the natural action of $T$ on itself. Alternatively an affine toric variety can be described as a normal algebraic monoid, whose unit 
group is a torus $T$. It is determined by a rational convex polyhedral cone, and can be constructed, starting with this cone, \cite{Fu}, \cite{Ne}.\\
In this subsection we state a generalization of this construction, starting with a not necessarily finitely generated rational convex cone. 
The results are similar as in the classical case, but some of the classical proofs can not be generalized. 
The proofs of these results, as well as some other results, can be found in \cite{M}. Here we omit the proofs, because certainly they can be also 
found elsewhere in the literature.
%
%
%
\subsubsection*{Rational convex cones and saturated submonoids of a lattice}
%
%
%
Let $L$ be a lattice of finite rank. Define the real vector space $V:=L\otimes_{\mathbb Z}\R\,$. Identify 
$L$ with $L\otimes 1$.\\
For a subset $Y\subseteq V$ denote by $\ccp{Y}$ / $\cc{Y}$ / $\spann{Y}$ the pointed convex cone / convex cone / subspace
generated by $Y$. The property ``finitely generated'' is defined in the obvious way.
\begin{Def}\mb{}\label{Fa1} A pointed convex cone / convex cone / subspace $X$ of $V$ is called rational, if there exists a subset of $L$, which
generates $X$.\\
A subsemigroup / submonoid / subgroup $S$ of $L$ is called saturated, if for all $x\in L$, and $n\in\N$ we have: $\quad nx\in S 
\quad\Rightarrow\quad x\in S$  
\end{Def}
The next theorem generalizes the Lemma of Gordan and its converse: 
\begin{Theorem}\mb{}\label{Fa2} The following maps are inverse, and respect the property ``finitely generated'':
\begin{eqnarray*}
  \{ \,\mb{rational pointed convex cones of }V\,\} && \{\,\mb{saturated subsemigroups of }L\,\} \\
   X \qquad \qquad\qquad&\mapsto & \qquad\qquad\quad X\cap L \\
   \ccp{S} \quad\qquad\qquad&\mapsot& \qquad\qquad\qquad S 
\end{eqnarray*}
Similar things hold if ``rational pointed convex cones'', ``saturated subsemigroups'' and ``$\ccp{S}$'' are replaced by\\
\hspace*{3em} ``rational convex cones'', ``saturated submonoids'' and ``$\cc{S}$'' ,\\
or if they are replaced by \\
\hspace*{3em}``rational subspaces'', ``saturated subgroups'' and ``$\spann{S}$''.\vspace*{1ex}\\
We call $X$ and $S$ associated, if they correspond under these maps.
\end{Theorem} 
Let $X\subseteq V$ be a convex cone. Recall the definitions of faces of $X$, the face lattice $Fa(X)$, exposed faces of $X$, the
 relative interior, and the hull of faces, \cite{Ro}. 
We denote a convex cone of the form $X-F$, where $F$ is a face of $X$, a dual face of $X$.\vspace*{1ex}\\
The algebraic forms of these definitions can be imitated for submonoids of a lattice:\\
Let $M$ be a submonoid of $L$. A subset $F\subseteq M$ is called a face of $M$, if $F$ is a submonoid, and $M\subseteq F$
a semigroup ideal, (i.e., $M\setminus F =\emptyset$, or else $M\setminus F\neq\emptyset$ and 
$M\setminus F +M\subseteq M\setminus F$).\\
The set of faces of $M$, denoted by $Fa(M)$, together with the inclusion of faces, is a lattice. The intersection of faces
coincides with the lattice intersection.\\
If $F$ is a face of $M$, then the set of faces of $F$ is given by
\begin{eqnarray*}
  \Fa{F} &=& \Mklz{G\in\Fa{M}\,}{\, G\subseteq F}\;.
\end{eqnarray*}
Define the relative interior of $F$ as 
\begin{eqnarray*}
  ri F &:=& F\setminus\bigcup_{G\in \Fa{F},\,G\neq F} G\;.
\end{eqnarray*}
The set $\Mklz{ri F\,}{\,F\in\Fa{M}}$ is a partition of $M$. $F$ is the smallest face containing an element $m\in ri F$.\\
The hull of a face $F$ is the group generated by $F$, which is equal to $F-F$. Note that $F=M\cap (F-F)$.\\
We denote a monoid $M-F$, where $F$ is a face of $M$, a dual face of $M$. Its set of faces is given by
\begin{eqnarray*}
  \Fa{M-F} &=& \Mklz{ G-F\,}{\, G\in \Fa{M} \mb{ with } G\supseteq F}\;.
\end{eqnarray*} 
\begin{Theorem}\mb{}\label{Fa3}\\
a) Let $X\subseteq V$ be a rational convex cone. Every face, relative interior of a face, hull of a face, and dual face of $X$ is 
rational.\vspace*{0.5ex}\\
b) Let $M\subseteq H$ be a saturated submonoid. Every face, relative interior of a face, hull of a face, and dual face of $M$ is 
saturated.\vspace*{0.5ex}\\
c) Let $X$ and $M$ be associated. By the maps corresponding to cones and submonoids of the last proposition, the 
face lattices $(\Fa{X},\subseteq)$ and $(\Fa{M},\subseteq)$ are isomorphic. The relative interiors, hulls, and 
dual faces of associated faces are associated.
\end{Theorem}
Our main example will be the Tits cone $X$, which is a P-rational convex cone in $\h_{\R}^*$. Its faces will be given explicitely. Because 
of their facetial structures, it is also not difficult to verify the last theorem directly.
%
%
%
\subsubsection*{Generalized affine toric varieties}
%
%
%
Let $M$ be a saturated submonoid of a lattice of finite rank. Let $\K$ be a field, $|\K|=\infty$. The set of homomorphisms of monoids 
\begin{eqnarray*}
 \ti{M} &:=& \mb{Hom}\left(\,(M,+\,)\,,\,(\K\,,\cdot\,)\,\right) 
\end{eqnarray*} 
gets the structure of an abelian monoid, by multiplying the homomorphisms pointwise. The monoid algebra $\KK{M}$ is identified with 
the coordinate ring $\KK{\ti{M}}$, identifying $\sum_{m\in M} c_m m$ with the function 
\begin{eqnarray*}
   (\sum_{m\in M} c_m\, m ) (\al) &:=& \sum_{m\in M} c_m\, \al(m) \quad,\quad \al\in\ti{M}\;\;.
\end{eqnarray*}  
The monoid structure of $\ti{M}$ induces a coalgebra structure on $\KK{\ti{M}}$.
$\ti{M}$ maps bijectively to the $\K$-valued points of $\KK{\ti{M}}=\KK{M}$, the inverse map given by restricting the $\K$-valued points onto 
$M\subseteq \KK{M}$.\\
We equip $\ti{M}$ with the Zariski topology induced by its coordinate ring.\vspace*{1ex}\\
If $M$ is a subgroup, then $\ti{M}= \mb{Hom}(M,\K)=\mb{Hom}(M,\K^\times)$ is a torus. 
\begin{Prop}\mb{}\label{Fa4} Let $M$ be a saturated submonoid of a lattice of finite rank. For $F\in\Fa{M}$ set
\begin{eqnarray*} 
     T(F)  &:=&  \Mklz{\al\in\ti{M}\,}{\,\al^{-1}(\K^\times)=F }\;\;.
\end{eqnarray*} 
Let $e(F)\in T(F)$ be the element given by
\begin{eqnarray*}
 e(F) m \;:=\; \left\{\begin{array}{cl}
                  1 & \mb{ if }\;m\in F \\ 
                  0 & \mb{ if }\;m\in M\setminus F
                  \end{array}  \right.\quad .
\end{eqnarray*}
1) $T(M)$ is the unit group of $\ti{M}$, $E:=\Mklz{e(F)}{F\in\Fa{M}}$ is the set of idempotents of
$\ti{M}$, and we have $\ti{M}= T(M)\, E = E\, T(M) $.\vspace*{1ex}\\ 
2) $T(F)$ is a subgroup of $\ti{M}$ with unit $e(F)$, isomorphic to the torus $\widetilde{F-F}$, an isomorphism
$\Phi_F:\,\widetilde{F-F}\to T(F)$ given by
\begin{eqnarray*}
    \Phi_F(\al)(m)\;:=\;\left\{ \begin{array}{cl}
                                  \al(m) & m\in F \\
                                     0   & m\in M\setminus F
                                 \end{array}\right. &\;\;,\;\; & \al\in \widetilde{F-F}\;\;.       
\end{eqnarray*}
The partition of $\ti{M}$ into $T(M)$-orbits is given by the tori $T(F)$, $F\in\Fa{M}$. \vspace*{0.7ex}\\
For $m\in M$ denote by $D(m)$ the principal open set $\{\al\in\ti{M}\,|\,\al(m)\neq 0\}$. We have
\begin{eqnarray*}
  \overline{T(F)} &=&   \bigcup_{G\in\Fa{F}} T(G)\;\;, \\
         D(m)    &=&   \bigcup_{G\in \Fa{M}\,,\,G\supseteq F } T(G)  \quad \mb{ where }\quad  m\in ri\,F\;\;.
\end{eqnarray*}
In particular, $T(M)$ is principal open and dense in $\ti{M}$.
\end{Prop}
{\bf Remark:} To use later, note that we get a surjective homomorphism of groups $\Psi_F:\,\widetilde{M-M}\to T(F)$ by
\begin{eqnarray*}
   \Psi_F(\al)(m)\;:=\;\left\{ \begin{array}{cl}
                                 \al(m) & m\in F \\
                                    0   & m\in M\setminus F
                                \end{array}\right. &\;\;,\;\; & \al\in \widetilde{M-M}         \;\;.
\end{eqnarray*}
%
%
%
\newpage\section{The monoid $\hat{G}$ and its structure}                                                                                 
%
%
%
To define the monoid $\GD$ we need the faces of the Tits cone. In the following subsection we state the description of these faces, 
as well as some other results, which are used later to investigate the structure of $\GD$.
%
%
%
\subsection{The face lattice of the Tits cone\label{FacesoftheTitscone}}
%
%
%
Induced by the action on the Tits cone $X$, the Weyl group $\We$ acts on the face lattice of $X$, which we denote by $(\,\RkX,\subseteq\,)$.
In this subsection we describe $(\,\RkX,\subseteq\,)$ together with its $\We$-action.\vspace*{1ex}\\
For $\emptyset\neq\Th \subseteq I$ we denote by $\Th^0$, resp. $\Th^\infty$ the set of indices 
corresponding to the sum of the components of $A_\Th$ of finite, resp. nonfinite type. We set
$\emptyset^0:=\emptyset^\infty:=\emptyset$.
\begin{Def}\mb{}\label{GD1} 
A subset $\Th\subseteq I$ is called special if $\Th=\Th^\infty$.
\end{Def}
Part 1 a), c), and 2) of the following theorem are due to E. 
Looijenga, \cite{Loo}, Lemma 2.2, and Corollary 2.3. Part 1 b) has been given by P. Slodowy in \cite{Sl1}, Kapitel 6.2, Korollar 2. An alternative 
proof can be found in \cite{M}.\vspace*{1ex}\\
For $J\subseteq I$ set $J^\bot:=\Mklz{i\in I}{a_{ij}=0\;\mb{ for all }j\in J}$. Recall from the classification of the generalized Cartan matrices, 
that a set $\emptyset\neq \Th\subseteq I$ is special if and only if $\left(\sum_{i\in\Th}\N\,h_i\right)\cap (-\Cq^\vee)\neq \emptyset$.

\begin{Theorem}\mb{}\label{GD2}
Let $\Th$ be special and set $R(\Th):=\We_{\Th^\bot}\FTq$.\vspace*{0.5ex}\\
1a) $R(\Th)$ is an exposed face of $X$. If $\Th\neq\emptyset$, then for every element $c\in\left(\sum_{i\in\Th}\N\,h_i\right)\cap (-\Cq^\vee)$ we have
\begin{eqnarray*}
  R(\Th) &  =        & X\cap \Mklz{\la\in \h_\R^*}{\la(c)=0}\;\;,\\
   X     & \subseteq & \Mklz{\la\in \h_\R^*}{\la(c)\geq 0}\;\;.
\end{eqnarray*}
b) The relative interior of $R(\Th)$ is given by 
\begin{eqnarray*}
 ri \,R(\Th) &=& \We_{\Th^\bot} \bigcup_{\Th^f\subseteq \Th^{\bot}\atop \Th^f=(\Th^f)^0} F_{\Th\cup\Th^f}\;\;.
\end{eqnarray*}
c) We have $\;\spann{R(\Th)} \:=\: \Mklz{\la\in \h_\R^*}{\la(h_i)=0\,\mb{ for all }\,i\in\Th}$.\vspace*{1ex}\\
2) The centralizers and normalizers of the $\We$-action:
\begin{eqnarray*}
    Z_{\cal W}(R(\Th)) &=&\Mklz{\sigma\in\We}{\,\sigma\la=\la\mb{ for all }\la\in R(\Th)}\;\,=\,\;\We_{\Th}\;\;, \\
      N_{\cal W}(R(\Th))&=&\Mklz{\sigma\in\We}{\sigma R(\Th)\;\,=\;\,R(\Th)}\;=\;\We_{\Th\cup \Th^\bot}\;\;.
\end{eqnarray*}
\end{Theorem}
The next corollary is an easy conclusion of part b) of the last theorem. Its first part has been given by P. Slodowy in \cite{Sl1}, Kapitel 6.2, 
Korollar 3. 
\begin{Cor}\mb{}\label{GD4} Every face of $X$ is $\We$-conjugate to exactly one of the faces $R(\Th)$, $\Th$ special.\\
A face $R$ is of the form $R(\Th)$, $\Th$ special, if and only if $ri\,R\,\cap \,\Cq\neq \emptyset$.
\end{Cor}
Due to this corollary, every face $R$ can be parametrized in the form $R=\sigma R(\Th)$ with an uniquely determined special set $\Th$, 
which we call the type of $R$, and $\sigma\in\We$. The element $\sigma$ is uniquely determined, if we restrict to the minimal coset 
representatives $\We^{\Th\cup\Th^\bot}$ of $\We/\We_{\Th\cup\Th^\bot}$.\vspace*{1ex}\\ 
Inclusion of faces can be rewritten into conditions involving the parametrisations:
\begin{Prop}\mb{}\label{GD5}
Let $\sigma,\sigma'\in\We$, $\Th,\Th'$ be special. The following statements 
are equivalent:
\vspace*{0.5ex}\\
\hspace*{0.65em} i) $\,\quad \sigma'R(\Th')\:\subseteq\:\sigma R(\Th)$
\vspace*{0.5ex}\\
\hspace*{0.65em} ii) $\quad \Th'\:\supseteq\:\Th\quad$ and $\quad
\sigma^{-1}\sigma'\:\in\:\We_{\Th^\bot}\We_{\Th'} $
\end{Prop}
\Proof
`$i)\Rightarrow ii)$': Because $\sigma R(\Th)$ contains facets of type $\Th'$, we conclude 
$\Th'\supseteq\Th$. Comparing the facets of type $\Th'$, we find $\sigma'F_{\Th'}\subseteq\sigma\We_{\Th^\bot}F_{\Th'}$,
which is equivalent to $\sigma^{-1}\sigma'\in\We_{\Th^\bot}\We_{\Th'}$.\vspace*{0.5ex}\\
`$ii)\Rightarrow i)$':
There exist $\,\tau_1\in\We_{\Th^\bot}$, $\tau_2\in\We_{\Th'}$, such that $\sigma^{-1}\sigma'=\tau_1\tau_2$. Because of $\Th'\supseteq\Th$, we have 
$\We_{{\Th'}^\bot}\subseteq\We_{\Th^\bot}$ and $\overline{F}_{\Th'}\subseteq\FTq $. Therefore
\begin{eqnarray*}
  \sigma^{-1}\sigma'R(\Th')\;=\;\tau_1\We_{\Th'^\bot}\overline{F}_{\Th'}\;\subseteq\;\We_{\Th^\bot}\FTq \;=\;R(\Th)\;\;.
\end{eqnarray*}
\Ende
It turns out to be difficult to rewrite the intersection of faces into conditions involving the 
parametrisations. But there is the following easy case:
If $\Th$, $\Th'$ are special, then also $\Th\cup\Th'$ is special, and we have $R(\Th)\cap R(\Th')
= R(\Th\cup\Th')$.\vspace*{1ex}\\
{\bf Some examples}:\\
{\bf 1)} If $A$ is of finite type the Tits cone is a linear space, and we have $\RkX=\{X\}$.\vspace*{0.5ex}\\
{\bf 2)} If $A$ is of affine or strongly hyperbolic type, then $\RkX=\{c,\,X\}$, where $c$ is the edge of the Tits cone.\vspace{0.5ex}\\
{\bf 3)} If $A$ is of indefinite type, but not strongly hyperbolic, then it is not difficult to see, that there exist infinitely many faces. 
As an example consider the hyperbolic generalized Cartan matrix  
\begin{eqnarray*}A&=&
\left(\begin{array}{ccc}  2 & -2 & 0 \\
                        -2 & 2 & -1 \\
                         0 & -1 & 2   
\end{array}\right)\;\;.\end{eqnarray*}
The form $(\;\mid\;)$ restricted to $\h_\R^*$ has signature (1,2). The Tits cone $X$ consists
of one component of the open cone $\Mklz{\la\in \h_\R^*}{(\la\mid\la)< 0}$ together with all 
closed, $P$-rational, isotropic half-lines on the boundary, compare \cite{FF}, \cite{Sl2} or \cite{M}.\\
The special sets are $\{1,2,3\}$, $\{1,2\}$, $\emptyset$. The corresponding faces of $X$ of these
types are the edge $c=\{0\}$, every closed, $P$-rational, isotropic half-line on the boundary, and 
$X$ itself.\\\\
The following two propositions are not difficult to prove, using the theorem, the corollary, and the proposition of above:
\begin{Prop}\mb{}\label{GD6}
Let $\emptyset\neq J\subseteq I$. Identify the Weyl group $\We(A_J)$ with the parabolic subgroup 
$\We_J$ of $\We\,$. The map
\begin{eqnarray*}
  \Upsilon_J:\;\; {\cal R}(X_{A_J})\, & \to     & \RkX      \\
   \tau R_{A_J}(\Th) & \mapsto & \tau R(\Th)
\end{eqnarray*}
is an embedding of lattices, and has the following image:
\begin{eqnarray*}
  \RkX_J &:=&  \Mklz{R\in\RkX}{R\supseteq R(J^\infty)} \\
  &\,=& \Mklz{\tau R(\Th)}{\Th\subseteq J \mb{ special },\;\tau\in\We_J} \\
  &\,=& \Mklz{\tau R(\Th)}{\Th\subseteq J^\infty \mb{ special },\;\tau\in
  \We_{J^\infty}}
\end{eqnarray*}
\end{Prop}
We set $\RkX_\emptyset:=\{X\}$.
\begin{Prop}\mb{}\label{GD7}
Let $A$ be decomposable, $A=A_1\oplus A_2$. With the notions of Subsection \ref{section2}, the embeddings 
$\Upsilon_{I_1}\,$, $\Upsilon_{I_2}$ are given by
\begin{eqnarray*}
\begin{array}{cccc}
  \Upsilon_{I_1}: & {\cal R}(X_{A_1}) & \to     & \RkX    \\
                  &       R    & \mapsto & R + X_{A_2}
\end{array} &\;\:,\:\; &
\begin{array}{cccc}
  \Upsilon_{I_2}: & {\cal R}(X_{A_2}) & \to     & \RkX    \\
                  &      R      & \mapsto & X_{A_1} + R
\end{array}\;,
\end{eqnarray*}
and $\RkX$ is the direct product of the sublattices $\RkX_{I_1}$ and $\RkX_{I_2}$. 
\end{Prop}
%
%
%
\subsection{The definition of the monoid $\hat{G}$}
%
%
%
The Tits cone $X$ is $P$-rational. The corresponding saturated submonoid $X\cap P$ of $P$ is related 
to the set of weights of the admissible highest weight modules $L(\La)$, $\La\in P^+$, by
\begin{eqnarray*}
  \bigcup_{\La\in P^+} P(L(\La)) &=& X\cap P\;\;.
\end{eqnarray*}
\begin{Def}\label{GD9}\mb{}
For a face $R$ of the Tits cone $X$ define a projection operator
\begin{eqnarray*}
 e(R) \;\in\; \Endo 
\end{eqnarray*}
by $\;
 e(R)v_\la\;:=\;\left\{\begin{array}{cl}
  v_\la & \la\in R\cap P \\
  0     & \la\notin R\cap P
 \end{array}\right.\quad,\quad v_\la\in L(\La)_\la,\;\la\in P(\La),\;\La\in P^+$.\vspace*{1ex}\\ Denote
its image by $B(R)$, and its kernel by $K(R)$. Set
\begin{eqnarray*}
  E \,:=\,\Mklz{e(R)\,}{\,R \,\mb{ a face of }\,X} &\;,\;&  E_{sp}\,:=\,\Mklz{e(R(\Th))\,}{\,\Th \,\mb{ special }}\;\;. 
\end{eqnarray*}
\end{Def}
{\bf Remark:} Fix $\La\in F_J\cap P$, $J\subseteq I$. It is not difficult to check that $e(R)$ does not act as zero on $L(\La)$ 
if and only if the type of $R$ is contained in $J$.\vspace*{1ex}\\ 
The following proposition can be proved easily: 
\begin{Prop}\mb{}\label{GD10}\\
1) $E$ is in the obvious way a monoid isomorphic to $(\,\RkX\,,\,\cap\,)$.\vspace*{0.5ex}\\
2) The monoid $N$ acts by conjugation on E. Explicitely, let $R$ be a face of the Tits cone $X$ and $n_\sigma\in N$. Then
\begin{eqnarray*}
  n_\sigma e(R) n_\sigma^{-1} &=& e(\sigma R) \;\;.
\end{eqnarray*}
Furthermore $n_\sigma B(R)=B(\sigma R)$, and $n_\sigma K(R)=K(\sigma R)$.
\end{Prop}
The Kac-Moody group $G$ acts faithfully on $\bigoplus_{\La\in P^+}L(\La)$. We identify $G$ with the corresponding subgroup of $\Endo$.
\begin{Def}\label{GD11} $\GD$ is the submonoid of $\Endo$ generated by $G$ and $E$.
\end{Def}
{\bf Remark:} The elements of $E$ are defined with respect to the Cartan subalgebra $\h$ of $\g$. 
The monoid $\GD$ does not depend on this choice: Because $G$ acts transitively on 
the Cartan subalgebras by the adjoint action, the set of projection operators defined with 
respect to $g\h$ is $gEg^{-1}$, ($g\in G$).\vspace*{1ex}\\
{\bf Some examples:} There are the following easy cases, which could be checked directly from the definitions:\\
{\bf 1)} If $A$ is of finite type, then $\GD=G$.\\
{\bf 2)} If $A$ is of affine type, then $\GD=(G'\cup \{e(c)\})\rtimes T_{rest}$. The element $e(c)$ is the zero of the monoid $G'\cup \{e(c)\}$.\\
{\bf 3)} If $A$ is of strongly hyperbolic type, then $\GD=G\cup \{e(c)\}$. The element $e(c)$ is the zero of the monoid $\GD$.\\
If $A$ is of indefinite type and not strongly hyperbolic, then due to the infinitely many faces of the Tits cone, the monoid $\GD$ is much more 
complicated.\\\\
Fix nondegenerate contravariant symmetric bilinear forms $\kBl$ on all modules $L(\La)$, $\La\in P^+$, and extend 
to a form on $\bigoplus_{\La\in P^+}L(\La)$, also denoted by $\kBl$, by requiring $L(\La)$ and $L(\La')$ to be orthogonal for 
$\La\neq\La'$.\\ 
It is easy to check, that the projections $e(R)$, $R\in\RkX$, are selfadjoint. Therefore taking the adjoint gives an involution $*$ of 
$\GD$, which extends the Chevalley involution of $G$. We also call this involution Chevalley involution. 
%
%
%
%
%
\subsection{Formulas for computations in $\hat{G}$}
%
%
%
The elements of $\GD$ are given by expressions of the form
\begin{eqnarray}\label{expression}
  g_1 e(R_1)\cdots g_m e(R_m) g_{m+1} &\mb{where}& \begin{array}{c}
                                             g_1,\ldots g_{m+1}\in G\\
                                             R_1,\ldots R_m\in \RkX
                                           \end{array},\; m\in \N\;.
\end{eqnarray}
It would be painful to investigate the structure of $\GD$ working with expressions of long length
$m$. The first thing to do is to look for some formulas, which allow to reduce the length of such
an expression (m=1 can be reached), and which allow to decide if two such expressions give the same 
element. These are given in Theorem \ref{GD18}, which is in essential of the following form:\vspace*{1ex}\\
Let $g,g'\in G$ and $R,R'\in \RkX$. Then:
\begin{eqnarray*}
  ge(R)\;=\;e(R')g'  &\iff & \mb{conditions on } g,g',R,R' \;\;.
\end{eqnarray*}
To prove this theorem we first determine some normalizers and centralizers, which are also useful at 
other places. To see why they are important look first at Proposition \ref{GD16}.
\begin{Def}\mb{}\label{GD12}
Let $V$ be a subspace of $\bigoplus_{\La\in P^+}L(\La)$, and let $Y$ be a submonoid of $\hat{G}$. Set
\begin{eqnarray*}
  \No{Y}{V} &:=& \Mklz{ x\in Y \,}{\, xV\:=\:V }\;\;, \\
  \Mb{Y}{V} &:=& \Mklz{ x\in Y \,}{\, xV\:\subseteq V }\;\;,\\
  \Mp{Y}{V} &:=& \Mklz{ x\in Y \,}{\, xV\:\supseteq V }\;\;,\\ 
  \Ze{Y}{V} &:=& \Mklz{ x\in Y \,}{\, \forall v\in V\;:\quad xv\:=\:v }\;\;. 
\end{eqnarray*}
\end{Def}
{\bf Remarks:} Let $Y$ be a subgroup of $G$.\vspace*{0.5ex}\\ 
1) If $\Mb{Y}{V}$ or $\Mp{Y}{V}$ is a group, then $\Mb{Y}{V}=\Mp{Y}{V}=\No{Y}{V}$.\vspace*{0.5ex}\\
2) Denote by $V^\bot$ the orthogonal complement of $V$. If $V^{\bot\bot}=V$, then we have 
\begin{eqnarray*}
  \quad N_Y(V^\bot)\;=\;N_{Y^*}(V)^*\;\;,\;\; M^\subset_Y (V^\bot)\;=\;M^\subset_{Y^*}(V)^*\;\;,\;\; 
  M^\supset_Y(V^\bot)\;=\;M^\supset_{Y^*}(V)^*\;\;.
\end{eqnarray*}
Because of $B(R)^\bot=K(R)$ and $K(R)^\bot=B(R)$, we can take $V=B(R),\,K(R)$.\\
\begin{Theorem}\mb{}\label{GD14}
Let $R=\tau R(\Th)$ be a face of the Tits cone, $\al\in\rW$. We have
\begin{eqnarray*}
\No{U_\al}{B(R)} &=& \Mb{U_\al}{B(R)}\;\,=\;\,\Mp{U_\al}{B(R)} \\
 &=& \left\{\begin{array}{ccl}
        U_\al & \mb{if} & \al\:\in \:\tau\left(\prW\cup\We_{\Theta}\Th\cup
        \We_{\Theta^\bot}\Th^\bot\right)  \\ 
        \{1\} & \mb{else} &
\end{array}\right.\;\;,\\[1ex]
 \Ze{U_\al}{B(R)} &=& \left\{\begin{array} {ccl}
    U_\al  &  \mb{if}  & \al\in \tau\Bigl(\,\left(\prW\setminus 
    \We_{\Theta^\bot}\Th^\bot\right)\cup\We_{\Theta}\Th\,\Bigr) \\
    \{1\}  & \mb{else} &
 \end{array}\right.\;\;.
\end{eqnarray*}
\end{Theorem}
\Proof
If $x_\al\in{\bf g}_\al$ and $v_\mu \in L(\La)_\mu\cap B(R)$, then
\begin{eqnarray*}
  \left(\exp x_\al\right) v_\mu &=& \sum_{k\in {\mathbb N}_0,\;\mu+k\al\,
  \in\, P(\Lambda)} \frac{x_\al^k}{k!}\,v_\mu \;\;. 
\end{eqnarray*}
To prove the theorem we have to investigate the half-strings $(\mu+\Nn\al)\cap P(\La)$ appearing in such sums.\vspace*{1ex}\\
$\bullet$ At first we examine the strings $(\mu+\Z\al)\cap P(\La)$ with $\mu\in P(\La)\cap R$:\vspace*{1ex}\\
If $\al\in\spann{R}\cap P$ then clearly:
\begin{eqnarray}\label{string1}
\al\mb{-string through }\mu\;\subseteq R\cap P\;\;.
\end{eqnarray}
If $\al\notin\spann{R}\cap P$ then:
\begin{eqnarray}
 \lefteqn{ \al\mb{-string through }\mu }
 \hspace*{3em}\nonumber\\
  &=& \left\{\begin{array}{lcc}
     \{ \mu,\,\mu-\al,\ldots\mu-\mu(h_\al)\al \} &\mb{ if } & \mu(h_\al) > 0 \\
                   \{ \mu \}                     &\mb{ if } & \mu(h_\al) = 0 \\
     \{ \mu,\,\mu+\al,\ldots\mu-\mu(h_\al)\al \} &\mb{ if } & \mu(h_\al) < 0 
   \end{array}\right.\;\;.\label{string2}
\end{eqnarray}
As an example we prove the case $\mu(h_\al)=0$, the other two cases are proved by similar arguments.
Suppose there exists an integer $z\in \Z\setminus\{0\}$, such that $\mu+z\al\in X$. Then also $\sigma_\al(\mu+z\al)\,=\,\mu-z\al\in X$. 
Because $R$ is a face of $X$, and 
\[ (\mu+z\al)+(\mu-z\al)\;\,=\;\,2\mu\;\,\in\;\, R\]
we get $\mu\pm z\al \in R$. Therefore
\[ (\mu+z\al)-(\mu-z\al)\;\,=\;\,2z\al\;\,\in\;\,\spann{R}\;, \]
contradicting $\al\notin\spann{R}$.\vspace*{1ex}\\
$\bullet$ Next we show:\vspace*{0.5ex}\\
1) The following statements are equivalent:\vspace*{0.5ex}\\
\hspace*{0.75em} i) $\;\;\forall\begin{array}{c}
\La\in P^+ \\
\mu\in R\cap P(\La)\end{array}\;:\;\;\;
(\mu+\Nn\al)\cap P(\La)\,\subseteq\, R\cap P$\vspace*{0.5ex}\\
\hspace*{0.75em} ii) $\;$Either \\
\hspace*{1.6em} $\;\quad \al\in \spann{R}\cap P\quad$.\\ 
\hspace*{1.9em} $\;$or\\ 
\hspace*{1.6em} $\;\quad \al\notin \spann{R}\cap P\quad $ and $\quad \forall\;\mu \in R\cap P\;:\;\; \mu(h_\al)\geq 0 \;\;$.\vspace*{0.7ex}\\
\hspace*{0.75em} iii) $\;\al\in \tau\left(\prW\cup\We_{\Theta}\Th\cup\We_{\Theta^\bot}\Th^\bot\right)\;\;$.\vspace*{1ex}\\
2) The following statements are equivalent:\vspace*{0.5ex}\\
\hspace*{0.75em} i) $\;\;\forall\begin{array}{c}
\La\in P^+ \\
\mu\in R\cap P(\La)\end{array}\;:\;\;\;(\mu+\Nn\al)\cap P(\La)\,=\, \{\mu\} \;\;$.\vspace*{0.5ex}\\
\hspace*{0.75em} ii) $\;\quad \al\notin \spann{R}\cap P\quad $ and $\quad \forall\; \mu\in R\cap P\;:\;\; \mu(h_\al)\geq 0 \;\;$.\vspace*{0.7ex}\\
\hspace*{0.75em} iii) $\;\al\in \tau\Bigl(\left(\prW\setminus\We_{\Theta^\bot}\Th^\bot\right)\cup\We_{\Theta}\Th\Bigr) \;\;$.\\\\
Taking into account $\bigcup_{\La\,\in\, P^+} P(\La)=X\cap P$, the equivalence of 1 i) and 1 ii) can be read off from (\ref{string1}) 
and (\ref{string2}).\vspace*{1ex}\\
If 2 ii) holds, then also 1 ii) holds, and therefore also 1 i). To show 2 i) suppose there exist $\La\in P^+$, $\mu\in R\cap P(\La)$, 
$n\in \N$, such that $\mu+n\al\in R\cap P$. Then we have $\al\in \spann{R}\cap P$, which contradicts 2 ii).\vspace*{1ex}\\
If 2 i) holds, then also  1 i) holds, and therefore also 1 ii). To show 2 ii) suppose $\al\in\spann{R}\cap P$. Then we have for all 
$\mu\in R\cap P$:
\begin{eqnarray*}
 \mu-\mu(h_\al)\al &=& \sigma_\al\mu\;\,\in\;\, \spann{R}\cap X \;\,=\;\,R \;\;,\\
 \sigma_\al \mu +\mu(h_\al)\al &=&\mu\;\,\in\;\, R \;\;.
\end{eqnarray*}
By using 2 i) and $\bigcup_{\La\,\in\, P^+} P(\La) = X\cap P$, we conclude $-\mu(h_\al)\leq 0$, $\mu(h_\al)\leq 0$, and therefore $\mu(h_\al)=0$. 
We have shown that $\sigma_\al$ fixes every element of $R$, and also every element of $\spann{R}$, which contradicts $\al\in\spann{R}$.\vspace*{1ex}\\
It is sufficient to show the equivalence of 1 ii), 1 iii) and the equivalence of 2 ii), 2 iii) only for $R=R(\Th)$.
For these equivalences it is sufficient to show the following statements:
\vspace*{0.5ex}\\
a) $\;\al\in\We_\Th \Th \quad\iff\quad \forall\; \mu\in R(\Th)\cap P\;:\;\; \mu(h_\al)=0 \;\;$.\vspace*{0.5ex}\\
b) $\;\al\in \We_{\Th^\bot}\Th^\bot \quad\iff\quad \al\in \spann{R(\Th)}\cap P \;\;$.\vspace*{0.5ex}\\
c) $\al\in\Delta_{re}^+\setminus\left(\We_\Th \Th \,\cup\, 
\We_{\Th^\bot} \Th^\bot\right) \quad\iff\quad $\\
\hspace*{1em} $\;\al\notin\spann{R(\Th)}\cap P\quad\mb{ and }\quad \forall\;\mu
\in R(\Th)\cap P\:: \:\;\mu(h_\al)\geq 0\quad\mb{ and }$\\
\hspace*{1em} $\;\exists\;\mu\in R(\Th)
\cap P \::\:\;\mu(h_\al)>0 \;\;$.\vspace*{1ex}\\
a) is valid due to $\We_\Th \,= \,Z_{\cal W}(R(\Th))\,$. To show b) note that
\begin{eqnarray}\label{NoWspannRTh}
  \We_\Th\,,\,\We_{\Th^\bot}\;\subseteq\;\We_{\Th\cup\Th^\bot}\;=\;N_{\cal W}(R(\Th))\;=\;
  N_{\cal W}(\spann{R(\Th)})\;,
\end{eqnarray}
where the last equality follows because for $\sigma\in N_{\cal W}(\spann{R(\Th)})$ we have
\begin{eqnarray*}
  \sigma R(\Th) &= & \left(\sigma R(\Th)-\sigma R(\Th)\right)\cap X \;\,=\,\;
  \sigma\left(R(\Th)-R(\Th)\right)\cap X \\
   &=&\left(R(\Th)-R(\Th)\right)\cap X \;\,=\;\, R(\Th) \;\;.
\end{eqnarray*}
By using $\Th^\bot\subseteq \spann{R(\Th)}$ and (\ref{NoWspannRTh}), we find $\We_{\Th^\bot}\Th^\bot\subseteq \spann{R(\Th)}\cap P$.\vspace*{0.5ex}\\
If $\al\in\spann{R(\Th)}\cap P$, then $\sigma_\al$ leaves $\spann{R(\Th)}$ invariant. Due to (\ref{NoWspannRTh}) we find 
$\al\in\We_{\Th\cup\Th^\bot}(\Th\cup\Th^\bot)=\We_\Theta \Th\dot{\cup}\We_{\Theta^\bot}\Th^\bot$.\\
Suppose there exist $\sigma\in\We_\Th$, $i\in \Th$, such that $\sigma\al_i\in\spann{R(\Th)}$. Because of (\ref{NoWspannRTh}) we find 
$\al_i \in \spann{R(\Th)}$, which contradicts $\al_i(h_i)=2$.\vspace*{1ex}\\
To c): We first show `$\Rightarrow$': Due to b) we find $\al\notin\spann{R(\Th)}\cap P$. Because of  a) there exists an element 
$\mu\in R(\Th)\cap P$ such that $\mu(h_\al)\neq 0$.\\
Note that due to \cite{KP3}, Lemma 2.1, we have 
\begin{eqnarray*}
\prW\setminus \We_{\Th\cup\Th^\bot}(\Th\cup\Th^\bot) &=&
\bigcap_{ \eta\,\in\,{\cal W}_{\Th\cup\Th^\bot} }\eta\prW\;\;.
\end{eqnarray*}
Therefore this is a $\We_{\Th\cup\Th^\bot}$-invariant set of positive roots. For every $\mu\in R(\Th)\cap P$ there exist 
$\sigma\in \We_{\Th\cup\Th^\bot}$, $\ti{\mu}\in \overline{F}_\Th \cap P$, such that $\mu=\sigma\ti{\mu}$ and we find: 
\[ \mu(h_\al)\;\,=\;\,\sigma\ti{\mu}(h_\al)\;\,=\;\,\ti{\mu}(h_{\sigma^{-1}   \al}) \;\, \geq\;\, 0\] 
Next we show `$\Leftarrow$': Because of a) and b) we find $\al\notin\We_\Th \Th \cup \We_{\Th^\bot} \Th^\bot$. 
Suppose $\al\in\Delta_{re}^-\setminus\left(\We_\Th \Th \cup\We_{\Th^\bot} \Th^\bot\right)$. Then $(-\al)\in 
\Delta_{re}^+\setminus\left(\We_\Th \Th \cup \We_{\Th^\bot} \Th^\bot\right)$, and due to c) `$\Rightarrow$' there
exists an element $\mu\in R(\Th)\cap P$, such that $\mu(h_{-\al})>0$. Inserting $h_{-\al}=-h_{\al}$ leads to a contradiction.\vspace*{1ex}\\
$\bullet$ Now we can proof the first statement of the theorem. Because of 
Remark 1) following Definition \ref{GD12} it is sufficient to show
\begin{eqnarray*}
  \Mb{U_\al}{B(R)} &=& \left\{\begin{array}{ccl}
        U_\al & \mb{if} & \al\:\in \:\tau\left(\prW\cup\We_{\Theta}\Th\cup
        \We_{\Theta^\bot}\Th^\bot\right)  \\ 
        \{1\} & \mb{else} &
\end{array}\right.\;\;.
\end{eqnarray*}
Let $\al\in \tau\left(\prW\cup\We_{\Theta}\Th\cup\We_{\Theta^\bot}\Th^\bot\right) $. We have to show the inclusion `$\supseteq$'.
If $x_\al\in {\bf g_\al}$, $v_\mu\in L(\La)_\mu\subseteq B(R)$, then 
\begin{eqnarray*}
  \left(\exp x_\al\right) v_\mu \;\,=\sum_{k\in {\mathbb N}_0,\; \mu+k\al\,
  \in\, P(\Lambda)} \underbrace{ \frac{x_\al^k}{k!}\,v_\mu }_{ \in\,
  L(\Lambda)_{\mu+ k\al} }  \;,
\end{eqnarray*}
and due to the equivalence of 1 iii) and 1 i) we have $P(\La)\cap(\mu+\Nn\al)\subseteq R\cap P$. Therefore 
$(\exp x_\al)v_\mu \in B(R)$.\vspace*{1ex}\\
Let $\al\notin\tau\left(\prW\cup\We_{\Theta}\Th\cup\We_{\Theta^\bot}\Th^\bot\right) $. We have to show the inclusion `$\subseteq$'.
Because of the equivalence of 1) i) and 1 iii) there exist $\La\in P^+$, $\mu\in P(\La)\cap R$, such that 
$P(\La)\cap (\mu+\Nn\al)\not\subseteq P\cap R$.\\
Now ${\bf s_\al} = \g_\al\oplus[\g_\al,\g_{-\al}]\oplus\g_{-\al}$ is isomorphic to $sl(2,\F)$, and there exists a decomposition of $L(\La)$ 
in a direct sum of irreducible, finite dimensional ${\bf s_\al}$-modules, which are $\h$-invariant. Because the $\al$-string through $\mu$ 
is finite, there is a ${\bf s_\al}$-modul $V$ among these modules, such that the $\al$-string through $\mu$ restricted to $\F h_\al$ gives 
the set of weights of $V$.\\
Suppose there exists an element $x_\al\in \gal\setminus\{0\}$ such that $\exp (x_\al)\in \Mb{U_\al}{B(R)}$. Then for 
$v_\mu\in V\cap L(\La)_\mu\setminus\{0\}$ we have
\begin{eqnarray*}
  B(R) &\ni &
  \left(\exp x_\al\right) v_\mu \,\;= \sum_{k\in {\mathbb N}_0,\; \mu+k\al\,
  \in\, P(\Lambda)} \underbrace{ \frac{x_\al^k}{k!}\,v_\mu }_{ \in\,
  L(\La)_{\mu+k\al}\setminus \{0\} } \;,
\end{eqnarray*}
which is a contradiction.\vspace*{1ex}\\
$\bullet$ The second statement of the theorem is proved similar to the first. Use 
the equivalence of 2 i) and 2 iii) instead of the equivalence of 1 i) and 1 iii).\\
\Ende
\begin{Theorem}\mb{}\label{GD15}
Let $R=\tau R(\Th)$ be a face of the Tits cone $X$. We have:\vspace*{1.5ex}\\
\hspace*{0.65em} 1) a) $\;\Mb{T}{B(R)}\;=\;\Mp{T}{B(R)}\;=\;\No{T}{B(R)}
\;=\;\qquad T\;$.\vspace*{0.5ex}\\
\hspace*{1.9em} b) $\quad\Ze{T}{B(R)} \,\;=\,\; \qquad \tau\, T_\Theta\,
\tau^{-1} \;$.\vspace*{1.5ex}\\
\hspace*{0.65em} 2) a) $\;\Mb{N}{B(R)}\;=\;\Mp{N}{B(R)}\;=\;\No{N}{B(R)}\;=\;
  \tau\,\ti{N}_{\Theta\cup\Theta^\bot}\,\tau^{-1} \;$.\vspace*{0.5ex}\\
  \hspace*{1.9em} b) $\quad\Ze{N}{B(R)} \,\;=\,\; \qquad\tau\,N_\Theta\,
\tau^{-1} \;$.\vspace*{1.5ex}\\
\hspace*{0.65em} 3) a) $\;\Mb{G}{B(R)}\;=\;\Mp{G}{B(R)}\;=\;\No{G}{B(R)}\;=\;
  \tau\, P_{\Theta\cup\Theta^\bot}\,\tau^{-1} \;$.\vspace*{0.5ex}\\
  \hspace*{1.9em} b) $\quad\Ze{G}{B(R)} \,\;=\,\; \tau\,\left(\,G_\Theta\ltimes 
U^{\Theta\cup\Theta^\bot}\,\right)\,\tau^{-1} \;$.\\
\end{Theorem}
\Proof
Because of Remark 1) following Definition \ref{GD12} we have to show 1a), 2a), and 3a) only for $\Mb{Y}{B(R)}$. Because of
\begin{eqnarray*}
  \Mb{Y}{B(\tau R(\Th))} &=& \Mb{Y}{n_\tau B(R(\Th))}\,\;=\,\;n_\tau \Mb{Y}{B(R(\Th))}n_\tau^{-1}\;\;, \\
    \Ze{Y}{B(\tau R(\Th))} &=& \Ze{Y}{n_\tau B(R(\Th))} \;\,=\;\, n_\tau\,\Ze{Y}{B(R(\Th))}n_\tau^{-1} \;\;,
\end{eqnarray*}
we can restrict to the cases $R=R(\Th)$, $Y=T,\,N,\,G$.\vspace*{1ex}\\
Obviously {\bf 1a)} holds. {\bf To 1b):}
Write $t\in T$ in the form $t=\prod_{i=1}^{2n-l}t_i(s_i)$ where $s_i\in \F^\times$. Then $t\in \Ze{T}{B(R(\Th))}$ if and only if 
$\prod_{i}s_i^{\la(h_i)}=1$ for all $\la\in R(\Th)\cap P$. Because of  $R(\Th)=\{\,\la\in X \,|\, \la(h_k)=0\,,\, k\in \Th\,\}\ni \La_i $ for all 
$i\notin \Th$, this is equivalent to $s_i=1$ for all $i\notin \Th$. \vspace*{1ex}\\
{\bf To 2a):}
Let $n_\sigma\in N$. Comparing the decompositions into weight spaces, we conclude that $n_\sigma B(R(\Th))\subseteq B(R(\Th))$ is equivalent 
to $\sigma R(\Th) \subseteq R(\Th)$. Due to Proposition \ref{GD5} this is equivalent to $\sigma\in \We_{\Th\cup\Th^\bot}$.\vspace*{1ex}\\
{\bf To 2b):}
$N_\Th$ is generated by $T_\Th $ and the elements $n_i(1)=\exp(e_i)\exp(-f_i)\exp(e_i)$, $i\in\Th$. Due to 1b) and 
Theorem \ref{GD14} we get $\Ze{N}{B(R(\Th))}\supseteq N_\Th$.\vspace*{0.5ex}\\
To show the reverse inclusion write $n_\sigma\in\Ze{N}{B(R(\Th))}$ in the form
\begin{eqnarray*}
  n_\sigma\;=\; t n_{i_1}(1)\cdots n_{i_k}(1)  &\mb{ where } & t\in T\,,\;\; \sigma\,=\,\sigma_{i_1}\cdots\sigma_{i_k} \,\in\,\We\; \mb{ reduced}\;\;.
\end{eqnarray*}
Because $n_\sigma$ fixes the weight spaces of $B(R(\Th))$, we have $\sigma\in Z_{\cal W}(R(\Th))=\We_\Th$, from which we conclude 
$i_1,\,\ldots,\, i_k\in\Th$, compare \cite{Hu2}, Chapter 5.5.\\
Because $n_{i_1}(1),\,\ldots\,\,n_{i_k}(1),\,n_\sigma$ are elements of the group $\Ze{N}{B(R(\Th))}$, we get 
$t\in \Ze{N}{B(R(\Th))}\cap T = \Ze{T}{B(R(\Th))}=T_\Th$.\vspace*{1ex}\\
{\bf To 3a):}
$U^+$ is generated by $U_\al$, $\al\in \prW$. Using Theorem \ref{GD14} we get 
\begin{eqnarray}\label{Asterix}
 u \,B(R(\Th)) \;=\; B(R(\Th)) \;& \mb{for all} &\; u\in U^+\;\;.
\end{eqnarray}
Due to the Bruhat decomposition of $G$, an element $g\in G$ can be written in the form $g = un\ti{u}$ with $u,\ti{u}\in U^+$, $n\in N$. 
Using (\ref{Asterix}) and 2a) we find
\begin{eqnarray*}
 g \in \Mb{G}{B(R(\Th))} \;\iff\; n\underbrace{\left(\ti{u} B(R(\Th))\right)}_{=\,B(R(\Th))}\:\subseteq\:
 \underbrace{u^{-1} B(R(\Th))}_{=\, B(R(\Th))}\;\iff\;n\in \tilde{N}_{\Th\cup\Th^\bot}\;\;.
\end{eqnarray*}
Therefore $\Mb{G}{B(R(\Th))}=U^+\ti{N}_{\Th\cup\Th^\bot}U^+ = P_{\Th\cup\Th^\bot}$.\\\\
{\bf To 3b):}
$G_\Th$ is generated by $U_{\al}$, $\al\in\pm\Th$. By using Theorem \ref{GD14} we get
\begin{eqnarray*}
G_\Th\;\subseteq\; \Ze{G}{B(R(\Th))}\;\;.
\end{eqnarray*}
The group $U^{\Th\cup\Th^\bot}$ is the normal subgroup of $P_{\Th\cup\Th^\bot}$ generated by the root groups $U_\al$, 
$\al\in\prW\setminus(\We_\Th\Th\cup\We_{\Th^\bot}\Th^\bot)$.
Due to Theorem \ref{GD14} these root groups are subgroups of the group $\Ze{G}{B(R(\Th))}$, which 
is normal in $\No{G}{B(R(\Th))}= P_{\Th\cup\Th^\bot}$. Therefore we conclude
\begin{eqnarray*}
     U^{\Th\cup\Th^\bot}\;\,\subseteq\;\,\Ze{G}{B(R(\Th))}\;\;.
\end{eqnarray*}
Now let $g\in\Ze{G}{B(R(\Th))}$. Because of $\Ze{G}{B(R(\Th))}\subseteq\No{G}{B(R(\Th))}=P_{\Th\cup\Th^\bot}$, and because of the decompositions
\begin{eqnarray*}
   P_{\Th\cup\Th^\bot}  &=& \ti{G}_{\Th\cup\Th^\bot}\ltimes U^{\Th\cup\Th^\bot}\;\;,\\
   \ti{G}_{\Th\cup\Th^\bot}&=& (G_\Th\times G_{\Th^\bot})\rtimes (T_{I\setminus (\Th\cup\Th^\bot)}T_{rest})\;\;,
\end{eqnarray*}
there exist elements $x\in G_{\Th^\bot}T_{I\setminus (\Th\cup\Th^\bot)}T_{rest}$, $y\in G_\Th\ltimes U^{\Th\cup\Th^\bot}$ such that $g=xy$, and
furthermore $x=gy^{-1}\in \Ze{G}{B(R(\Th))}$.\vspace*{0.5ex}\\
We have to show $x=1$. Because of $U^+\cap U^-=\{1\}$, it is sufficient to show $x,x^*\in U_{\Th^\bot}$.\\
$\al)$ At first we show $x\in U^{\Th^\bot}$:
Because of the Bruhat decomposition of $G_{\Th^\bot}T_{I\setminus(\Th\cup\Th^\bot)}T_{rest}$, we can write $x$ as a product
\begin{eqnarray*}
  x\;\,=\;\,un_\sigma \ti{u} &\mb{ where } & u,\,\ti{u}\in U_{\Th^\bot}\,,\;\;n_\sigma\in N_{\Th^\bot} T_{I\setminus (\Th\cup\Th^\bot)}T_{rest} \;\;.
\end{eqnarray*}
By applying $x$ on $v_\La\in L(\La)_\La\setminus\{0\}$, $\La\in\FTq\cap P$, we get
\begin{eqnarray*}
   v_\La\,\;=\,\;u n_\sigma\ti{u}v_\La\;\,=\,\;\underbrace{n_\sigma v_\La}_{\neq 0}+\sum_{q\in Q^+}v_{\sigma\La+q}\;\;.
\end{eqnarray*}
Comparing the components of the weight spaces we find 
\begin{eqnarray}\label{nsigvLa}
  n_\sigma v_\La &=& v_\La\;\;.
\end{eqnarray}
Because we have $\sigma\La=\La$ for all $\La\in\FTq\cap P$, we conclude $\sigma\in\We_\Th\cap\We_{\Th^\bot}=\{1\}$. Write 
$n_\sigma = n_1\in T_{I\setminus (\Th\cup\Th^\bot)}T_{rest}T_{\Th^\bot}$ in the form
\begin{eqnarray*}
  n_\sigma\;=\;\prod_{i\notin\Th,\,i=1,\ldots 2n-l} t_i(s_i) &\;\,,\;\, & s_i\in\F^\times\;\;,
\end{eqnarray*} 
and insert in (\ref{nsigvLa}). Since $\La_k\in \FTq\cap P$ for all $k\in \{1,\,\ldots,\, 2n-l\}\setminus\Th$, we find 
\begin{eqnarray*}
    1\;\,=\,\;s_k &\mb{ for all } & k\,\in\,\{1,\,\ldots,\, 2n-l\}\setminus \Th\;\;.
\end{eqnarray*}
Therfore we get $n_\sigma=1$, and $x=u\ti{u}\in U_{\Th^\bot}$.\vspace*{0.5ex}\\
$\beta)$ To show $x^*\in U_{\Th^\bot}$ it is sufficient to show $x^*\in\Ze{G}{B(R(\Th))}$, because then we can apply $\al)$ to $x^* $.\\
Let $v\in B(R(\Th))$. Because of $x\in \Ze{G}{B(R(\Th))}$, we find 
\begin{eqnarray*}
 \kkB{x^*v-v}{v'}   \;\,=\,\; \kkB{v}{xv'-v'} \,\;=\,\; 0 &\mb{ for all }& v'\in B(R(\Th))\;\;.
\end{eqnarray*}
We get $x^*v-v\in (B(R(\Th)))^\bot = K(R(\Th)) $.\\ 
Because of 3 a) we have $x^*\in (U_{\Th^\bot})^*\subseteq \No{G}{B(R(\Th))}$, and therefore we also get $x^*v-v\in B(R(\Th))$.\\ 
Because $B(R(\Th))$ and $K(R(\Th))$ intersect trivially, we conclude $x^*v=v$.\\
\End\\
%
%
%
%
%
%
\begin{Prop}\mb{}\label{GD16}
Let $R$, $S$ be faces of the Tits cone. Let $Y$ ba a subgroup of $G$, and $y\in Y$. We have:
\begin{eqnarray}
  y e(R)\;=\;e(S) \quad &\iff &\quad R\:=\:S\;\;\mb{ and }\;\;y\in \Ze{Y}{B(R)} \;\;.\quad\label{Glzenl}\\
  e(R) y\;=\;e(S) \quad &\iff &\quad R\:=\:S\;\;\mb{ and }\;\;y\in \Ze{Y}{B(R)}^* \;\;.\quad\label{Glzenr}\\[1.5ex]
  y e(R)y^{-1}\;=\;e(R) \quad &\iff &\quad y\in \No{Y}{B(R)}\cap\No{Y^*}{B(R)}^* \;\;.\quad\label{Glnor}
\end{eqnarray}
\end{Prop}
{\bf Remarks:}\\ 
1) For $R\:=\:\tau R(\Th)$ we get due to Theorem \ref{GD15}$\,$:\vspace*{1.5ex}\\
$\begin{array}{c||c|c|c|}
Y & \biggl. T\biggr. & N & G \\ \hline \hline
\Ze{Y}{B(R)} & \biggl.\tau\, T_\Th\, \tau^{-1} \biggr.& \tau\, N_\Th\, \tau^{-1} &
\tau\,\left(G_\Th\ltimes U^{\Th\cup\Th^\bot}\right)\, \tau^{-1} \\ \hline 
\left.\No{Y}{B(R)}\cap \No{Y^*}{B(R)}^*\right. &\;\;T\;\; & 
\left.\tau\, \ti{N}_{\Th\cup
\Th^\bot}\,\tau^{-1} \right.& \biggl.\tau\,\ti{G}_{\Th\cup\Th^\bot}\, \tau^{-1}\biggr. 
\end{array}$\vspace*{1.5ex}\\
2) For $R=S=R(\Th)$, and $Y=G$ this has been claimed by D. Peterson.\vspace*{1ex}\\ 
\Proof
We first show (\ref{Glzenl}), then  (\ref{Glzenr}) follows by application of the 
involution $*\,$. Let $ye(R)=e(S)$. The kernels $K(R),K(S)$ of 
the two linear maps are the same. Comparing the decompositions of $K(R)$ and $K(S)$ 
into weight spaces, we find $P\,\cap\,(X\setminus R)=P\,\cap\,(X\setminus S)$. Therefore $R=S$.\\    
Now $ye(R)=e(R)$ is equivalent to $ye(R)v=e(R)v$ for all $v\in\bigoplus_{\La\in P^+}L(\La)$, which is in turn equivalent to 
$y\in\Ze{Y}{B(R)}$.\vspace*{0.5ex}\\
Next we show show (\ref{Glnor}). The linear projections $y e(R) y^{-1}$, $e(R)$ are equal if and only if their images $y B(R)$, $B(R)$, 
and their kernels $y K(R)$, $K(R)$ are equal. This is equivalent to $y\in \No{Y}{B(R)}\cap\No{Y}{K(R)}$, and by Remark 2) following 
Definition \ref{GD12}, we get $\No{Y}{K(R)} = \No{Y^*}{B(R)}^*$.\vspace*{1ex}\\
\End\\
Let $\Th$ be special. Because of the Levi decomposition $P_{\Th\cup\Th^\bot}=\ti{G}_{\Th\cup\Th^\bot}\ltimes U^{\Th\cup\Th^\bot}$,
and the decomposition $\ti{G}_{\Th\cup\Th^\bot}=(G_\Th\times G_{\Th^\bot})\rtimes T_\Th^{co}$, where 
$T_\Th^{co}:=T_{I\setminus (\Th\cup\Th^\bot)}T_{rest}$, we get the following decompositions:
\begin{eqnarray*}
 P_{\Th\cup\Th^\bot} \;\:&=& \left(G_{\Th^\bot}\rtimes T_\Th^{co}\right)\ltimes
 \left(G_\Th\ltimes U^{\Th\cup\Th^\bot}\right)\;\;,\\
 \left(P_{\Th\cup\Th^\bot}\right)^* &=& \left(G_\Th\ltimes U^{\Th\cup\Th^\bot}\right)^*
 \rtimes\left(G_{\Th^\bot}\rtimes T_\Th^{co}\right)\;\;.
\end{eqnarray*}
The projections belonging to these semidirect products are denoted by
\begin{eqnarray*}
 p_\Th\,:\; P_{\Th\cup\Th^\bot}  \;\;  & \to & G_{\Th^\bot}\rtimes T_\Th^{co} \;\;,\\
 p^*_\Th\,:\; \left(P_{\Th\cup\Th^\bot}\right)^* & \to & G_{\Th^\bot}\rtimes T_\Th^{co}\;\;.
\end{eqnarray*}
It is easy to see, that we have $p^*_\Th = * \circ p_\Th\circ * $, and  $p^*_\Th(x) = p_\Th(x)$ for all $x\in P_{\Th\cup\Th^\bot}\cap 
(P_{\Th\cup\Th^\bot})^* = \ti{G}_{\Th\cup\Th^\bot}$.\\
\begin{Prop}\mb{}\label{GD17}\\
If $x\in P_{\Th\cup\Th^\bot}$ then
$\qquad x e(R(\Th))\,=\, p_\Th(x) e(R(\Th)) \,=\, e(R(\Th))p_\Th(x)\;\;$.\vspace*{0.5ex}\\
If $x\in \left(P_{\Th\cup\Th^\bot}\right)^*$ then
$\quad e(R(\Th)) x\;=\; e(R(\Th)) p^*_\Th(x)\;=\; p^*_\Th(x)e(R(\Th))\;\;$.\vspace*{1ex}\\
In particular we have
\begin{eqnarray*}
  U^+ e(R(\Th)) &=& U^+_{\Th^\bot}e(R(\Th)) \,\;=\;\, e(R(\Th))U^+_{\Th^\bot}    \;\;,\\
  e(R(\Th)) U^- &=& e(R(\Th)) U^-_{\Th^\bot} \,\;=\;\, U^-_{\Th^\bot} e(R(\Th))  \;\;.
\end{eqnarray*}
\end{Prop}
{\bf Remark:}
D. Peterson claimed $U^+ e(R(\Th))\subseteq e(R(\Th)) U^+$, and $e(R(\Th))U^-\subseteq U^-e(R(\Th))$.\vspace*{1ex}\\
\Proof
Write $x\in P_{\Th\cup\Th^\bot}$ in the form 
\begin{eqnarray*}
 x\;=\; p_\Th(x) a &\mb{ with }& a\in G_\Th\ltimes U^{\Th\cup\Th^\bot}\,,\;
 \,p_\Th(x)\,\in\, G_{\Th^\bot}T^{co}_\Th \;\;.
\end{eqnarray*}
By using (\ref{Glzenl}) and (\ref{Glnor}) of Proposition \ref{GD16} we get
\begin{eqnarray*}
  x e(R(\Th)) \;\,=\,\; p_\Th(x) \underbrace{a e(R(\Th))}_{=\:e(R(\Th))} \,\;=\,\; 
                                 \underbrace{p_\Th(x) e(R(\Th)) p_\Th(x)^{-1}}_{=\:e(R(\Th))} p_\Th(x)\;\;.
\end{eqnarray*}
Because of the decompositions $U^+ = U_{\Th\cup\Th^\bot}\ltimes U^{\Th\cup\Th^\bot}$ and $U_{\Th\cup\Th^\bot} = U_{\Th^\bot}\times U_\Th$, 
we have $p_\Th(U^+)=U_{\Th^\bot}$.\vspace*{1ex}\\
The remaining equations follow by application of the involution $*$.\\
\Ende\\
The last two propositions are special cases of the following theorem, which is the main theorem for computations in $\GD$. 
In the next subsection we will see, that this theorem allows to reduce any computation in $\GD$ to computations in $G$ and $\RkX$.
\begin{Theorem}\mb{}\label{GD18}\\
a) Let $x,y\in G$ and $\Th,\Xi$ special. The following statements are equivalent:\vspace*{1ex}\\
\hspace*{1.3em} i) $\;\;\quad x e(R(\Th))\;=\; e(R(\Xi))y \;$.\vspace*{0.5ex}\\
\hspace*{1.3em} ii) $\;\Th\:=\:\Xi\quad $ and
$\quad x\in P_{\Th\cup\Th^\bot}\,,\;y\in\left(P_{\Th\cup
  \Th^\bot}\right)^*\;$ with $\;p_\Th(x)\:=\:p_\Th^*(y)\;$.\vspace*{1ex}\\
b) Let $R$ be a face of the Tits cone, and $n_\sigma\in N$. Then
\begin{eqnarray*}
  n_\sigma e(R) n_\sigma^{-1}\;\,=\;\, e(\sigma R)\;\;.
\end{eqnarray*}  
\end{Theorem}
\Proof
b) has already been shown, and a), `$ii) \Rightarrow i)$' is an easy consequence of the last 
proposition. For the reverse direction we show at first $\Th=\Xi$:
Comparing the images of $x e(R(\Th))$ and $e(R(\Xi))y$ we get
\begin{eqnarray*}
  x B(R(\Th)) \;\,=\;\,B(R(\Xi))\;\;.
\end{eqnarray*}
Write $x$ in the form $x=un_\sigma \ti{u}$ with $ u,\ti{u}\in U^+$, $n_\sigma\in N$, and insert in this equation. 
Since $U^+\subseteq \No{G}{B(R(\Th))},\,\No{G}{B(R(\Xi))}$ we find $n_\sigma B(R(\Th)) =B(R(\Xi))$. 
Comparing the decompositions into weight spaces, 
we get $\sigma R(\Th)\cap P =R(\Xi)\cap P$. From this follows $\Th = \Xi$.\\
Now we can show the remaining statements of ii). Comparing the images and kernels of $x e(R(\Th))$ and $e(R(\Th))y$, we find
\begin{eqnarray*}
  x B(R(\Th))\;\,=\;\,B(R(\Th)) &\:\;,\;\: &  K(R(\Th))\;\,=\,\; y^{-1} K(R(\Th)) \;\;.
\end{eqnarray*}
Due to Theorem \ref{GD15} and Remark 2) following Definition \ref{GD12}, we get $x\in P_{\Th\cup\Th^\bot}$, $y\in (P_{\Th\cup\Th^\bot})^*$.
Using the last proposition and Proposition \ref{GD16} we find:
\begin{eqnarray*}
  p_\Th(x) e(R(\Th))\;=\; x e(R(\Th))\;=\; e(R(\Th))y \;=\; p_\Th^*(y) 
  e(R(\Th)) \qquad\\
  \Rightarrow\quad  \underbrace{\left(p_\Th(x)\right)^{-1}p_\Th^*(y)}_{\in\,
  T_\Th^{co} G_{\Th^\bot} }\:\in\: G_\Th\ltimes U^{\Th\cup\Th^\bot}
  \quad\Rightarrow\quad \left(p_\Th(x)\right)^{-1}p_\Th^*(y)\:=\:1\;\;.
\end{eqnarray*}
\End\\
The next corollary is an easy consequence of the last theorem:  
\begin{Cor}\mb{}\label{GD19}
Let $R,\,S$ be faces of the Tits cone $X$. Let $x,\,y\in G$. Then:
\begin{eqnarray*}
 x e(R)\;=\; e(S) y \quad &\Rightarrow &\quad \exists \;\tau\in\We\;:\;\;\;S\;=\;\tau R \;\;.
\end{eqnarray*}
\end{Cor}
%
%
%
\subsection{The unit regularity of $\hat{G}$}
%
%
%
For a monoid $M$ denote by $M^\times$ its unit group, and by $Idem(M)$ its set of idempotents. $M$ is called unit regular, 
if $M=M^\times Idem(M)=Idem(M)\, M^\times$.
\begin{Theorem}\mb{}\label{GD20}\\
a) The Kac-Moody group $G$ is the unit group of $\GD$. \vspace*{1ex}\\
b) The set of idempotents of $\GD$ is given by
\begin{eqnarray*}
  \mb{Idem}(\GD) &=& \Mklz{ g e(R) g^{-1} }{g\in G\,,\; R\in \RkX} \\ 
  &=& \Mklz{ g e(R(\Th)) g^{-1}}{ g\in G\,,\; \Th \mb{ special }}\;\;.
\end{eqnarray*}
c) We have $\GD \,=\, G \,E\, G  \,=\, G \,E_{sp}\, G\,$. 
In particular $\GD$ is unit regular.
\end{Theorem}
{\bf Remark:}
P. Slodowy already guessed $\GD=G E_{Sp} G$, D. Peterson claimed $\GD=G E_{Sp} G$, and $\mb{Idem}(\GD)=
\Mklz{ g e(R(\Th)) g^{-1}}{ g\in G\,,\; \Th \mb{ special}}$.\\\\
\Proof\\
$\bullet$ To prove c) it is sufficient to show $\GD \subseteq G E_{sp} G$. Due to the definition of $\GD$, and the formula $e(\sigma R(\Th))=n_\sigma 
e(R(\Th)) n_\sigma^{-1}$, an element of $\GD$ can be written in the form
\begin{eqnarray*}
  g_1 e(R(\Th_1))\,\cdots\, g_p e(R(\Th_p)) &\;\mb{ where } \;& g_1,\ldots,\, g_p\in G\,,\;\;\Th_1,\ldots,\,\Th_p\;\mb{ special}\;\;.
\end{eqnarray*}
We can transform this element in an element of $G E_{sp} G$ by applying $(p-1)$ times the following step, which uses the Birkhoff 
decomposition of $G$, Proposition \ref{GD17}, and Theorem \ref{GD18} b):
\begin{eqnarray*}
 \lefteqn{ e(R(\Th))u^- n_\sigma u^+ e(R(\ti{\Th}))\,\;=\,\; 
 p_\Th^* (u^-)e(R(\Th))n_\sigma e(R(\ti{\Th}))p_{\ti{\Th}}(u^+) }\\
  &=&  p_\Th^* (u^-)\, e\left( R(\Th)\cap \sigma R(\ti{\Th}) \right)n_\sigma 
  p_{\ti{\Th}}(u^+)\,\;=\,\;   p_\Th^* (u^-)n_\tau e(R(\Xi))n_\tau^{-1} 
  n_\sigma p_{\ti{\Th}}(u^+)\\
  &&\mb{where }\quad u^-\in U^-\,,\;n_\sigma\in N\,,\;u^+\in U^+\,,\mb{ and }\;\;\Th,\,\ti{\Th}\mb{ special }\;.
\end{eqnarray*}
$\bullet$ Obviously we have $(\GD)^\times \supseteq G$. To show the reverse inclusion let $ge(R)h$ be a unit. Because $(\GD)^\times$ is a group 
containing $G$, we get 
\begin{eqnarray*}
   e(R)\,\;=\,\;g^{-1}\left(ge(R)h\right)h^{-1}\;\,\in\,\;(\hat{G})^\times\;\;.
\end{eqnarray*}
Because of $e(R)^2=e(R)$ we conclude 
\begin{eqnarray*}
   e(R)\,\;=\,\;e(R)^2 e(R)^{-1}\;\,=\,\;e(R) e(R)^{-1}\;\,=\,\;1\;\;.
\end{eqnarray*}
Therefore $ge(R)h = gh \in G$. \vspace*{1ex}\\
$\bullet$ Obviously $\mb{Idem}(\GD)  \supseteq  \Mklz{ g e(R(\Th)) g^{-1}}{ g\in G\,,\; \Th\mb{ special} }$. To show the reverse inclusion let 
$ge(R(\Th))h$ be idempotent. Then we have
\begin{eqnarray*}
 e(R(\Th))\,h\,g \,e(R(\Th)) &=& e(R(\Th))\;\;.
\end{eqnarray*}
Using the Birkhoff decomposition of $G$, we can write $hg$ in the form
\begin{eqnarray*}
  hg\;=\; v n_\sigma u &\mb{ with }& v\in U^-\,,\;n_\sigma\in N\,,\,u\in U^+\;\;.
\end{eqnarray*}
To cut short the notation set $p := p_\Th$, $p^* := p_\Th^*$. Due to Proposition 
\ref{GD17} and Theorem \ref{GD18} b), we find
\begin{eqnarray}
  e(R(\Th)) &=& e(R(\Th))\,v \,n_\sigma \,u \,e(R(\Th))\;=\; 
  p^* (v)e(R(\Th))\,n_\sigma \,e(R(\Th))p(u) \nonumber \\
  &=&  p^* (v) \,e\left( R(\Th)\cap \sigma R(\Th) \right) n_\sigma p(u) \;\;.\label{penp}
\end{eqnarray}
Due to Corollary \ref{GD19} the faces $R(\Th)$ and $R(\Th)\cap\sigma R(\Th)$ are $\We$-conjugate. In particular they are of the same dimension. 
Because $R(\Th)\cap\sigma R(\Th)$ is contained in $R(\Th)$ these two faces coincide. We conclude $R(\Th)\subseteq\sigma R(\Th)$. Applying the same 
argument once more, we find $R(\Th)=\sigma R(\Th)$. Therefore $\sigma\in\We_{\Th\cup\Th^\bot}$.\\ 
Inserting in (\ref{penp}) and using Theorem \ref{GD18} a), we get 
\begin{eqnarray*}
  1\;\,=\,\;   p(p^*(v))p^*(n_\sigma p(u))\,\;=\,\;p^*(v)p^*(n_\sigma) p(u)\;\;.
\end{eqnarray*}
Now we have
\begin{eqnarray*}
  ge(R(\Th))h &=& (gu^{-1})u e(R(\Th))v n_\sigma (gu^{-1})^{-1}\\
  &=& (gu^{-1})p(u)e(R(\Th))p^*(v)p^*(n_\sigma)p(u)p(u)^{-1}(gu^{-1})^{-1} \\
  &=& (gu^{-1})e(R(\Th))(gu^{-1})^{-1}\;\;.
\end{eqnarray*}
\End\\
{\bf Remark:} With part c) of the last theorem we have reached the description of the elements of $\GD$ by expressions of the
form (\ref{expression}) of short length. Using Theorem \ref{GD18} we can decide if two such 
expressions give the same element. We have even reached more. Using the length reduction step given in the proof of c), 
we are able to compute the product of two such expressions.
But we have to work with the projections $p_\Th\,$, which can be difficult.\\ 
%
%
%
\subsection{The Weyl monoid $\hat{\cal W}$ and the monoids $\hat{T}$, $\hat{N}$}
%
%
%
In this subsection we first introduce and investigate the Weyl monoid $\hat{\We}$, which plays the same role for the monoid $\GD$ as the Weyl group 
$\We$ for the Kac-Moody group $G$. It has similar structural properties as the Renner monoid of a reductive algebraic group. In our 
further investigations we will see, that it is really the analogue of a Renner monoid.\vspace*{1ex}\\
The Weyl group acts on the monoid (\,$\RkX\,,\,\cap$\,). The semidirect product $\RkX\rtimes\We$ consists of the set 
$\RkX\times\We$ with the structure of a monoid given by 
\begin{eqnarray*}
  (R,\sigma)\cdot(S,\tau) &:=& (R\cap\sigma S,\sigma \tau)\;\;.
\end{eqnarray*}
It is easy to see that we get a congruence relation on $\RkX\rtimes\We$ by
\begin{eqnarray*}
   (R,\sigma) \sim (R',\sigma') &:\iff & R\:=\:R'\quad\mb{and}\quad\sigma'\sigma^{-1}\in Z_{\cal W}(R)\;\;.
\end{eqnarray*}
We denote the congruence class of $(R,\sigma)$ by $\ve{R}\sigma$.
\begin{Def} \label{GD21} We call the monoid $\hat{\We} := (\RkX\rtimes\We)/\sim $ the Weyl monoid.  
\end{Def}
It is easy to check, that we get embeddings of monoids
\begin{eqnarray*}
  \begin{array}{ccc}      
    \RkX &\to     &\WeD\\
      R  &\mapsto &\ve{R}1
  \end{array} 
  & \quad , \quad &
  \begin{array}{ccc}
         \We &\to     & \WeD   \\                       
      \sigma &\mapsto & \ve{X}\sigma
  \end{array}\;\;. 
\end{eqnarray*}
We denote the element $\ve{R}1$ by $\ve{R}$, and the image of $\RkX$ by ${\cal E}$. We identify $\We$ with its image in $\WeD$. 
\begin{Prop}\label{GD22}
$\WeD$ is a unit regular monoid, with unit group $\We$ and set of idempotents ${\cal E}$.
\end{Prop}
\Proof
Obviously $\ve{X}\sigma$ is a unit with inverse $\ve{X}\sigma^{-1}$. Now let $\ve{R}\sigma \in (\WeD)^\times$. Then there exists an 
element $\ve{S}\tau\in\WeD$, such that
\begin{eqnarray*}
  \ve{X} &=& \ve{R}\sigma\cdot\ve{S}\tau \,\;=\;\, \ve{R\cap\sigma S}\sigma\tau \;\; .
\end{eqnarray*}
Therefore we get $X=R\cap\sigma S$. Because $X$ is the biggest element of $(\RkX,\subseteq)$, we conclude $R=X$.\\
Obviously $\ve{R}$ is idempotent. Now let $\ve{R}\sigma$ be idempotent. Because of
\begin{eqnarray*}
   \ve{R}\sigma \;\,= \;\, (\ve{R}\sigma)(\ve{R}\sigma)\:\,=\,\:\ve{R\cap\sigma R}\sigma^2\;\;,
\end{eqnarray*}
we find $\sigma^2\sigma^{-1}=\sigma \in Z_{\cal W}(R)$. From this follows $\ve{R}\sigma=\ve{R}$.\\
Now the unit regularity is obvious by the definition of $\WeD$.\\
\End 
A monoid $M$ is called an inverse monoid, if for every element $m\in M$ there exists a unique element $m^{inv}\in M$, such that 
$m \,m^{inv} \,m=m$ and $m^{inv} \,m \,m^{inv}=m^{inv}$. The map $\mb{}^{inv}:M\to M$ is an involution extending the inverse map of the unit group.  
\begin{Prop}
$\WeD$ is an inverse monoid.
\end{Prop}
\Proof It is easy to check that $\ve{\sigma^{-1}R}\sigma^{-1}$ is an inverse of $\ve{R}\sigma\in \WeD$. Furthermore the idempotents
of $\WeD$ commute. Due to \cite{How}, Theorem 5.1.1, $\WeD$ is an inverse monoid.\\
\End
The following two propositions are not difficult to prove:
\begin{Prop} For $J\subseteq I$ set ${\cal E}_J:=\Mklz{\ve{R}}{R\in \RkX_J}$. Then
\begin{eqnarray*}
  \WeD_J \;\,:=\;\, \We_J {\cal E}_J\;\,=\;\,{\cal E}_J \We_J
\end{eqnarray*}
is a submonoid of $\WeD$. For $J,\,K\subseteq I$ we have $\,\WeD_J\subseteq \WeD_K$ if and only if $J\subseteq K$.
\end{Prop}
We call $\WeD_J$ a standard parabolic submonoid of $\WeD$, and every conjugate $\sigma \WeD_J \sigma^{-1}$, $\sigma\in\We$, a parabolic 
submonoid of $\WeD$. 
\begin{Prop}
$\WeD$ acts faithfully on the Tits cone by 
\begin{eqnarray*}
  \left(\sigma\ve{R}\right) \la \;\,:=\;\, \left\{\begin{array}{ccc}
                                       \sigma\la &\mb{ if }& \la\in R \\
                                         0 &\mb{ if }& \la\in X\setminus R 
                                      \end{array}\right.   &\quad,\quad& \sigma\ve{R}\,\in\,\WeD\;\;.
\end{eqnarray*} 
The parabolic submonoid $\WeD_J$ is the stabilizer of any element of the facet $F_J$, as well as the stabilizer of $F_J$ as a whole, $J\subseteq I$.
\end{Prop} 
Next we introduce the monoids $\TD$ and $\ND$. In our further investigations we will see, that $\TD$ can be described as an affine generalized toric 
variety.
\begin{Def}\mb{}\label{GD23}
Denote by $\TD$ the monoid generated by $T$ and $E$, and by $\ND$ the monoid generated by $N$ and $E$.
\end{Def}
It is not difficult to prove the following proposition, describing the structure of $\TD$ and $\ND$:
\begin{Prop}\mb{}\label{GD24}
$\TD$ is an abelian inverse unit regular monoid with unit group $T$ and set of idempotents $E\,$.\\
$\ND$ is an inverse unit regular monoid with unit group $N$ and set of idempotents $E\,$.
\end{Prop}
The Weyl group is isomorphic to $N/T$. There is a similar description for the Weyl monoid.
We get a congruence relation on $\ND$ as follows:
\begin{eqnarray*}
\quad\hat{n}\:\sim\:\hat{n}'\quad:\iff\quad\hat{n}T\:=\:\hat{n}'T\quad\iff \quad\hat{n}'\:\in\:\hat{n}T\quad\iff\quad\hat{n}\in\hat{n}'T
\end{eqnarray*}
Set $\ND/T:=\ND/\sim$. In a similar way define $\TD/T$.
\begin{Prop}\mb{}\label{GD25}
The monoid $\ND/T$ is isomorphic to $\WeD$ by:
\begin{eqnarray*}
  \kappa\,:\quad\WeD\;\:\;    & \to     &   \;\;\ND/T       \\
       \sigma \ve{R}        & \mapsto & n_\sigma e(R)T
\end{eqnarray*}
The set of idempotents ${\cal E}$ of $\WeD$ is mapped onto $\TD/T$.
\end{Prop}
\Proof
To show that $\kappa$ is well defined and injective let $R=\tau R(\Th)$, and $R'\in \RkX$. Because $\TD$ is abelian, because of  
Theorem \ref{GD16} (\ref{Glzenl}), and $Z_{\cal W}(R)=\tau \We_{\Th}\tau^{-1}$, we find:
\begin{eqnarray*}
  \lefteqn{\,n_\sigma e(R)T\:=\:\ti{n}_{\sigma'}e(R')T \quad \iff \quad 
  \exists\:t\in T \::\;\;n_\sigma^{-1}\ti{n}_{\sigma'}te(R')\:=\:e(R) }\\
  &\iff & \quad R\:=\:R'\quad\mb{and}\quad \exists\:t\in T\::\;\;n_\sigma^{-1} 
  \ti{n}_{\sigma'}t\in\tau N_\Th\tau^{-1} \\
  &\iff & \quad R\:=\:R'\quad\mb{and}\quad\sigma^{-1}\sigma' \:\in\:
  \tau\We_\Th\tau^{-1}\quad\iff\quad \sigma \ve{R}\:=\:\sigma'\ve{R'}
\end{eqnarray*}
Due to the last proposition, $\kappa$ is surjective. It is a homomorphism of monoids because of
$\kappa(e(X)) = T$, and
\begin{eqnarray*}
  \kappa(\sigma \ve{R}\cdot\tau \ve{S}) &=& \kappa(\sigma\tau \,\ve{\tau^{-1}R\cap S})
      \;\,=\;\,n_{\sigma\tau}T e(\tau^{-1}R\cap S) \\
     &=& \left(n_\sigma e(R) n_\tau e(S)\right) T \,\;=\;\,\kappa(\sigma \ve{R})
      \cdot\kappa(\tau \ve{S}) \;\;.
\end{eqnarray*}
\End
%
%
%
\subsection{Some double coset partitions of $\hat{G}$}
%
%
%
In this subsection we determine for some subgroups of $G$ easy representative systems of the corresponding double 
coset partitions of $\GD$.
\begin{Prop}\label{GD26} We have
\begin{eqnarray*}
  \GD &=& \dot{\bigcup_{\Th\;special}} G e(R(\Th)) G  \;\;.
\end{eqnarray*}
\end{Prop}
\Proof This follows by combining Theorem \ref{GD20} c) and Theorem \ref{GD18} a).\\
\End
\begin{Prop}\mb{}\label{GD27}
Let $\epsilon\in \{+,-\}$. We have
\begin{eqnarray*}
  \GD &=& \dot{\bigcup_{R\,\in\,{\cal R}(X)}} G e(R)B^\epsilon \,\;=\;\,
          \dot{\bigcup_{R\,\in\,{\cal R}(X)}} B^\epsilon e(R)G  \;\;.   
\end{eqnarray*}
\end{Prop}
\Proof
Let $\Th$ be special. By using the Bruhat and Birkhoff decompositions of $G$, Proposition 
\ref{GD17}, and Theorem \ref{GD18} b), we find
\begin{eqnarray*}
  G e(R(\Th))G &=& \bigcup_{\sigma\,\in\,\cal W} Ge(R(\Th))B^-n_\sigma B^\epsilon \,\;=\;\,
                    \bigcup_{\sigma\,\in\,\cal W} Ge(R(\Th))n_\sigma B^\epsilon \\
  &=& \bigcup_{\sigma\,\in\,\cal W} Ge(\sigma^{-1}R(\Th))B^\epsilon\,\;=\;\,\bigcup_{R \;of\;type\;\Th} Ge(R)B^\epsilon \;\;.
\end{eqnarray*}
To show that the last union is disjoint, let $Ge(\sigma R(\Th))B^\epsilon = Ge(\tau R(\Th))B^\epsilon$.
Then there exist elements $g\in G$, $b\in B^\epsilon$ such that
\begin{eqnarray*} 
  gn_\sigma e(R(\Th))n_\sigma^{-1}b &=& n_\tau e(R(\Th))n_\tau^{-1}\;\;.
\end{eqnarray*}
Using Theorem \ref{GD18} a), we conclude $b n_\tau\,\in\,n_\sigma 
(P_{\Th\cup\Th^\bot})^*$. Because of the generalized Birkhoff and Bruhat decompositions
\begin{eqnarray*}
  G &=& \dot{ \bigcup_{x\,\in\,{\cal W}/{\cal W}_{\Th\cup\Th^\bot}} }
        B^\epsilon x (P_{\Th\cup\Th^\bot})^* \;\;,
\end{eqnarray*}
we get $\tau\in\sigma\We_{\Th\cup\Th^\bot}=\sigma N_{\cal W}(R(\Th))$. Therefore $ \tau R(\Th)=\sigma R(\Th)$.\vspace*{1ex}\\
Taking into account the previous proposition we have proved the first equation. The second follows by application of $*$.\\  
\End 
\begin{Theorem}\mb{}\label{GD28} There are the Bruhat and Birkhoff decompositions 
\begin{eqnarray*}
  \GD \;\,=\;\, \dot{\bigcup_{\hat{w}\,\in\,\hat{\cal W}}} B^\eps\hat{w}B^\delta \,\;=\;\,
  \dot{\bigcup_{\hat{n}\,\in\,\hat{N}}} U^\eps\hat{n}U^\delta \;\;,\qquad\qquad\qquad\qquad\qquad\qquad\\
        \qquad\qquad\mb{ where }\quad  (\epsilon,\delta)\;=\;\underbrace{(+,+)\;,\;(-,-)}_{Bruhat}\;,\;\underbrace{(+,-)\;,\;(-,+)}_{Birkhoff}\;\;.
\end{eqnarray*}
\end{Theorem}
{\bf Remark:} D. Peterson claimed $\GD=U^\pm\ND U^\pm$.\vspace*{1ex}\\ 
\Proof
Let $\Th$ be special. We denote by $\We^{\Th\cup\Th^\bot}$ the minimal coset representatives of $\We/\We_{\Th\cup\Th^\bot}$. Due to 
\cite{Hu2}, Section 5.12, and \cite{K2}, Lemma 3.11, they can be  
characterized by $\We^{\Th\cup\Th^\bot}=\Mklz{\sigma\in\We\,}{\,\sigma\al_i \in \prW \;\mb{ for all }\;i\in\Th\cup \Th^\bot }$. 
\\
Because of $N_{\cal W}(R(\Th))=\We_{\Th\cup\Th^\bot}$, the faces of $X$, which are conjugated to $R(\Th)$, are of the form $\sigma R(\Th)$, 
$\sigma\in\We^{\Th\cup\Th^\bot}$. Due to the previous proposition we have
\begin{eqnarray*}
  \GD &=& \dot{\bigcup_{\Th\,sp.\,,\,\sigma\,\in\,{\cal W}^{\Th\cup\Th^\bot}}} G e(\sigma R(\Th)) B^\delta\;\;.    
\end{eqnarray*}
Now we transform the term $G e(\sigma R(\Th)) B^\delta$, using Proposition \ref{GD16}, \ref{GD17}, Theorem \ref{GD18} b), the generalized 
Bruhat and Birkhoff decompositions for $G$, and the Bruhat and Birkhoff decompositions for $G_{\Th^\bot}$:  
\begin{eqnarray*}
  G e(\sigma R(\Th)) B^\delta &=&  G e(R(\Th)) n_\sigma^{-1} B^\delta \;\,=\;\,
  \bigcup_{\tau\in {\cal W}^{\Th\cup\Th^\bot}} B^\eps n_\tau P_{\Th\cup\Th^\bot} e(R(\Th)) n_\sigma^{-1} B^\delta\\
  &=& \bigcup_{\tau\in {\cal W}^{\Th\cup\Th^\bot}} B^\eps n_\tau G_{\Th^\bot}T^{co}_\Th e(R(\Th)) n_\sigma^{-1} B^\delta\\
  &=& \bigcup_{\tau\in {\cal W}^{\Th\cup\Th^\bot}} B^\eps n_\tau U^\eps_{\Th^\bot}N_{\Th^\bot}U^\delta_{\Th^\bot}e(R(\Th)) 
      n_\sigma^{-1} B^\delta\\
  &=& \bigcup_{\tau\in {\cal W}^{\Th\cup\Th^\bot}} B^\eps \underbrace{n_\tau U^\eps_{\Th^\bot}n_\tau^{-1}}_{\subseteq\, U^\eps} 
      n_\tau N_{\Th^\bot}e(R(\Th))n_\sigma^{-1} \underbrace{n_\sigma U^\delta_{\Th^\bot}n_\sigma^{-1}}_{\subseteq\, U^\delta} B^\delta\\
  &=& U^\eps N e(R(\Th)) n_\sigma^{-1} U^\delta\;\,=\;\,\bigcup_{\hat{n}\,\in\,N e(\sigma R(\Th))} U^\eps \hat{n} U^\delta\;\;.   
\end{eqnarray*}
Before we show that the last union is disjoint, note that the depth of $\la\in P(\La)$ is defined as   
\begin{eqnarray*}
  d_\La(\la) &:=& ht(\La-\la)\;\,\in\;\,\Nn\;.
\end{eqnarray*}
Let $U^\eps n_\tau e(\sigma R(\Th))U^\delta  =  U^\eps\ti{n}_\eta e(\sigma R(\Th)) U^\delta$. Then there 
exist elements $u\in U^\eps$, $\ti{u}\in U^\delta$, such that
\begin{eqnarray}\label{esrthtiu}
   e(\sigma R(\Th))\ti{u}^{-1} &=& n_\tau^{-1}u \ti{n}_\eta e(\sigma R(\Th))\;\;.
\end{eqnarray}
Due to Theorem \ref{GD18} a) we have $u\ti{n}_\eta\in n_\tau(\sigma P_{\Th\cup\Th^\bot}\sigma^{-1})$. 
Taking into account the generalized Bruhat and Birkhoff decompositions
\begin{eqnarray*}
  G &=& \dot{\bigcup_{x\,\in\,{\cal W}/\sigma{\cal W}_{\Th\cup\Th^\bot}\sigma^{-1} }} 
         U^\eps x \left(\sigma P_{\Th\cup\Th^\bot}\sigma^{-1}\right)\;\;,
\end{eqnarray*}
we conclude
\begin{eqnarray*}
  \tau^{-1}\eta\;\in\;\sigma\We_{\Th\cup\Th^\bot}\sigma^{-1}\;\;.
\end{eqnarray*}
Because of $\bigcup_{\La\in P^+}L(\La)=X\cap P$, we may choose an element $\La\in P^+$, such that 
$P(\La)\cap\sigma\We_{\Th^\bot}F_\Th\neq\emptyset$. 
Note that this set is invariant under $\tau^{-1}\eta$. Furthermore choose an element $\la\,\in\, (P(\La)\,\cap\,\sigma\We_{\Th^\bot}F_\Th\,)$ 
such that
\begin{eqnarray*}
  d_\La (\la) &=& \min\Bigl(d_\La\left(P(\La)\,\cap\,\sigma\We_{\Th^\bot}F_\Th\:\right)\Bigr)\quad\mb{ if } \quad\delta \;\,= \;\,+\;\;,\\
  d_\La (\tau^{-1}\eta \la) &=& \min\Bigl(d_\La\left(P(\La)\,\cap\,\sigma\We_{\Th^\bot}F_\Th\:\right)\Bigr) \quad\mb{ if } \quad\delta \;\,= \;\,-\;\;.
\end{eqnarray*}
If we apply both sides of (\ref{esrthtiu}) on $v_\la\in L(\La)_\la\setminus\{0\}$, we get
\begin{eqnarray*}
  n_\tau^{-1}u\ti{n}_\eta v_\la &=& e(\sigma R(\Th))\ti{u}^{-1}v_\la\;\;.
\end{eqnarray*}
Comparing the components of the weight space $L(\La)_{\tau^{-1}\eta\la}$, we conclude that there exists a $q\in Q^\eps_0$ such that
\begin{eqnarray}\label{laplusq}
  \tau^{-1}\eta\la &=& \la + q \;\;.
\end{eqnarray} 
Because of the minimality of $d_\La(\la)$ for $\delta=+$, and the minimality of $d_\La (\tau^{-1}\eta \la)$ for $\delta=-$, we find $q=0$. 
Inserting in (\ref{laplusq}) we conclude
\begin{eqnarray}\label{tinvsinTh}
 \tau^{-1}\eta &\in & \sigma\We_\Th\sigma^{-1}\;\,=\;\,Z_{\cal W}(\sigma R(\Th))\;\;.
\end{eqnarray}
Applying both sides of (\ref{esrthtiu}) on any $w_\mu\in B(\sigma R(\Th))_\mu$, $\mu\in \sigma R(\Th)\cap P$, we get 
\begin{eqnarray*}
  e(\sigma R(\Th))\ti{u}^{-1}w_\mu &=& n_\tau ^{-1}u\ti{n}_\eta w_\mu\;\;.
\end{eqnarray*}
Comparing the components of the weight $\mu$, taking into account (\ref{tinvsinTh}), we find
\begin{eqnarray*}
  w_\mu &=& n_\tau^{-1} \ti{n}_\eta w_\mu\;\;.
\end{eqnarray*}
Therefore $n_\tau^{-1}\ti{n}_\eta \in Z_N(B(\sigma R(\Th)))$. Due to Proposition \ref{GD16} we have
\begin{eqnarray*}
 n_\tau e(\sigma R(\Th)) &=& \ti{n}_\eta e(\sigma R(\Th))\;\;.
\end{eqnarray*}
\End
%
%
%
\subsection{Constructing $\hat{G}$ from the root groups data system}
%
%
%
The definition of the monoid $\GD$ makes use of the admissible highest weight representations.  
In the following theorem we describe $\GD$, using only the root groups data system of $G$, and the monoid $(\RkX,\cap)$, which we identify 
with the isomorphic monoid $E$.\vspace*{1ex}\\
Denote by $G\sqcup E$ the free product (= coproduct) of the monoids $G$ and $E$. We identify $G$, $E$ with their images under the canonical
injections $G\hookrightarrow G\sqcup E$, $E\hookrightarrow G\sqcup E$.
\begin{Theorem}\mb{}\label{GD29}
Let $\sim$ be the congruence relation on $G\sqcup E$ generated by the following relations:\vspace*{0.5ex}\\
1) For every $ R\in\RkX\,,\;n_\sigma\in N\,$: $\qquad n_\sigma e(R) 
  n_\sigma^{-1}\:\sim\:e(\sigma R)$ \vspace*{0.5ex}\\
2) For every special set $\Th\subseteq I$, $x\in U_\al$, $\al\in\rW$:\\
  \hspace*{0.65em}$\begin{array}{clccc}
  a) \mb{ If} & \al\in\We_{\Th\cup\Th^\bot}(\Th\cup\Th^\bot) \,: 
  & xe(R(\Th))x^{-1} &\sim & e(R(\Th))  \\
  b) \mb{ If} & \al\in
  (\prW \setminus \We_{\Th^\bot}\Th^\bot) \cup \We_\Th\Th  \,:
  & xe(R(\Th)) &\sim & e(R(\Th))\\
  c) \mb{ If} & \al\in (\nrW \setminus \We_{\Th^\bot}\Th^\bot) 
  \cup \We_\Th\Th \,: & e(R(\Th))x &\sim & e(R(\Th))
  \end{array}$\vspace*{0.5ex}\\
Then $\GD$ is isomorphic to $ (G\sqcup E) / \sim\,$.
\end{Theorem}
\Proof
Let $\Th$ be special. At first we show that the following relations hold:\vspace*{0.5ex}\\
  \hspace*{0.65em}$\begin{array}{clccc}
  \alpha) \mb{ For every } & x\in\ti{G}_{\Th\cup\Th^\bot} : 
  & xe(R(\Th))x^{-1} &\sim & e(R(\Th))\\
  \beta) \mb{ For every } & x\in G_\Th\ltimes U^{\Th\cup\Th^\bot} : 
  & xe(R(\Th)) &\sim & e(R(\Th))\\
  \gamma) \mb{ For every } & x\in \left(G_\Th\ltimes U^{\Th\cup\Th^\bot}\right)^* :  
  & e(R(\Th))x &\sim & e(R(\Th)) 
  \end{array}$\\
\hspace*{0.85em} $\delta)$ For every $x\in P_{\Th\cup\Th^\bot}$, $y\in (P_{\Th\cup\Th^\bot})^*$ with 
$p_\Th(x)=p_\Th^*(y)$:
\begin{eqnarray*}
 \quad xe(R(\Th)) \;\sim\; e(R(\Th))y
\end{eqnarray*}  
The relations $\alpha$) are valid, because of 1), 2a), and the definition of 
$\ti{G}_{\Th\cup\Th^\bot}$.\\
An element $x\in G_\Th\ltimes U^{\Th\cup\Th^\bot}$ is a product of factors of the form\vspace*{0.5ex}\\
\hspace*{0.75em} (1) $\;u_\al$ with $\al\in\We_\Th\Th\,$,\vspace*{0.5ex}\\
\hspace*{0.75em} (2) $\;u_{\al_1}u_{\beta_1} \cdots u_{\al_m}u_{\beta_m} \cdot 
u_\gamma \cdot u_{\beta_m}^{-1}u_{\al_m}^{-1}\cdots 
u_{\beta_1}^{-1}u_{\al_1}^{-1}\;$ with $\;\al_1,\ldots\al_m\in$\\
\hspace*{0.75em} $\prW\setminus
\We_{\Th^\bot}\Th^\bot\,$, $\beta_1,\ldots\beta_m\in\prW\cap
\We_{\Th^\bot}\Th^\bot\,$, $\gamma\in\prW\setminus(\We_\Th\Th\cup
\We_{\Th^\bot}\Th^\bot)\,$,\vspace*{0.5ex}\\
where $u_\delta\in U_\delta$. To simplify the notation of the expressions (2), factors $u_\delta$ belonging to the same root 
$\delta$, at different distance from the central factor $u_\gamma$, may be different.\vspace*{1ex}\\
For an expression of the form (1), we have according to 2b): 
\begin{eqnarray*}
u_\al e(R(\Th))\sim e(R(\Th))\;\;.
\end{eqnarray*}
For an expression of the form (2), we have according to 2a), 2b):
\begin{eqnarray*}
\lefteqn{ u_{\al_1}u_{\beta_1} \cdots u_{\al_m}u_{\beta_m}\cdot 
 u_\gamma \cdot u_{\beta_m}^{-1}u_{\al_m}^{-1}\cdots 
 u_{\beta_1}^{-1}u_{\al_1}^{-1}\, e(R(\Th)) }\\
 &\sim & e(R(\Th))\,u_{\beta_1} \cdots u_{\beta_m} \cdot 
 u_{\beta_m}^{-1}\cdots u_{\beta_1}^{-1}\,\;=\,\; e(R(\Th)) \;\;.
\end{eqnarray*}
Therefore the relations $\beta$) are valid. Similarly $\gamma$) follows from 2a), 2c).\\
Using $\alpha)$, $\beta)$, and $\gamma)$, the relations $\delta)$ are proved in the same way as Proposition \ref{GD17}.\vspace*{0.5ex}\\ 
Now we can show that $(G\sqcup E)/\sim$ is isomorphic to $\hat{G}$.
Due to Theorem \ref{GD14}, Proposition \ref{GD16}, Theorem \ref{GD18} b), and the definition of $\GD$, 
there exists a surjective homomorphism of monoids
\begin{eqnarray*}
  \phi\,:\;(G\sqcup E)/\sim &\to & \hat{G}
\end{eqnarray*}
with $\phi([x])=x$ for all $x\in G$, $x\in E$.\\ 
$\phi$ is also injective. As in the proof of Theorem \ref{GD20} c), we can show
\begin{eqnarray*}
  (G\sqcup E)/\sim &=& \Mklz{ [g] [e(R(\Th)) ] [h] }
  { g,h\in G\,,\;\Th \mb{ special} }
\end{eqnarray*}
by using $\delta)$ and 1). Because of Theorem \ref{GD18} a), and $\delta$), we have
\begin{eqnarray*}
 ge(R(\Th))h \;=\; \ti{g}e(R(\ti{\Th}))\ti{h}\hspace*{25ex}\\
 \Rightarrow\quad \Th\:=\:\ti{\Th}\;\mb{ and }\;\begin{array}{c}
  \ti{g}^{-1} g\in P_{\Th\cup\Th^\bot} \\
  \ti{h}h^{-1}\in(P_{\Th\cup\Th^\bot})^* 
 \end{array} \quad\mb{with}\quad p_\Th(\ti{g}^{-1}g)\:=\:p_\Th^*
 (\ti{h}h^{-1}) \\
 \Rightarrow \quad [g] [e(R(\Th))] [h] \;=\; 
              [\ti{g}] [ e(R(\ti{\Th})) ] [\ti{h}] \;\;.\hspace*{20ex}
\end{eqnarray*}
\Ende
{\bf Remarks:} 1) The theorem allows to decide, if an (anti-)automorphism of $G$ can be extended to an (anti-)automorphism of $\GD$, 
compare \cite{M} for examples.\vspace*{0.5ex}\\ 
2) The theorem indicates how to generalize the monoid $\GD$ to groups with a system of root groups data, which have been defined in \cite{Ti}. 
%
%
%
%
%
%
\subsection{The action of $\GD$ on the admissible modules of $\cal O$} 
%
%
Recall that we denote by ${\cal O}_{adm}$ the full subcategory of $\cal O$, whose objects consist of admissible $\g$-modules. Due 
to the complete reducibility theorem, compare \cite{K2}, Corollary 10.7 or \cite{MoPi}, Section 6.5, 
every module of ${\cal O}_{adm}$ is a sum of highest weight modules $L(\La)$, $\La \in P^+$. Denote the corresponding 
category of $G$-modules also by ${\cal O}_{adm}\,$.
\begin{Prop}\label{GD36}\mb{}\\
1) Every $G$-module $V$ of ${\cal O}_{adm}$ extends to a $\GD$-module, the action of the idempotent $e(R)$, $R\in\RkX$, 
given by
\begin{eqnarray*}
   e(R)v_\la\;\,=\,\;\left\{\begin{array}{cl}
  v_\la & \la\in R  \\
  0     & \la\in X\setminus R
 \end{array}\right.\quad,\quad v_\la\in V_\la\,,\;\la\in P(V)\;\;.
\end{eqnarray*}
This extension is compatible with sums, tensor products, and submodules.\vspace*{0.5ex}\\
2) Every $G$-homomorphism between two $G$-modules of ${\cal O}_{adm}$ is also a $\GD$-ho\-mo\-mor\-phism between
the corresponding $\GD$-modules. 
\end{Prop}
\Proof
Due to its definition the monoid $\GD$ acts on $L(\La)$, $\La \in P^+$. Choose a 
decomposition of $V$ into a direct sum of such modules to define an action of $\GD$ on $V$. This action is independent 
of the chosen decomposition and it is given in 1).\\
Let $V$, $W$ be two $G$-modules of ${\cal O}_{adm}$. We have $P(V),\,P(W)\subseteq X\cap P$. If $R$ is a 
face of $X$, we find for $v_\la\in V_\la$, $w_\mu\in W_\mu$:
\begin{eqnarray*}
  (e(R)v_\la)\otimes (e(R)w_\mu) \;=\; \left\{\begin{array}{ccc}
                v_\la\otimes w_\mu & if & \la\in R \,\mb{ and } \,\mu\in R \\
                   0 &else &   \end{array}\right\}\\
  =\;\left\{\begin{array}{ccc}
                v_\la\otimes w_\mu & if & \la+\mu\in R \\
                   0 &else &   \end{array}\right\}
     \;=\; e(R) (v_\la\otimes w_\mu)\;\;.
\end{eqnarray*}
Therefore the extension is compatible with tensor products.\\
Let $\phi:V\to W$ a $G$-homomorphism. For 2) it is 
sufficient to show $\phi\circ e(R) = e(R)\circ\phi$, which can easily be checked on the weight 
spaces of $V$.\\
\End
\begin{Prop}\label{GD36+1} The monoid $\GD$ acts faithfully on $\bigoplus_{i=1}^{2n-l} L(\La_i)$.
\end{Prop}
\Proof Let $i\in\{n+1,\,\ldots,\,2n-l\}$. The modules $(L(\La_i),\pi_{\La_i})$ and $(L(-\La_i),\pi_{-\La_i})$ are one dimensional. By checking on the 
elements of $G$ and $E$, which generate $\GD$, it is easy to see, that for every $x\in\GD$ we have
\begin{eqnarray*}
   \pi_{\La_i}(x) \;\,=\;\, c(x)\, id_{L(\La_i)} &\mb{ where } & c(x)\in \F^\times\;\;,
\end{eqnarray*}
and furthermore 
\begin{eqnarray*}
   \pi_{-\La_i}(x) \;\,=\;\, \frac{1}{c(x)}\, id_{L(-\La_i)}\;\;.
\end{eqnarray*}
For $i\in\{n+1,\,\ldots,\,2n-l\}$ and $n_i\in \N$ we define $L(\La_i)^{\otimes\, (-n_i)}:=L(-\La_i)^{\otimes n_i}$. For $\La\in P^+$ we define 
$L(\La)^{\otimes\, 0}:=L(0)$.\\
Now let $x,x'\in\GD$ such that $\pi_{\La_i}(x)=\pi_{\La_i}(x')$ for all $i=1,\,\ldots,\,2n-l$. An element $\La$ of $P^+$ is of the form 
\begin{eqnarray*}
   \La\;\,=\,\;\sum_{i=1}^{2n-l} m_i\,\La_i &\mb{ where }& m_1,\,\ldots,\,m_n\in\Nn \;\mb{ and }\;m_{n+1},\,\ldots,\,m_{2n-l}\in \Z\;\;. 
\end{eqnarray*}
The module $\bigotimes_{i=1}^{2n-l} L(\La_i)^{\otimes m_i}$ is an module of ${\cal O}_{adm}$. Its submodule generated by 
$\bigotimes_{i=1}^{2n-l} (L(\La_i)_{\La_i})^{\otimes m_i}$ is an irreducible 
highest weight module of highest weight $\La$. Because of the last proposition we find $\pi_\La(x)=\pi_\La(x')$.\\
\End
%
%
%
%
%
\subsection{The submonoids $\GD_J$ $(J\subseteq I)$}
%
%
%
%
%
Let $\emptyset\neq J\subseteq I$. Choose an optimal realization for the generalized Cartan submatrix $A_J$, consisting of saturated 
sublattices $H(A_J)\subseteq H$ and $P(A_J)\subseteq P$, such that 
\begin{eqnarray*}
  (h_j)_{j\in J} \;\,\subseteq H(A_J)\;\, \quad \mb{ and }\quad (\al_j)_{j\in J} \;\,\subseteq \;\, P(A_J)\;\;. 
\end{eqnarray*}
We have $P=P(A_J)\oplus H(A_J)^\bot$. Denote by $q_J:P\to P(A_J)$ the corresponding projection.\\
We have an embedding of groups 
\begin{eqnarray}\label{emb1}
   G(A_J) &\to& G\;\;.
\end{eqnarray}
Due to Proposition \ref{GD6} we also have an embedding of lattices ${\cal R}(X_{A_J}) \to \RkX$, which induces an embedding of monoids 
\begin{eqnarray}\label{emb2}
   E(A_J) & \to & E\;\;,
\end{eqnarray}
whose image we denote by $E_J$. Our aim is to show, that the embeddings (\ref{emb1}), (\ref{emb2}) induce an embedding of the monoid 
$\widehat{G(A_J)}$ into $\GD$, and to describe the part $\GD_J$ of the image, which is independent of the chosen sublattices $H(A_J)$ 
and $P(A_J)$. The results of this subsection will be used later.\vspace*{1ex}\\
The next proposition is well known.
\begin{Prop}\label{GD30-1} Let $\La\in P^+$, and regard the $\g$-module $L(\La)$ as a $\g(A_J)$-module.\vspace*{0.5ex}\\
1) Define an equivalence relation on $P(\La)$ by 
\begin{eqnarray*}
   \la\sim\mu  &:\iff&    \la-\mu\in Q_J \;\;.
\end{eqnarray*}
For every equivalence class $C$, the space $V_C:=\bigoplus_{\la\in C} L(\La)_\la$ is an admissible $\g(A_J)$-module of the category 
${\cal O}(\g(A_J))$.\vspace*{0.5ex}\\ 
In particular the $\g(A_J)$-module $L(\La)$ decomposes in a direct sum of admissible irreducible highest weight modules.\vspace*{0.5ex}\\
2) The $\g(A_J)$-module
\begin{eqnarray*}
   L_J(\La) &:=& U(\n_J^-)L(\La)_{\La}
\end{eqnarray*}
is an admissible irreducible highest weight module of highest weight $q_J(\La)$. 
Set $P(L_J(\La)):= \Mklz{\la\in P(\La)}{L_J(\La)\cap L(\La)_{\la}\neq \{0\} }$. Then 
\begin{eqnarray*}
  L_J(\La)\;\,=\;\, \bigoplus_{\la\in P(L_J(\La))} L(\La)_\la 
\end{eqnarray*}
is the weight space decomposition of $L_J(\La)$, where $L(\La)_{\la}$ is the weight space of weight $q_J(\la)$.\vspace*{0.5ex}\\
The module $L_J(\La)$ is an isotypical component of $V_{[\La]}$.
\end{Prop}
\begin{Def}\mb{}Let $\emptyset\neq J\subseteq I$.\label{GD30}\\
\hspace*{2.5ex} $\GD_J$ is the monoid generated by $G_J$ and $E_J$.\\
\hspace*{2.5ex} $\ND_J$ is the monoid generated by $N_J$ and $E_J$.\\
\hspace*{2.5ex} $\TD_J$ is the monoid generated by $T_J$ and $E_J$.\vspace*{1ex}\\
Set $\GD_\emptyset:=\ND_\emptyset:=\TD_\emptyset:=E_\emptyset:=\{\,e(X)\,\}$.
\end{Def}
{\bf Remarks:}\\ 
1) Similarly to Proposition \ref{GD25} we have $\WeD_J \cong \ND_J/T_J$.\vspace*{0.5ex}\\
2) It is easy to show that we have
\begin{eqnarray*}
  \GD \:=\: \GD_I\rtimes T_{rest}\quad,\quad
  \ND \:=\: \ND_I\rtimes T_{rest}\quad,\quad
  \TD \:=\: \TD_I\times T_{rest} \;\;.
\end{eqnarray*}
\begin{Prop}\label{GD31} Let $\emptyset\neq J\subseteq I$. The embeddings (\ref{emb1}) and (\ref{emb2}) induce an embedding of $\widehat{G(A_J)}$ into 
$\GD$. 
 Restricted to $\widehat{G(A_J)}_J$, this embedding is independent of the chosen sublattices, its image given by $\GD_J$.
\end{Prop}
\Proof Due to Propositions \ref{GD30-1} and \ref{GD36}, the monoid $\widehat{G(A_J)}$ acts on $L(\La)$, $\La\in P^+$. It is easy to check, that 
this action is compatible with the embeddings (\ref{emb1}), (\ref{emb2}). Therefore we get a morphism of monoids 
\begin{eqnarray}\label{GAJmorph} 
    \widehat{G(A_J)}\to \GD\;\;.
\end{eqnarray}
Note that the elements of $\GD_I$ act as identity on the modules $L(\La_i)$, $i=n+1,\,\ldots,\,2n-l$. Due to Proposition \ref{GD36+1}, the monoid 
$\GD_I$ acts 
faithfully on the sum $\bigoplus_{i\in I}L(\La_i)$. Therefore $\widehat{G(A_J)}_J$ acts faithfully on $\bigoplus_{j\in J}L_J(\La_j)$, and the 
restriction of the morphism (\ref{GAJmorph}) to $\widehat{G(A_J)}_J$ is injective. Obviously the image of $\widehat{G(A_J)}_J$ is given by 
$\GD_J$.\vspace*{1ex}\\
Because $T(A_J)$ acts faithfully on $\bigoplus_{\La\in P^+}L(\La)$, we also get an embedding of $T(A_J)$ into $T$. 
Denote by $T_{rest,\, J}$ the image of $T(A_J)_{rest}$.\vspace*{1ex}\\ 
The morphism (\ref{GAJmorph}) maps $\widehat{G(A_J)}=\widehat{G(A_J)}_J\rtimes T(A_J)_{rest}$ onto $\GD_J T_{rest,\,J}$. It is easy to check, that 
we have $\GD_J \rtimes T_{rest,\,J}$. Therefore the morphism is injective.\\
\End 
\begin{Cor}\mb{}\label{GD33} Let $\emptyset\neq J\subseteq I$.\\ 
a) The unit group and the set of idempotents of $\GD_J$ are given by 
\begin{eqnarray*}
 \;(\GD_J)^\times\;=\;G_J  \;&\;,\;&\; \mb{Idem}(\GD_J) \:=\:\Mklz{ g e(R) g^{-1} }{g\in G_J\,,\; R\in \RkX_J}\;\;.
\end{eqnarray*}
\hspace*{0.8em} We have $\GD_J \:=\: G_J E_J G_J $. In particular $\GD_J$ is unit regular. \vspace*{1ex}\\
b) There are the Bruhat- and Birkhoff decompositions 
\begin{eqnarray*}
  \GD_J &=& \dot{\bigcup_{\hat{w}\,\in\,\hat{\cal W}_J}} \,B_J^\eps\,\hat{w}\,B_J^\delta \qquad,\qquad \eps,\delta\in \{\,+\,,\,-\,\}\;\;.
\end{eqnarray*}
\end{Cor}
\Proof The statements of the corollary can be proved easily for $\GD_I$, by using the corresponding statements for $\GD$, and Remark 2) following 
Definition \ref{GD30}.\\
Now the statements for $\GD_J$ follow from the last proposition.\\
\End
\begin{Prop}\label{GD34} Let $\emptyset\neq J\subseteq I$. The monoid $\GD_J$ has a zero if and only if $J$ is special. 
In this case $e(R(J))$ is the zero of $\GD_J$.
\end{Prop}
{\bf Remark:} By using $\GD = \GD_I\rtimes T_{rest}$, we conclude that $\GD$ has a zero if and only if $A$ is nonsingular and has no component 
of finite type. In this case the zero is given by $e(c)$, where $c=\{0\}$ is the edge of the Tits cone.\\\\
\Proof Due to Proposition \ref{GD31} it is sufficient to prove this proposition only for $\GD_I$.\\
Let $I$ be special. Then the face $R(I)$ coincides with the edge $c$ of the Tits cone. We show that $e(c)$ is the zero of $\GD_I$, by checking this 
for $G_I$ and $E$: 
Due to Proposition \ref{GD16} we have $G_I e(c)= e(c)= e(c) G_I$. Because of $R\cap c=c$, we have $e(R) e(c)=e(c)= e(c) e(R)$ for all 
$R\in \RkX$.\vspace*{1ex}\\
Now suppose that $z$ is the zero of $\GD_I$. Due to the last corollary, $z$ is of the form $z=g e(R(\Th)) h$, with $g,h\in G_I$ and $\Th$ special. 
Because $z$ is the zero we get $z=g^{-1}z h^{-1}=e(R(\Th))$. We also have $T_I z = z$. Due to Proposition \ref{GD16} we get $T_I\subseteq T_\Th$, from 
which follows $\Th=I$. In particular $I$ is special.\\
\End
Now we can give the normalizers and centralizers of the highest weight spaces of the admissible highest weight modules, which we will need 
later at several places.  
\begin{Prop}\mb{}\label{GD37}
Let $J\subseteq I$. For $\La\in P^+\cap F_J$ we have:\vspace*{1ex}\\
a) $\quad N_{\hat{G}}(L(\La)_\La)\,=\,B\WeD_J B\;$.\vspace*{1ex}\\
b) $\quad Z_{\hat{G}}(L(\La)_\La)\,=\, U^J\GD_J X_\La U^J\,=\, U\ND_J X_\La U\;$, where \\
$\hspace*{0.2em} X_\La\,:=\, \Mklz{\prod_{i=1,\,\ldots,\, 2n-l,\,i\notin J}t_i(s_i)\in T_{I\setminus J}T_{rest}}
{\prod_{i=1,\,\ldots,\, 2n-l,\,i\notin J} s_i^{\La(h_i)}\,=\,1\,}\;$.
\end{Prop}
\Proof
The proofs of a) and b) are similar, using $N_{G}(L(\La)_\La) = B \We_J B = U\ti{N}_J U$ and $Z_{G}(L(\La)_\La) = U N_J X_\La U$. As an example 
we proof a):\\ 
Due to the Bruhat decomposition of $\GD$ an element $\hat{g}\in\GD$ can be written in the form $\hat{g}=u n_\sigma e(R) \ti{u}$, where 
$u,\,\ti{u}\in U$, $n_\sigma \in N$, and  $R\in \RkX$. By using $\La\in F_{J^\infty\cup J^0}\subseteq ri(R(J^\infty))$ we find
\begin{eqnarray*}
\lefteqn{\hat{g} (L(\La)_\La) \;=\; L(\La)_\La  \;\iff\;  n_\sigma e(R) L(\La)_\La \;=\;L(\La)_\La\;\iff}\\
  && \La\in R \quad\mb{and}\quad n_\sigma L(\La)_\La \;=\;L(\La)_\La 
  \;\iff\; R\supseteq R(J^\infty) \quad \mb{and}\quad n_\sigma \in \ti{N}_J  \;\;.
\end{eqnarray*}
\End
Set $L_\emptyset(\La):=L(\La)_\La$, $\La\in P^+$. Later we will also need the following proposition.
\begin{Prop}\mb{}\label{Gq4}
Let $\La\in P^+$ and $J\subseteq I\,$. Then $U^J$ fixes the points of $L_J(\La)$.
\end{Prop}
\Proof 
The case $J=\emptyset$ is obvious, let $J\neq\emptyset$. Choose $v_\La\in L(\La)_\La\setminus\{0\}$. Due to the definition of $L_J(\La)$,
it is sufficient to show the following statements:\vspace*{1ex}\\
\hspace*{0.2em} $S_0: \quad\qquad\qquad u\, v_\La \:=\: v_\La \quad \mb{for all}\quad u\in U^J  \;\;.$\vspace*{1ex}\\ 
\hspace*{0.2em} $S_n,\;n\in\N\,:\;\; u\, x_1\cdots x_n v_\La\:=\:x_1\cdots x_n v_\La\quad 
\mb{for all}\quad u\in U^J,\;x_1,\,\ldots,\, x_n\in \n_J^-  .$\vspace*{1ex}\\
The statement $S_0$ is valid. Now the induction step from $S_n$ to $S_{n+1}$, $n\in\Nn$:\\
Fix an element $u\in U^J$. Let $x_1,\ldots,\, x_n\in \n_J^-$. Let $m\in \N\,$ and 
$y_i\in \g_{\beta_i}$, $\beta_i\in (\W_J)^-_{re} $, $i=1,\ldots,\, m$. For $t_1,\ldots,\, t_m\in\F$ set
\begin{eqnarray*}
 u(t_1,\ldots,\, t_m) &:=& \exp(t_1 y_1)\cdots \exp(t_m y_m) \;\,\in\;\,P_J\;\;.
\end{eqnarray*}
Because $U^J$ is normal in $P_J$, we have $u(t_1,\ldots,\, t_m)^{-1}\,u\,
u(t_1,\ldots,\, t_m)\,\in\, U^J\,$. Using the induction assumption, we find
\begin{eqnarray*}
 \lefteqn{u\, u(t_1,\ldots,\, t_m)\,x_1\cdots x_n v_\La} \\
 &=& u(t_1,\ldots,\, t_m)\,u(t_1,\ldots,\, t_m)^{-1}\,u\,u(t_1,\ldots,\, t_m)\,
 x_1 \cdots x_n v_\La \\
 &=& u(t_1,\ldots,\, t_m)\,x_1\cdots x_n v_\La\;\;.
\end{eqnarray*}
Because the $y_1,\ldots,\, y_m$ act locally nilpotent, the left and right side of this equation is polynomial 
in $t_1,\ldots,\, t_m$. Since $|\F\,|=\infty $, the coefficients of the monomial $t_1\cdot t_2\ldots \cdot t_n$ on the left 
and on the right are equal. We get   
\begin{eqnarray*}
 u\, y_1\cdots y_m x_1\cdots x_n v_\La \;=\; y_1\cdots y_m x_1\cdots x_n v_\La\;\;.
\end{eqnarray*}
This equation is also valid, if we permute $y_1,\ldots y_m$. We conclude
\begin{eqnarray*}
 \lefteqn{u\, [y_1[y_2\cdots[y_{m-1}, y_m]\cdots]]\, x_1\cdots x_n v_\La}\\
  &=& [y_1[y_2\cdots[y_{m-1}, y_m]\cdots]]\, x_1\cdots x_n v_\La \;\;.
\end{eqnarray*}
The Lie algebra $\n_J^-$ is generated by  $\g_\beta$, $\beta\in(\W_J)^-_{re}$. Therefore we get 
\begin{eqnarray*}
  u\, y x_1\cdots x_n v_\La \;=\; y x_1\cdots x_n v_\La \quad &\mb{for} &
  \quad y,\,x_1,\,\ldots,\, x_n\in \n_J^- \;\;.
\end{eqnarray*}
\End
%
%
%
%
\subsection{The monoid $\GD$ for a decomposable matrix $A$}
%
%
%
%
Let the generalized Cartan matrix $A$ be decomposable, $A=A_{I_1}\oplus A_{I_2}$. Choose optimal realizations $H_1\supseteq\Mklz{h_i}{i\in I_1}$, 
$P_1\supseteq\Mklz{\al_i}{i\in I_1}$ and $H_2\supseteq\Mklz{h_i}{i\in I_2}$, $P_2\supseteq\Mklz{\al_i}{i\in I_2}$
of $A_{I_1}$ and $A_{I_2}$ as described in Subsection \ref{section2}. Recall that $H=H_1\oplus H_2$ and $P=P_1\oplus P_2$.\\
Due to Proposition \ref{GD31} we get embeddings of $\widehat{G(A_{I_1})}$ and $\widehat{G(A_{I_2})}$ into $\GD$. Denote by $\GD_1$ and 
$\GD_2$ its images.
\begin{Prop}\mb{}\label{GD35} We have $\GD=\GD_1 \times \GD_2$.
\end{Prop}
\Proof
The Kac-Moody algebras $\g(A_{I_1})$ and $\g(A_{I_2})$ embed into $\g$. To simplify the notation we identify these Kac-Moody algebras 
with its images. Similarly we identify $\widehat{G(A_{I_1})}$ and $\widehat{G(A_{I_2})}$ with $G_1$ and $G_2$.\vspace*{1ex}\\
If $L_i(\La_i)$ is an irreducible highest weight module of $\g(A_{I_i})$ of highest weight $\La_i\in P_i$, $i=1,2$, then 
\begin{eqnarray}\label{L1L2}
   L_1(\La_1)\otimes L_2(\La_2)
\end{eqnarray} 
is an irreducible highest weight module of $\g=\g(A_{I_1})\oplus \g(A_{I_2})$ of highest weight $\La_1+\La_2$.\\
Because of $P^+=P_1^+ \oplus P_2^+$, all admissible irreducible highest weight modules of $\g$ are obtained in this way.\vspace*{1ex}\\ 
We have $G=G(A_{I_1})G(A_{I_2})$, and due to Proposition \ref{GD7} we have $E=E(A_{I_1}) E(A_{I_2})$. Now let $g_1\in G(A_{I_1})$, 
$g_2\in G(A_{I_2})$, and $e_1\in E(A_{I_1})$, $e_2\in E(A_{I_2})$. The elements $g_1 g_2$, and 
$e_1 e_2$ act on a decomposable element $v_1\otimes v_2$ of the module (\ref{L1L2}) by 
\begin{eqnarray*}
    g_1 g_2 (v_1\otimes v_2)\;\,=\;\, g_1 v_1\otimes g_2 v_2   &\quad,\quad &      e_1 e_2 (v_1\otimes v_2)\;\,=\;\, e_1 v_1\otimes e_2 v_2  \;\;.
\end{eqnarray*} 
From this follows easily, that the elements of $\widehat{G(A_{I_1})}$ commute with the elements of $\widehat{G(A_{I_2})}$ in $\GD$, and 
$\GD = \widehat{G(A_{I_1})}\widehat{G(A_{I_2})}$.\vspace*{1ex}\\
Now let $\hat{g}_1,\hat{g}_1'\in \widehat{G(A_{I_1})}$ and $\hat{g}_2,\hat{g}_2'\in \widehat{G(A_{I_2})}$, such that 
$\hat{g}_1\hat{g}_2=\hat{g}_1'\hat{g}_2'$. Regarding this equation on the modules $L_1(\La_1)\otimes L_2(0)$, $\La_1\in P_1^+$, we find 
$\hat{g}_1=\hat{g}_1'$. Similarly we get $\hat{g}_2=\hat{g}_2'$.\\
\End
%
%
%
%
\newpage\section{An easy algebraic geometric setting}
%
%
%
%
In this section we develop an easy, infinite dimensional algebraic geometry analogous to elementary 
algebraic geometry, generalizing the topologized coordinate rings of $L(\La)$, and of 
the Kostant cones $G (L(\La)_\La)$, $\La\in P^+$, given in \cite{KP2}, Section 3A.\vspace*{1ex}\\
This algebraic geometric setting is useful, because the monoids, coordinate rings, and Lie algebras, which are obtained by generalized 
Tannaka-Krein reconstructions, fit into this context in a natural way. For an example you may look at the next section.\vspace*{1ex}\\
There is a technical point. In the next section, we characterize the monoid, which is obtained by a generalized Tannaka-Krein reconstruction from 
the category ${\cal O}_{adm}$ and its category of restricted duals, as a 
Zariski closure $\Gq$ in a very big space. We have $\GD\subseteq \Gq$, and one of our main aims in Section \ref{sectionGqGD} is to show 
equality. But this is a long way. We have to work with $\GD$ before we really know, that $\GD$ is closed. For this reason
we have to formulate our algebraic geometric concepts also for nonclosed varieties.\vspace*{1ex}\\
Most proofs in this section are easy, or parallel to elementary algebraic geometry. They are omitted except the proofs of some propositions,
which will be used later in central position. Some more details can be found in \cite{M}.\vspace*{1ex}\\
In this section $\K$ denotes a field with $|\K\,|=\infty\,$, and all vector spaces are vector spaces over $\K\,$. 
%
%
%
%
\subsection{Varieties and pnc-varieties\label{s4.1}}
%
%
\subsubsection*{The category of (pnc-)varieties:}
The category of (pnc-)varieties is a full subcategory of the category, whose objects are topological algebras $\KK{A}$, consisting of  
$\K$-valued functions defined on a set $A\neq \emptyset$, $\K$ endowed with the discrete topology, a basis of neighborhoods of the zero given 
by a filter of ideals ${\cal F}_A\,$, (which means a filter in the lattice of ideals of $\KK{A}$).\\
A morphism of two such objects $(A,\KA,\FA)$, $(B,\KB,\FB)$ is a map 
$\phi:A\to B$, such that the comorphism $\phi^*:\KB\to\KA$ exists and is 
continuous, i.e., $(\phi^*)^{-1}\left(\FA\right)\subseteq \FB $.\vspace*{1ex}\\
In our notation we emphasize the filter of ideals more than the topology, because some of the following constructions involve rather operations 
with the ideals of the filter than with the sets of the topology.\vspace*{1ex}\\
The varieties of our algebraic geometry are constructed as follows:\vspace*{1ex}\\
Let $V$ be a vector space, $V^*$ its dual, and $X\subseteq V^*$ a subspace, which separates the points of $V$.\\ 
The algebra generated by $X$ gives a coordinate ring $\KXV$ on $V$, which is isomorphic to the symmetric algebra in $X$. The filter of ideals 
$\FV$ is given by the filterbase of ideals  
  \begin{eqnarray*}
    \Mklz{\:\Vi{U}\;}{\; U \mb{ a finite dimensional subspace of } V.\;}\;,
  \end{eqnarray*}
where $I(U)$ denotes the ideal of functions of $\KXV$ vanishing on $U$.\vspace*{1ex}\\
A set $\emptyset\neq A\subseteq V$ is equipped with a coordinate ring together with a filter of ideals by restriction:
\begin{eqnarray*}
   \KK{A} \;:=\; \KXV\res{A} &\quad,\quad &
    {\cal F}_A \;:=\; \{\,I\res{A}\,\mid \,I\in {\cal F}_V\,\}\;\;.
\end{eqnarray*}
Denote by $\Spm \KA$ the set of 1-codimensional ideals of $\KA$, and define
\begin{eqnarray*}
  \FSpm \KA &:=& {\cal F}_A\cap \Spm \KA \;.
\end{eqnarray*}
(If $\K$ is endowed with the discrete topology, these ideals are the kernels of the 
continuous homomorphisms of algebras from $\KA$ to $\K\,$.) We have a bijective map
\begin{eqnarray*}
   \begin{array}{ccc} \overline{A} &\to     & \FSpm\KA \\
                              a    &\mapsto & I_{\overline{A}}(a)\res{A}\end{array}\;, 
\end{eqnarray*} 
where $\overline{A}$ is the Zariski closure of $A$, and $I_{\overline{A}}(a)$ the 
ideal of functions of $\KK{\overline{A}}$, which vanish in $a$. In particular the map 
\begin{eqnarray*}
        A   &\to     & \FSpm\KA \\
        a   &\mapsto &  \quad\, I_A(a)
\end{eqnarray*}
is bijective if and only if $A$ is Zariski closed in $V$. Due to this fact, the 
following definitions are meaningful:
\begin{Def}\mb{}\label{V1}
An algebra of functions topologized by a filter of ideals isomorphic to such an object
$(A,\KA,\FA)$ is called a pnc-variety. It is called a variety if $A$ is Zariski 
closed. (``pnc'' means ``possibly not closed''.)
\end{Def}
{\bf Remarks:}\\
1) Due to the filter of ideals every pnc-variety can be completed to a variety.\vspace*{1ex}\\
2) If $V$ is of infinite dimension there are many possibilities to choose a subspace $X$ of $V^*$, which separates the points of $V$, 
and therefore there are many possible coordinate rings of $V$. This flexibility will be of great importance.\vspace*{0.5ex}\\ 
The situation is different if we restrict to a finite dimensional subspace $U$ of $V$. Because $X\res{U}$ separates the points of 
$U$, there is only the possibility $X\res{U}\,=U^*$, and $\KK{U}$ is the classical coordinate ring of $U$.
The topology on $\KK{U}$ is discrete because of $I(U)\res{U}\,=\{0\}\in{\cal F}_U$.\\
A map between two discrete topological spaces is always continuous. Therefore the
category of varieties includes the category of the classical affine algebraic 
varieties.
\subsubsection*{Tangent spaces and tangent maps:}
Define the tangent space of a (pnc-)variety $(A,\KA,\FA)$ at a point $a\in A$ by
\begin{eqnarray*}
      T_a A &:=& \left\{\:\delta\in\mb{Der}_a\KK{A}\,\mid \,\exists\, I\in \FA:\; I\subseteq 
     \:\mb{kernel}\:\delta\:\right\}\;\;.
\end{eqnarray*}
(If $\K$ is endowed with the discrete topology, these are the derivations in $a$ which are 
continuous.)\\
If $(A,\KA,\FA)$, $(B,\KB,\FB)$ are pnc-varieties, $a\in A$, and $\phi: A\to B$ is a morphism, we get
a linear map, the tangent map at $a$, by:
\begin{eqnarray*}
 T_a\phi:\;\;T_a A &\to & T_{\phi(a)}B \\
  \delta\quad  &\mapsto & \delta\circ\phi^*
\end{eqnarray*}
{\bf Example:}
Let $V$ be a vector space, and $(V,\KXV,\FV)$ the variety given by a subspace $X\subseteq V^*$,
which separates the points of $V$.\\ 
The tangent space $T_a V$ at $a\in V$ can be identified with $V$ by means of the linear
bijective map $V\to T_a V\,$, which assigns to an element $v\in V$ the derivation $\delta_v\in T_a V$ defined by
  \begin{eqnarray*}
    f(a+tv)\:=\:f(a)+t\delta_v (f) +O(t^2)&\quad,\quad & t\in \K\;,\quad f\in\KXV\;\;.
  \end{eqnarray*}
In this way, the tangent space $T_a V$ is independent of the chosen point separating subspace $X\subseteq V^*$.
\subsubsection*{Various constructions:}
It is not difficult to see, that the following constructions do not leave the category of (pnc-)varieties:\vspace*{0.5ex}\\
{\bf 1) Substructures:} Let $(B,\KB,\FB)$ be a pnc-variety and $a\in A\subseteq B\,$. The pnc-variety $(A,\KA,\FA)$ is 
defined in the obvious way by restricting the functions of 
$\KB$ onto $A\,$. If $B$ is a variety and $A$ is Zariski closed in $B$, then $A$ is also a variety.\\ 
The inclusion map $\iota: A\to B$ is a morphism, and its tangent map 
$T_a\iota:\,T_a A\to T_a B$ is injective with the image
\begin{eqnarray*} 
T_a\iota(T_a A)\; =\;\Mklz{\delta\in T_a B}{\mb{ kernel }\delta\supseteq I_B(A)}\;. 
\end{eqnarray*}
(Here $I_B(A)$ denotes the vanishing ideal of $A$ in $\KB\,$.)\vspace*{1ex}\\  
{\bf 2) Products:} Let $(C,\KC,\FC)$, $(D,\KD,\FD)$ be (pnc-)varieties, and $c\in C$, $d\in D$.
The algebra $\KC\otimes\KD$ topologized by the filterbase of ideals
\begin{eqnarray*}
     \Mklz{(I_C\otimes 1\, ,\,1\otimes I_D)\;}{\;I_C\in \FC\,,\,I_D\in \FD}
  \end{eqnarray*}
(the bracket ( ) stands for ``the ideal generated by'') gives in the obvious way 
a (pnc-)variety $(\CD, \KCD, \FCD)$ on $\CD\,$, which is, together with the projections $pr_C$, $pr_D$, the 
product of $(C,\KC,\FC)\,$, $(D,\KD,\FD)$ in the category of (pnc-)varieties.\\
The map $T_c\, pr_C\times T_d\, pr_D:\,T_{(c,d)}(C\times D) \to T_c C \times T_d D$ is bijective. Explicitely, the inverse of $\delta_c\in T_c C$, 
$\delta_d\in T_d D$ is determined by
\begin{eqnarray*}
  \left( (T_c\, pr_C\times T_d\, pr_D)^{-1}(\delta_c,\delta_d)\right)(f\otimes 1) &=& \delta_c(f) \quad,\quad f\in\KC \;\;, \\
  \left( (T_c\, pr_C\times T_d\, pr_D)^{-1}(\delta_c,\delta_d)\right)(1\otimes g) &=& \delta_d(g) \quad,\quad g\in\KD \;\;.
\end{eqnarray*}   
{\bf 3) Principal open sets:} Let $(A,\KA,\FA)$ be a (pnc-)variety and $g\in\KA\setminus \{0\}$.  
The localization $\KA_g$ topologized by the filter of ideals
\begin{eqnarray*}
    \Mklz{ \,I^e\,}{\,I\in\FA}\;,&\mb{ where }& I^e:\;=\;\Mklz{ \frac{f}{g^n} }{ f\in I\;,\;n\in \Nn }\;\;, 
\end{eqnarray*} 
gives in the obvious way a 
(pnc-)variety $(D(g),\KK{D(g)},\FDg)$ on the principal open set $D(g):=\Mklz{ a\in A }{ g(a)\neq 0 }$.\\ 
The inclusion map $j:D(g)\to A$ is a morphism, and for $a\in D(g)$ the tangent map  
$ T_a j:\,T_a (D(g))\to T_a A $ is bijective. Its inverse map is given explicitely by
\begin{eqnarray*}
   (T_a j)^{-1}(\ti{\de}) (\frac{f}{g^n}) \;\,= \;\,\frac{\ti{\de}(f)g(a)^n-f(a)\ti{\de}(g^n)}{g(a)^{2n}} &,& \ti{\de}\in T_a A\;,\;
   \frac{f}{g^n}\in\KK{D(g)}\;\;.
\end{eqnarray*}
\subsubsection*{Some propositions:}
The following two propositions on products and principal open sets of subvarieties will be important later:
\begin{Prop}\mb{}\label{V2}\\
1) Let $A_1,\,\ldots,\, A_k$ be sub-pnc-varieties of a pnc-variety $B$. Let $m: A_1\times \ldots \times 
A_k \to B$ be a surjective map, such that for every $i\in\{1,\,2,\,\ldots,\,k\}$ and for every element $a_i\in A_i$, there exist elements
$a_j\in A_j$, $j\neq i$ with $m(a_1,\,\ldots,\, a_i,\,\ldots,\, a_k)=a_i$. Suppose further that 
the comorphism $ m^*:\KK{B}\to \KK{A_1} \otimes \ldots \otimes \KK{A_k}$ exists and is surjective.\\ 
Then $m$ and $m^*$ are bijective, and $m^{-1}$ is a morphism of pnc-varieties.\vspace*{0.5ex}\\
2) Let $A_1,\,\ldots,\, A_k$ be varieties, $B$ a pnc-variety, and $\phi : B\to A_1\times \ldots \times A_k $ a bijective morphism 
with bijective comorphism $\phi^*\,$. Then $B$ is also a variety. 
\end{Prop}
\Proof\\
{\bf To 1):} Clearly $m^*$ is injective, because $m$ is surjective. To show the injectivity of $m$ let
$a_1,a_1'\in A_1 ,\, \ldots ,\,a_k,a_k'\in A_k$, such that $m(a_1,\,\ldots,\, a_k)\,=\,m(a_1',\,\dots,\, a_k')$. Let 
$i\in \{1,\ldots k\}$. For all $f\in\KK{A_i}$ we have
\begin{eqnarray*}
  f(a_i)\;=\; \Bigl((m^*)^{-1}(1\otimes \dots\otimes 1\otimes \underbrace{f_i}_{i\mb{-}th\, place}
  \otimes \,1 \otimes \dots \otimes 1)\Bigr)\Bigl(\underbrace{m(a_1,\,\ldots,\, a_k)}_{=\,m(a_1',\,\ldots,\, a_k')}\Bigr)\;=\; f(a_i')\;\;.
\end{eqnarray*}
Because $\KK{A_i}$ separates the points of $A_i$, we get $a_i=a_i'$.\vspace*{0.5ex}\\
It remains to show that $m^{-1}$ is a morphism. Due to the definition of the filter of ideals for the product of the 
pnc-subvarieties, it is sufficient to show for $I\in \FB$ and $i\in\{1,\,\dots,\, k\}$ the inclusion
\begin{eqnarray*}
  \Bigl((m^{-1})^*\Bigr)^{-1}(I) &\supseteq & 1\otimes \ldots \otimes 1\otimes \underbrace{I\res{A_i}}_{i\mb{-}th\,place }\otimes \, 1 \otimes \ldots 
  \otimes 1\;\;.
\end{eqnarray*}
Let $\tau_i:A_i\to B$ be the inclusion map, and define $\al_i := pr_i\circ m^{-1}: B \to A_i$. 
Because of $\al_i\circ\tau_i = id$ we find
\begin{eqnarray*} 
  \tau_i^*(I) &=&
 \left((\al_i\circ \tau_i)^*\right)^{-1}\left(\,\tau_i^*(I)\,\right) \,\;=\;\, 
 (\al_i^*)^{-1}\Bigl((\tau_i^*)^{-1}\left(\tau_i^*(I)\right)\Bigr) \\
 &=& (\al_i^*)^{-1}(I)\,+\,\underbrace{(\al_i^*)^{-1}\left(\mb{kernel\,} \tau_i^*\right)}_{=\,\{0\}} \;\;.
\end{eqnarray*}
Using this equation and the definition of $\al_i$, we get
\begin{eqnarray*}
\lefteqn{ (m^{-1})^*\Bigl(\,1\otimes \ldots 1\otimes \underbrace{I\res{A_i}}_{i\mb{-}th\,place }
\otimes \, 1\ldots \otimes 1\,\Bigr) }\\ 
 &&=\;\,(m^{-1})^*\Bigl(\,pr_i^*\left(\tau_i^*(I)\right)\,\Bigr) \;\,=\;\,
 \al_i^*\Bigl(\,(\al_i^*)^{-1}(I)\,\Bigr)\;\,\subseteq\;\, I\;\;.
\end{eqnarray*}
{\bf To 2):} We have to show that the map $B\to \FSpm \KK{B}$ is surjective.
Let $J\in\FSpm\KK{B}$. Because $\phi$ is a morphism we have
\begin{eqnarray*}
  (\,\phi^*\,)^{-1}(J) &\in & \FSpm\Bigl(\,\KK{A_1}\otimes
\ldots\otimes\KK{A_k}\,\Bigr)\;.
\end{eqnarray*}
$A_1\times \ldots\times A_k$ is a variety because $A_1,\,\ldots ,\,A_k$ are varieties.
Therefore there exist $a_1\in A_1,\,\ldots, \,a_k\in A_k$ such that
\begin{eqnarray*} 
(\,\phi^*\,)^{-1}(J)\,=\,I_{A_1\times\cdots \times A_k}\left(\,(a_1,\,\ldots,\, a_k)\,\right)\;\;.
\end{eqnarray*}
Applying $\phi^*$ on both sides we get $J=I_B\left(\phi^{-1}(a_1,\,\ldots,\, a_k)\right)$.\\ 
\End
\begin{Prop}\mb{}\label{V3}
Let $(B,\KB,\FB)$ be a variety and $g\in\KB\,$. Let $A\subseteq B$ be a 
sub-pnc-variety, such that $\overline{A}=B$.\vspace*{1ex}\\
If $\left(D_A\left(g\res{A}\right),\K\,[\,D_A\left(g\res{A}\right)\,],
{\cal F}_{D_A\left(g\,\res{A}\right)}\right)$ is a variety, then we have $D_A\left(g\res{A}\right)=D_B(g)$.
\end{Prop}
\Proof
We only have to show $D_A\left(g\res{A}\right)\supseteq D_B(g)$. 
It is easy to see that the map 
\[ \begin{array}{cccc}
\phi: &    \KB_g      & \to     & \KA_{g\,\res{A}} \\
     & \frac{f}{g^n} & \mapsto & \frac{ f\,\res{A} }{\left(g\,\res{A}\right)^n}
\end{array} \] 
is a well defined, surjective homomorphism of algebras. It is injective 
because of $\overline{A}=B\,$. $\phi$ identifies the filter of ideals ${\cal F}_{D_B(g)}$ 
and ${\cal F}_{D_A(g\,\res{A})}\,$, and therefore also $\FSpm \K\,[\,D_B(g)\,] $ and 
$\FSpm \K\,[\,D_ A(g\res{A})\,]\,$.\\
Now let $x\in D_B(g)$. Because $D_A(g\res{A})$ is a variety, there exists an element $y\in D_A(g\res{A})$ such that
\begin{eqnarray*}
  \phi\left(I_{D_B(g)}(x)\right) &=& I_{D_A(g\,\res{A})}(y)\;.
\end{eqnarray*}
Applying $\phi^{-1}$ on both sides we find $I_{D_B(g)}(x)= I_{D_B(g)}(y)$. We conclude $\quad x = y \in D_A(g\res{A})$.\\
\End\\
Let $V$ be a vector space, and $(V,\,\KXV,\,\FV)$ the variety given by a subspace $X\subseteq V^*$,
which separates the points of $V$. Let $a\in A\subseteq V$.\vspace*{1ex}\\
$\bullet$ The Zariski closure of $A$ is given by $\overline{A}=\Vm{\Vi{A}}\,$. 
\vspace*{1ex}\\
$\bullet$ Identify $T_a A\,$, $T_a\overline{A}$ with a subspace of $T_a V $ using the tangential maps 
of the inclusions, and identify $T_a V$ with $V$ using the example of above. Then
\begin{eqnarray*}
  T_a \overline{A}\;=\;T_a A\;=\;\Mklz{v\in V\:}{\:\de_v\left(\Vi{A}\right)
  = 0\,,\quad \begin{array}{c}\de_v \textrm{ the derivation of }T_a V\\
                    \textrm{ belonging to } v\,.\end{array}  }\;\;.
\end{eqnarray*} 
If the vanishing ideal $\Vi{A}$ of $A$ in $\KXV$ is not known, but the coordinate ring $\KA$ is 
known very well, then nevertheless the Zariski closure and the tangent space can be determined:
\begin{Prop}\mb{}\label{V4}We have:\vspace*{0.65ex}\\
a) $\overline{A}\; =\;\Mklz{v\in V}{\exists\;\beta_v\in \mb{Alg-Hom}(\KK{A},\K) 
 \;\;\,\forall f\in X\:: \;\;\beta_v\left(f\res{A}\right)\,=\,f(v)}\;\;.
\vspace*{0.5ex}$\\
\hspace*{0.75em} Furthermore for all $v,v'\in\overline{A}\,$:
\qquad $\beta_v\,=\,\beta_{v'}\quad\iff\quad v\,=\,v'\;\;$.\vspace*{1ex}\\
b) $T_a\overline{A}\;=\;T_a A\;=\;\Mklz{v\in V}{\exists\:\epsilon_v\in 
\mb{Der}_a\KK{A} \;\;\,\forall f\in X\:: \;\;\epsilon_v\left(f\res{A}
\right)\,=\,f(v)}\;\;.\vspace{0.5ex}$\\
\hspace*{0.75em} Furthermore for all $v,v'\in T_a\overline{A}\,$:
\qquad $\epsilon_v\,=\,\epsilon_{v'}\quad\iff\quad v\,=\,v'\;\;$.
\end{Prop}
\Proof
We only show a), the proof of b) is analogous.
$v \in \overline{A} $ is equivalent to $\Vi{A} \subseteq \Vi{v}$. This is equivalent
to the fact, that the evaluation map ${\tilde{\beta}}_v: \KXV\to\K$ in $v$
factors to a homomorphism of algebras $\beta_v:\KA\to\K\,$, with $\beta_v(f\res{A})=
\ti{\beta}_v(f)=f(v)$ for all $f\in X\,$.\\
If $\beta_v=\beta_{v'}$, then $ v=v'$, because $X$ separates the points of $V$.
If $v=v'$, then $\beta_v=\beta_{v'}$, because $X\res{A}$ generates the algebra $\KA$.\\
\End
%
%
%
\subsection{Weak (pnc-)algebraic monoids\label{s4.2}}
%
%
%
%
Let $M$ be a monoid, and let $(M,\KK{M},{\cal F}_M)$ be a (pnc-)variety, such that for all $m\in M$ the right and left translations $r_m$, 
$l_m$ are morphisms of (pnc-)varieties. We get an injective linear map 
\begin{eqnarray*}
  \Psi_{l}:\, T_1 M &\to & Der(\KK{M})
\end{eqnarray*}
by $(\Psi_{l}(\delta)f)(m):= \delta(l_m^* f)$, $f\in\KK{M}$, $m\in M$. Its image is given by  
  \begin{eqnarray*}
       \Mklz{ \,\eps\in\mb{Der}(\KK{M})\; }{\; \begin{array}{cc}\forall\;m\in M\;: & l_m^*\circ\eps\,=\,\eps\circ l_m^* \\
      \exists\;J\in {\cal F}_M\;: & \eps(J)\,\subseteq\, \Vi{1} \end{array}}\;\;.
\end{eqnarray*}
Similarly, we get an injective linear map 
\begin{eqnarray*}
  \Psi_{r}:\, T_1 M &\to & Der(\KK{M})
\end{eqnarray*}
by $(\Psi_{r}(\delta)f)(m):= \delta(r_m^* f)$, $f\in\KK{M}$, $m\in M$. Its image is given by  
  \begin{eqnarray*}
      \Mklz{ \,\eps\in\mb{Der}(\KK{M})\; }{\; \begin{array}{cc}\forall\;m\in M\;: & r_m^*\circ\eps\,=\,\eps\circ r_m^* \\
      \exists\;J\in {\cal F}_M\;: & \eps(J)\,\subseteq\, \Vi{1} \end{array}}\;\;.
\end{eqnarray*}
%
We call $(M,\KK{M},{\cal F}_M)$ a weak (pnc-)algebraic monoid, if the images of $\Psi_{l}$ and $\Psi_{r}$ are 
Lie subalgebras of $Der(\KK{M})$, and if the Lie algebra structures on $T_1M$, obtained by pulling back, are opposite.\\ 
The tangent space $T_1 M$ equipped with the structure of the Lie algebra, which is obtained by pulling back by $\Psi_{l}$, 
is called the Lie algebra Lie(M) of $M$.\\
%
%
%
A morphism of weak (pnc-)algebraic monoids $(M,\KK{M},{\cal F}_M)$, $(N,\KK{N},{\cal F}_N)$ consists of a homomorphism of monoids 
$\phi:M\to N$, which is also a morphism of the (pnc-)\-va\-rie\-ties. The map 
  \begin{eqnarray*}
     Lie(\phi)\,:=\,T_1\phi: \;Lie(M) &\to & Lie(N)
  \end{eqnarray*}
is a homomorphism of Lie algebras.\vspace*{1ex}\\
Let $(M,\KK{M},{\cal F}_M)$ be a weak (pnc-)algebraic monoid with an involution of monoids $*:M\to M$, which is also an isomorphism of 
(pnc-)varieties. 
This involution induces an involution $*$ of the algebra $\KK{M}$ by $f^*(m):=f(m^*)$, $f\in\KK{M}$, $m\in M$. It is easy to check that we 
also get an involution $*$ of the Lie algebra $Lie(M)$ by $\delta^*:=\delta\circ *\,$, $\delta\in Lie(M)$.\vspace*{1ex}\\
{\bf Examples:} In the following examples a dot $\cdot$ marks the place where the arguments of a function have to be put in.\\
{\bf 1)} Let $V$ be a vector space, and $(V,\KXV,\FV)$ the variety given by a subspace $X\subseteq V^*$, which separates the points of $V$.\\   
For an endomorphism $\phi:V\to V$ denote by $\phi^*:V^*\to V^*$ its adjoint map. We get a subalgebra of $End(V)$ by 
\begin{eqnarray*}
   End_X(V) &:=& \Mklz{\phi\in End(V)}{ \phi^*(X)\subseteq X}\;\;.
\end{eqnarray*}
Equip $End_X(V)$ with the variety structure given by the point separating subspace
\begin{eqnarray*}
     X_{End} &:=& span\Mklz{f_{\al v}:=\al(\,\cdot\, v)}{ \al\in X\,,\,v\in V}
\end{eqnarray*}
of $(End_X(V))^*$. It is not difficult to check, that $End_X(V)$ together with the concatenation of linear maps is a weak algebraic monoid.\\ 
By the identification of $T_1 End_X(V)$ with $End_X(V)$, the Lie algebra $Lie(End_X(V))$ identifies with the Lie algebra associated to the associative 
algebra $End_X(V)$.\vspace*{1ex}\\
{\bf 2)} Let $\left(\,V\,,\,\kBl\,\right)$ be a vector space  together with a nondegenerate symmetric bilinear form, $char\,\K\neq 2$. Then 
\begin{eqnarray*}
   X &:=&\{\,f_v:=\kB{v}{\cdot\;}\,|\,v\in V\,\}\;\,\subseteq \;\,V^*
\end{eqnarray*}
is a subspace of $V^*$, which separates the points of $V$. The weak algebraic monoid $(\mb{End}_X(V),\KK{End_X(V)},{\cal F}_{End_X(V)})$ of the 
last example coincides with $(\mb{Adj}(V),\KAdjV,\FAdjV)$ defined by
\begin{eqnarray*}
   Adj(V)&:=&\Mklz{ \phi\in End (V) }{\mb{ The } \kBl\mb{-adjoint map } \phi^*\in End(V)\mb{ exists. }}  \;\;,  \\
X_{Adj(V)} &:= &\mb{span}\{\,f_{vw}:=\kB{v}{\cdot\; w}\res{Adj(V)}\,|\,v,w\in V\,\}\;\,\subseteq\;\,\left(Adj(V)\right)^*  \;\;.
\end{eqnarray*}
The adjoint map $*:Adj\to Adj$ is an involutive isomorphism of the varieties.\\\\
We omit the definition of an action of a weak (pnc-)algebraic monoid. But note that for every $v\in V$ the map 
$\psi_v: End_X(V) \to V$ defined by $\psi_v(\phi):=\phi v$ is a morphism.
%
%
%
%
%
%
%
%
%
%
\newpage
\section{A generalized Tannaka-Krein reconstruction}
%
%
%
The category ${\cal O}_{adm}$ of admissible $\g$-modules generalizes the category of finite dimensional representations of a semisimple 
Lie algebra, keeping the complete reducibility theorem. In this section we associate to ${\cal O}_{adm}$ and its category of restricted duals a 
monoid with coordinate ring, as well as the Lie algebra of this monoid, by a generalization of the Tan\-na\-ka\--Krein reconstruction.
We show that this monoid can be identified with a weak algebraic monoid. Because of this result, we can use the algebraic geometric methods 
developed in the last section to determine this monoid and its Lie algebra. In the next two sections we show, that the monoid can be identified 
with $\GD$, and its Lie algebra can be identified with the Kac-Moody algebra $\g$.\\
The generalized Tannaka-Krein reconstruction is not restricted to this particular example, and we will investigate this reconstruction in detail in 
\cite{M4}. But this example is special for the following reasons: There is a complete reducibility theorem for ${\cal O}_{adm}$. The restricted 
duals have nice properties. The Kac-Moody algebra is generated by elements, which act locally finite on all modules of ${\cal O}_{adm}$. 
Furthermore a dense part of this monoid, the Kac-Moody group, has already been investigated. Also the restriction of the coordinate ring 
of the monoid onto the Kac-Moody group has already been investigated.\vspace*{1ex}\\
${\cal O}_{adm}$ is a tensor category. To two $\g$-modules $V$, $W$ of ${\cal O}_{adm}$ we can assign a direct sum $V\oplus W$ and a tensor product 
$V\otimes W$, which are both modules of ${\cal O}_{adm}$. We also choose an one-dimensional module $V_0$ of ${\cal O}_{adm}$, on which $\g$ acts 
trivially.\vspace*{1ex}\\ 
The category of restricted duals of ${\cal O}_{adm}$ is obtained in the following way:\\
Let $(V,\pi_V)$ be a $\g$-module of ${\cal O}_{adm}$. Its restricted dual $V^{(*)}$ is defined as a linear space by 
\begin{eqnarray*}
   V^{(*)} &:=& \bigoplus_{\la\in P(V)} V_\la^*\;\,\subseteq \;\, V^*\;\;.
\end{eqnarray*}
It separates the points of $V$.
For every $x\in \g$ the dual map $\pi_V(x)^*: V^*\to V^*$ of $\pi_V(x): V\to V$ restricts to a map $\pi_V(x)^{(*)} : V^{(*)}\to V^{(*)}$,
which we call the restricted dual map. In this way $V^{(*)}$ gets the structure of a $\g^{op}$-module.\\
For every $\g$-homomorphism $\phi:V\to W$ of ${\cal O}_{adm}$, the restricted dual map $\phi^{(*)}:W^{(*)}\to V^{(*)}$ exists, and is a homomorphism 
of $\g^{op}$-modules.\vspace*{1ex}\\
Now let $V$, $W$ be $\g$-modules of ${\cal O}_{adm}$. Then in the obvious way we get linear embeddings 
\begin{eqnarray*}
   V^{(*)}\oplus W^{(*)} \;\,\to\;\,(V\oplus W)^*    &\quad,\quad &   V^{(*)}\otimes W^{(*)}\;\,\to\;\, (V\otimes W)^*  \;\;.
\end{eqnarray*}
Because the weight spaces are finite dimensional, and because of the conditions on the set of weights of the modules of ${\cal O}_{adm}$, the images of 
these embeddings coincide with $(V\oplus W)^{(*)}$ and $(V\otimes W)^{(*)}$.\vspace*{1ex}\\
Forgetting the $\g$-module structure of the objects of ${\cal O}_{adm}$, a natural transformation of the resulting category 
${\cal O}_{adm}^{fg}$ of linear spaces and linear maps consists of a family of linear maps
\begin{eqnarray*}
   m &=& (\,m_V\in End(V)\,)_{V\in {\cal O}_{adm}^{fg}}\;\;,
\end{eqnarray*} 
such that the diagram 
\begin{eqnarray*}
   V &\stackrel{m_V}{\to} & V \\
  \phi \downarrow \, &                    &  \,\downarrow \phi\\
   W &\stackrel{m_W}{\to} & W  
\end{eqnarray*}
commutes for all objects $V$, $W$, and all morphisms $\phi:V\to W$.
The natural transformations form a set, because the isoclasses of ${\cal O}_{adm}^{fg}$ form a set.\\
Induced by the algebras of endomorphisms $End(V)$, $V\in {\cal O}_{adm}^{fg}$, the set of natural transformations gets the structure of an 
associative $\F$-algebra with unit, therefore also the structure of a Lie algebra.\vspace*{1ex}\\
It is trivial to check the next proposition:
\begin{Prop} The set $M$ of natural transformations 
\begin{eqnarray*}
   m &=& (\,m_V\in End_{V^{(*)}}(V)\,)_{V\in {\cal O}_{adm}^{fg}}\;\;,
\end{eqnarray*} 
which satisfy the properties\vspace*{0.5ex}\\
\hspace*{1em} (1) $\quad m_{V\oplus W}\,=\,m_V\oplus m_W\;$ for all objects $V$, $W$,\vspace*{0.5ex}\\
\hspace*{1em} (2) $\quad m_{V\otimes W}\,=\, m_V\otimes m_W\;$ for all objects $V$, $W$,\vspace*{0.5ex}\\
\hspace*{1em} (3) $\quad m_{V_0}\, =\, id_{V_0}\;$,\vspace*{0.5ex}\\
is a submonoid of the monoid of natural transformations of ${\cal O}_{adm}^{fg}$.
\end{Prop}
{\bf Remarks:} 1) $M$ is the biggest monoid acting reasonably on the modules of ${\cal O}_{adm}$, compatible with the duals. We know already, 
that the monoid $\GD$ embeds onto a submonoid of $M$.\\
2) Recall that $m_V\in End_{V^{(*)}}(V)$ means the restricted dual map $m_V^{(*)}:V^{(*)}\to V^{(*)}$ exists. In this way $V^{(*)}$ is a 
module of the opposite monoid $M^{op}$.\\
3) It is easy to check, that the following generalization of (1) holds: Let $V_j$, $j\in J$, be objects of ${\cal O}_{adm}^{fg}$, such that 
$\bigoplus_{j\in J} V_j$ is also an object of ${\cal O}_{adm}^{fg}$. Then we have 
\begin{eqnarray*}
   m_{\,\bigoplus_{j\in J} V_j} &=& \bigoplus_{j\in J}\; m_{V_j}\;\;.
\end{eqnarray*}
The next proposition describes the coordinate ring $\FK{M}$ of the monoid $M$.
For $\phi\in V^{(*)}$, and $v\in V$ define a function $f_{\phi v}: M\to\F$, the matrix coefficient of $\phi$ and $v$, by
\begin{eqnarray*}
   f_{\phi v}(m)\;\,:=\;\, \phi(m_V v)\quad,\quad m\in M\;\;.
\end{eqnarray*}
\begin{Prop}\label{TKR2}
1) The set 
\begin{eqnarray*}
  \FK{M} &:=& \Mklz{f_{\phi v}\,}{\,\phi\in V^{(*)}\,,\, v\in V\,,\,\,V\in {\cal O}_{adm}^{fg} }
\end{eqnarray*}
is a subalgebra with unit of the algebra of functions on $M$. It separates the points of $M$. Left and right 
multiplications with elements of $M$ induce comorphisms of $\FK{M}$.\vspace*{1ex}\\
2) The monoid $M^{op}\times M$ acts on the algebra $\FK{M}$ by 
\begin{eqnarray*}
((m_1,m_2)f)(m):=f(m_1mm_2)\quad, \quad f\in\FK{M}\quad,\quad m_1,m_2,m\in M. 
\end{eqnarray*}
We get an $M^{op}\times M$-equivariant linear isomorphism 
\begin{eqnarray*}
    \bigoplus_{\La\in P^+} L(\La)^{(*)} \otimes L(\La)\to \FK{M}
\end{eqnarray*}
by assigning $\phi\otimes v$ the matrix coefficient $f_{\phi v}$, $\phi\in L(\La)^{(*)}$, $v\in L(\La)$, $\La\in P^+$.
\end{Prop}
\Proof\\
a) Let $V$, $W$ be objects of ${\cal O}_{adm}^{fg}$, and $\al:V\to W$ an isomorphism. Since the monoid $M$ consists of natural transformations, 
we have
\begin{eqnarray*}
   f_{\phi\, \al(v)}\;\,=\;\,f_{\al^{(*)}(\phi)\, v}\quad\mb{ for all }\quad \phi\in W^{(*)}\quad ,\quad v\in V\;\;.
\end{eqnarray*}
In particular the image of the map $V^{(*)}\otimes V \to \FK{M}$, which assigns $\phi\otimes v$ the matrix coefficient $f_{\phi v}$, depends 
only on the isomorphism class of $V$.\vspace*{1ex}\\
b) Next we show the properties of the set $\FK{M}$, which are stated in 1):\\ 
Let $V$, $W$ be objects of ${\cal O}_{adm}^{fg}$, and 
$v\in V$, $w\in W$, $\phi\in V^{(*)}$ and $\psi\in W^{(*)}$. Because of the properties (1) and (2) in the definition of $M$ we have
\begin{eqnarray*} 
   f_{\phi v} \; + \; f_{\psi w} \;\,=\;\, f_{\phi\oplus\psi\, v\oplus w} \quad,\quad
   f_{\phi v}f_{\psi w} \;\,=\;\, f_{\phi\otimes\psi\, v\otimes w} \;\;.
\end{eqnarray*}
Because of property (3) we have
\begin{eqnarray*}
   f_{\phi v} &=& \phi(v)\; 1\quad,\quad \phi\in V_0^{(*)}\,=\,V_0^*\quad,\quad v\in V_0 \;\;.
\end{eqnarray*}
Let $m,\,\ti{m}\in M$, and suppose $f(m)=f(\ti{m})$ for all $f\in\FK{M}$. Then for every object $V$ of ${\cal O}_{adm}^{fg}$ we have 
$\phi(m_V v)=\phi(\ti{m}_V v)$ for all $\phi\in V^{(*)}$, $v\in V$. Because $V^{(*)}$ is point separating on V, we find $m_V=\ti{m}_V$.\vspace*{1ex}\\
The left and right translations induce comorphisms, because for $v\in V$, $\phi\in V^{(*)}$, $V$ an object of ${\cal O}_{adm}^{fg}$, we have 
\begin{eqnarray*}
 l_m^* \,f_{\phi v} \;\,=\;\, f_{m_V^{(*)}(\phi)\, v}  & \quad,\quad & r_m^* \,f_{\phi v} \;\,=\;\, f_{\phi\, m_V(v)}\;\;.
\end{eqnarray*}
c) We get a linear map 
\begin{eqnarray*}
   \Phi_M :\,\bigoplus_{\La\in P^+} L(\La)^{(*)} \otimes L(\La)\to \FK{M}
\end{eqnarray*}
by assigning $\phi\otimes v$ the matrix coefficient $f_{\phi v}$, $\phi\in L(\La)^{(*)}$, $v\in L(\La)$, $\La\in P^+$.\\ 
This map is surjective:
Let $V$ be an object of ${\cal O}_{adm}^{fg}$. Due to the complete reducibility theorem, \cite{K2}, Theorem 10.7, there exists an isomorphism 
$\al:V\to \bigoplus_{j\in J} L(\La_j)$, where $\La_j\in P^+$. Let $v\in V$ and $\phi\in V^{(*)}$. Since 
$\al^{(*)}: \bigoplus_{j\in J} L(\La_j)^{(*)}\to V^{(*)}$ is bijective, there exist an element $\psi\in \bigoplus_{j\in J} L(\La_j)^{(*)}$ such 
that $\phi=\al^{(*)}(\psi)$. Due to a) we find
\begin{eqnarray*}
  f_{\phi v} \;\,=\;\, f_{\al^{(*)}(\psi)\, v}  \;\,=\;\, f_{\psi\, \al(v)} \;\;.
\end{eqnarray*}
This matrix coefficient is an element of the image of $\Phi_M$.\vspace*{1ex}\\
Denote by $\Phi_G:\bigoplus_{\La\in P^+} L(\La)^{(*)} \otimes L(\La)\to \FK{G}$ the isomorphism of the Peter and Weyl theorem for the algebra of 
strongly regular functions $\FK{G}$. Obviously we get a linear map $res:\FK{M}\to \FK{G}$ by restricting the functions of $\FK{M}$ onto $G$. 
Because of $\Phi_G = res\,\circ\, \Phi_M$, the map $\Phi_M$ is injective.\vspace*{1ex}\\  
The monoid $M^{op}\times M$ acts by homomorphism of algebras, because of $(m_1,m_2)f=(l_{m_1}^*\circ r_{m_2}^*)(f)$, $f\in\FK{M}$. The 
$M^{op}\times M$-equivariance is obvious.\\ 
\End
\begin{Cor}\label{TKR3} 1) If we identify the Kac-Moody group $G$ with the corresponding subgroup of $M$, we have:\vspace*{0.5ex}\\
$\bullet$ The Kac-Moody group $G$ is Zariski closed in $M$.\vspace*{0.5ex}\\
$\bullet$ The algebra of strongly regular functions $\FK{G}$ is isomorphic to the coordinate ring $\FK{M}$ by the restriction map.\vspace*{0.5ex}\\
2) The coordinate ring $\FK{M}$ is an integrally closed domain. In particular $M$ is irreducible. 
\end{Cor}
{\bf Remark:} In this way, the algebra of strongly regular functions is really the coordinate ring of the monoid $M$.\vspace*{1ex}\\
\Proof
With the notations of the proof of the last proposition, the restriction map  $res=\Phi_G \,\circ\, \Phi_M^{-1}$ is bijective. From this follows 
part 1) of the corollary.\\
Since $\FK{M}$ is isomorphic the algebra of strongly regular 
functions, part 2) follows from the properties of the algebra of strongly regular functions.\\
\End
The next Proposition gives the Lie algebra $Lie(M)$ of ($M$, $\FK{M}$), and the adjoint action of the unit group $M^\times$ on $Lie(M)$.
\begin{Prop} Let $Lie(M)$ be the set of natural transformations
\begin{eqnarray*}
   x &=& (\,x_V\in End_{V^{(*)}}(V)\,)_{V\in {\cal O}_{adm}^{fg}}
\end{eqnarray*}
which satisfy the following properties:\vspace*{0.5ex}\\
\hspace*{1em} (1) $\quad x_{V\oplus W}\,=\,x_V\oplus x_W\;$ for all objects $V$, $W$.\vspace*{0.5ex}\\
\hspace*{1em} (2) $\quad x_{V\otimes W}\,=\, x_V\otimes id_{W}+ id_V \otimes x_W\;$ for all objects $V$, $W$.\vspace*{0.5ex}\\
\hspace*{1em} (3) $\quad x_{V_0}\, =\, 0_{V_0}\;$.\vspace*{0.5ex}\\
\hspace*{1em} (4) There exists a derivation $\delta_x:\FK{M}\to \F$ in 1, such that $\delta_x(f_{\phi v})=\phi(x_V v)$\\  
\hspace*{3em} for all $\phi\in V^{(*)}$, 
$v\in V$, $V\in {\cal O}_{adm}^{fg}$.\vspace*{0.5ex}\\
Then Lie(M) is a Lie subalgebra of the Lie algebra of natural transformations of ${\cal O}_{adm}^{fg}$.
The unit group $M^\times$ of $M$ acts on $Lie(M)$ by conjugation. 
\end{Prop}
{\bf Remarks:}\\ 
1) The Lie algebra $Lie(M)$ is the biggest Lie algebra acting reasonably on the modules of ${\cal O}_{adm}$, compatible with 
the restricted duals, and with the coordinate ring $\FK{M}$.\\
2) Clearly we want to have the compatibility condition, that the Lie algebra of natural transformations, corresponding to the Kac-Moody algebra 
$\g$, is a subalgebra of $Lie(M)$. We will show in the last section that we have equality.\vspace*{1ex}\\ 
\Proof
Let $x,\,y\in Lie(M)$ and $m\in M^\times$. It is trivial to check, that (1), (2), and (3) are also valid for $[x,y]$ and $mxm^{-1}$.\\
Note that statement (4) is equivalent to the statement:\vspace*{0.5ex}\\
\hspace*{1em} (4') There exists a derivation $\ti{\delta}_x\in Der(\FK{M})$, such that 
$\ti{\delta}_x(f_{\phi v}) =f_{\phi\, x_V v}$\\ 
\hspace*{3em} for all $\phi\in V^{(*)}$, $v\in V$, $V\in {\cal O}_{adm}$.\vspace*{0.5ex}\\
The derivation $\ti{\delta}_x$ is the left invariant derivation of $\FK{M}$ corresponding to $\delta_x$.\\\
If $x,\,y\in Lie(M)$, then $[\ti{\delta}_x,\ti{\delta}_y]\in Der(\FK{M})$ has the pro\-per\-ties of (4') for $[x,y]$.\\
For $m\in M^\times $ the comorphisms $l_m^*$, $r_{m^{-1}}^*$ of the left and right translations are isomorphism of algebras. If $x\in Lie(M)$, then 
$\delta_x\circ l_m^*\circ r_{m^{-1}}^* $ has the properties of (4) for $mxm^{-1}$.\\ 
\End
Recall the definition of the weak algebraic monoid $\EndoX$, where $X=\bigoplus_{\La\in P^+}L(\La)^{(*)}\subseteq 
\left(\bigoplus_{\La\in P^+}L(\La)\right)^*$, given in the last section. It is easy to check, that
\begin{eqnarray*}
         \lefteqn{\grEndoX }\\     
         &:=& \Mklz{\,\phi\in\EndoX\;}{\;\phi \left(L(\La)\right)\subseteq L(\La),\;\La\in P^+\;}                                            
\end{eqnarray*}
is a closed submonoid. We will work with this weak algebraic monoid for the following reason: Equip $L(\La)$ with the variety structure given by 
$L(\La)^{(*)}\subseteq L(\La)^*$, $\La\in P^+$. Then $\grEndoX$ acts on $L(\La)$, $\La\in P^+$, but $\EndoX$ does not act.\vspace*{1ex}\\
Recall that we have identified the Kac-Moody group $G$ with a subgroup of $\Endo$. Actually we have $G \subseteq \grEndoX$, and the closure 
$\Gq$ gets the structure of a weak algebraic monoid. 
\begin{Theorem}
We get an embedding of monoids
\begin{eqnarray*}
 \Upsilon:\; M &\to & \grEndoX\;\;,
\end{eqnarray*}
with image $\Gq$, by
\begin{eqnarray*}
  \Upsilon(m)v:=m_{L(\La)}v \quad,\quad m\in M\quad,\quad v\in L(\La)\quad,\quad \La\in P^+\;\;.
\end{eqnarray*}
If we identify $M$ with $\Gq$ by $\Upsilon$, then the coordinate ring $\FK{M}$ identifies with the coordinate ring $\FK{\Gq}$, and 
the Lie algebra $Lie(M)$ identifies with $Lie(\Gq)$ by the map $x\mapsto \delta_x$.
\end{Theorem}
\Proof
First we show $\Upsilon(M)\subseteq \Gq$. Every element $m\in M$ determines an evaluation homomorphism of $\FK{M}$. Due to Corollary \ref{TKR3} the 
coordinate ring $\FK{M}$ is isomorphic to the algebra of strongly regular functions $\FK{G}$ by the restriction map. Therefore we also get an 
algebra homomorphism $\beta_m:\FK{G}\to \F$, which has the property 
\begin{eqnarray*}
    \beta_m(f_{\phi v}\res{G})\;\,=\;\, \phi(m_{L(\La)} v) &\mb{ for all }& \phi\in L(\La)^{(*)}\quad,\quad v\in L(\La)\quad,\quad \La\in P^+\;\;.
\end{eqnarray*}
Due to Proposition \ref{V4} we find $\Upsilon(m)\in\Gq$.\vspace*{1ex}\\
Denote by $\ti{\Upsilon}:M\to\Gq$ the map $\Upsilon$, whose image has been restricted to $\Gq$. We show that this map is surjective by 
constructing a map $\Omega: \Gq\to M$ with the property $\ti{\Upsilon}\circ \Omega = id_{\overline{G}}$:\\ 
Let $\hat{g}\in\Gq$. For every object $V$ of ${\cal O}_{adm}$ choose an isomorphism
\begin{eqnarray*}
 \gamma:V\to\bigoplus_{j\in J}L(\La_j)\;\;,
\end{eqnarray*}
and define a linear map $\Omega_V(\hat{g}):=\gamma^{-1} \circ\hat{g}\circ\gamma$. Because the adjoint maps of $\gamma$ and $\hat{g}$ exist, we have 
$\Omega_V\in End_{V^{(*)}}V$. The next formula shows that the map $\Omega_V(\hat{g})$ is also independent of the chosen isomorphism.\\ 
The algebra $\FK{G}$ is isomorphic to $\FK{\Gq}$ by the restriction map. Denote by $\al_{\hat{g}}:\FK{G}\to\F$ the homomorphism, which corresponds to 
the evaluation homomorphism of $\FK{\Gq}$ in $\hat{g}$. It is easy to check that we have
\begin{eqnarray*}
   \phi(\Omega_V(\hat{g})v) &=& \al_{\hat{g}}(f_{\phi v}\res{G})\quad \mb{ for all } \quad\phi\in V^{(*)}\quad,\quad v\in V\;\;.
\end{eqnarray*}
Using this formula, it is also not difficult to see, that the maps $\Omega_V(\hat{g})$, $V$ an object of ${\cal O}_{adm}^{fg}$, define an element 
$\Omega(\hat{g})$ of $M$.\vspace*{1ex}\\
It is easy to check, that the comorphism $\ti{\Upsilon}^*:\FK{\Gq}\to\FK{M}$ exists, and is given by the concatenation of the linear isomorphisms
\begin{eqnarray*}
  \FK{\Gq} \;\stackrel{res}{\to}\;\FK{G}\;\stackrel{\Phi_G^{-1}}{\to} \;\bigoplus_{\La\in P^+}L(\La)^{(*)}\otimes L(\La)\;
  \stackrel{\Phi_M}{\to}\;\FK{M}\;\;,
\end{eqnarray*}
in particular $\ti{\Upsilon}^*$ is an isomorphism.\vspace*{1ex}\\
The map $\ti{\Upsilon}$ is also injective. Let $m,m'\in M$, such that $\ti{\Upsilon}(m)=\ti{\Upsilon}(m')$. Then for all $\phi\in L(\La)^{(*)}$, 
$v\in L(\La)$, and $\La\in P^+$ we have
\begin{eqnarray*}
  \ti{\Upsilon}^*(f_{\phi v})(m)\;\,=\;\,\phi(\ti{\Upsilon}(m)v)\;\,=\;\,\phi(\ti{\Upsilon}(m')v)\;\,=\;\,\ti{\Upsilon}^*(f_{\phi v})(m')\;\;.
\end{eqnarray*}
Due to Proposition \ref{TKR2} the functions of $\FK{M}$ separate the points of $M$. Because $\ti{\Upsilon}^*$ is surjective, we find $m=m'$.
\vspace*{1ex}\\
Identify $Lie(\Gq)$ with a Lie subalgebra of the Lie algebra $\grEndoX$. It remains to show that we get an embedding of Lie algebras 
\begin{eqnarray*}
 \Upsilon':\; Lie(M) &\to & \grEndoX\;\;,
\end{eqnarray*}
with image $Lie(\Gq)$, by
\begin{eqnarray*}
  \Upsilon'(x)v:=x_{L(\La)}v \quad,\quad x\in Lie(M)\quad,\quad v\in L(\La)\quad,\quad \La\in P^+\;\;.
\end{eqnarray*}
This can proved in a similar way as the corresponding result for $\Upsilon$.\\
\End
Now we rewrite some of the last results by using the nondegenerate contravariant symmetric bilinear forms on the modules $L(\La)$, $\La\in P^+$. 
Clearly these depend on the Cartan subalgebra $\h$, which we choose by the construction of the Kac-Moody algebra. Nevertheless this is advantageous, 
because it will save work in later sections.\vspace*{1ex}\\
Recall that we have fixed a nondegenerate contravariant symmetric bilinear form $\kBl$ on every module $L(\La)$, $\La\in P^+$. We have extended these 
forms to a form on $\bigoplus_{\La\in P^+}L(\La)$, also denoted by $\kBl$, by requiring $L(\La)$ and $L(\La')$ to be orthogonal for $\La\neq \La'$.\\ 
The variety structure on $L(\La)$ of above is given by
\begin{eqnarray*}
  L(\La)^{(*)}\;\,=\;\, \Mklz{f_v:=\kB{v}{\cdot\;}\,}{\, v\in L(\La)}\;\,\subseteq\;\,(L(\La))^*\;\;.
\end{eqnarray*}
The weak algebraic monoid $\grEndoX$ can be identified with
\begin{eqnarray*}
         \lefteqn{\grAdj }\\     
         &:=& \Mklz{\,\phi\in\Endo\;}{\;\begin{array}{c}    
         \mb{The adjoint }\;\phi^* \;\mb{ exists and}\\    
         \phi \left(L(\La)\right)\subseteq L(\La),\;\La\in P^+\;.    
         \end{array}\,}  \;\;,                                          
\end{eqnarray*}
whose structure of a variety is given by
\begin{eqnarray*}
   \mb{span} \Mklz{f_{vw}:=\kkB{v}{\,\cdot\;w}\res{gr\mb{-}Adj}\,}{\,v,w\in L(\La)\,,\,\La\in P^+}\;\;.
\end{eqnarray*}
(As above a dot $\cdot$ marks the place where the arguments of a function have to be put in.)\vspace*{1ex}\\ 
The involution $*$ of $\grAdj$ restricts to the Chevalley involution on the the Kac-Moody group $G\subseteq \grAdj$. 
Because $G$ is invariant under the involution $*$, also the closure $\Gq$ is invariant under $*$.\\
We use this involution of $\Gq$, which we also call Chevalley involution, to get rid of the $\Gq^{op}$-action. We define an action 
$\pi$ of $\Gq\times\Gq$ on $\FK{\Gq}$, and an involution $*$ of $\FK{\Gq}$ by
\begin{eqnarray*}
  \begin{array}{ccc}
    \pi(x,y)\,f &:=& f(x^*\,\cdot\,y) \\
     f^* (x) &:=& f(x^*)
   \end{array}\;, &\mb{where} & x,y\in G\,,\quad f\in\FK{\Gq}\;\;.
\end{eqnarray*} 
Then the Peter and Weyl theorem of above takes the following form: We get a $\Gq\times\Gq$-equivariant linear bijective map
\begin{eqnarray*}
  \Psi:\;\;\bigoplus_{\La\,\in\, P^+} L(\La)\otimes L(\La) \;\,\to \;\,\FGq
\end{eqnarray*}
by assigning $ v\otimes w $ the function $f_{vw}\res{\overline{G}}=\kB{v}{\;\cdot\;w}\,\res{\overline{G}}$, $v,w\in L(\La)$, $\La\in P^+$. 
It identifies the sum of the switch-maps with the involution.\vspace*{1ex}\\  
{\bf Remark:} Due to the last theorem, our main aims in the next two sections are to determine the closure $\Gq\subseteq\grAdj$, as well 
as the Lie algebra $Lie(\Gq)$ of the weak algebraic monoid $\Gq$.
%
%
%
%
%
%
\newpage\section{The proof of $\Gq=\GD$ and some other theorems\label{sectionGqGD}}
%
%
%
%
In this section we prove that the closure $\Gq\subseteq \grAdj$ coincides with $\GD$. This shows that $\GD$ is the monoid associated to the category 
of admissible $\g$-modules of $\cal O$, and its category of restricted duals.\vspace*{1ex}\\ 
We also determine some other closures. In particular we show $\Tq=\TD$. Furthermore $\Tq$ is a generalized affine toric variety, 
and the action of the Weyl group on $T$ extends to an action on $\Tq$. We show that $U^\pm$ is closed. We show $\Nq=\ND$, which makes 
the analogy of the Weyl monoid $\WeD\cong\ND/T=\Nq/T$ with a Renner monoid even more tight. We also show that the dense unit group $G$ is open in $\Gq$.\\
All these results indicate that $\Gq=\GD$, together with its coordinate ring, is an analogue of a reductive algebraic monoid.\vspace*{1ex}\\
Although $\Gq$ does not act on the flag varieties realized in the projective spaces $\mathbb{P}(L(\La))$, $\La\in P^+$, we show that it acts on the 
corresponding affine cones.\vspace*{1ex}\\ 
Furthermore we show that the Kac-Peterson-Slodowy part of the $\F$-valued points of the algebra of strongly regular functions coincides with the part, 
which is given by $\GD$.\vspace*{1ex}\\
The hearts of the proofs of $\Gq=\GD$ and $\Nq=\ND$ are inductions over the cardinality of $J\subseteq I$,
showing $\overline{G_J}=\GD_J$ and $\overline{N_J}=\ND_J$. The results in the next subsections, describing 
something indexed by $J$, or indexed by ``$rest$'', are needed to prepare these proofs. 
%
%
%
\subsection{The coordinate rings and closures of $T$, $T_J$ ($J\subseteq I$), and $T_{rest}$}
%
%
%
First some notations and remarks. Let $J\subseteq I$. Define the following sublattices of $H$ and $P$:\vspace*{1.5ex}\\
$\begin{array}{ccc}
    H_J &:=& \Z\mb{-span}\Mklz{h_j}{j\in J}\;, \\
    P_J &:=& \Z\mb{-span}\Mklz{\La_j}{j\in J}\;,  
\end{array}$
$\begin{array}{ccc}
    H_{rest} &:=& \Mklz{h\in H}{\La_i(h)=0\;\mb{ for }\;i\in I}\;,  \\
    P_{rest} &:=& \Mklz{\la\in P}{\la(h_i)=0 \;\mb{ for }\;i\in I}\;. 
\end{array}$\vspace*{1.5ex}\\
Define the following projections:\vspace*{1ex}\\
\hspace*{2em} $\begin{array}{cccc}
  p_J: &   P & \to     & \quad P_J                    \\
       & \la & \mapsto & \sum_{j\in J} \la(h_j)\La_j
\end{array}\;\;\;\;$,
$\;\;\;\;\begin{array}{cccc}
  p_{rest}: &  P  & \to     & \;P_{rest} \\
            & \la & \mapsto & \la - p_I(\la)
\end{array}\;\;$.\vspace*{1ex}\\
Recall that the torus  $T$ of the Kac-Moody group can be described by the following two isomorphisms of 
groups:\vspace*{1.5ex}\\
\hspace*{2em} $\begin{array}{ccc}
 H\otimes_\Z\F^\times &\to & \quad\:T \\
 \sum_{\tau}h_\tau\otimes s_\tau &\mapsto & \prod_\tau t_{h_\tau}(s_\tau) 
\end{array}\;\;\;\;$,
$\;\;\;\;\begin{array}{ccc}
\mb{Hom}(P,\F^\times) &\to & \;T \\
     \al   &\mapsto & t(\al)
\end{array}\;\;$,\vspace*{1ex}\\
where $t_h(s)v_\la := s^{\la(h)}v_\la$ and $t(\al)v_\la := \al(\la)v_\la$, $v_\la\in L(\La)_\la$, $\La\in P^+$.\vspace*{1ex}\\
The group algebra $\FK{P}$ of the lattice $P$ can be identified with the classical coordinate ring of the torus $T$,
identifying $\sum c_\la e_\la\in \FK{P}$ with the function
\begin{eqnarray*}
  \left(\,\sum_\la \,c_\la \,e_\la\,\right)(\,t(\al)\,) \;\,:=\;\,\sum_\la \,c_\la \,\al(\la)
  &\quad,\quad & \al\in\mb{Hom}(P,\F^\times)\;\;.
\end{eqnarray*}\vspace*{0.1ex}\\
We have similar descriptions for the torus $T_J$, replacing $T$, $H$, $P$ by $T_J$, $H_J$, $P_J$, where now $t(\al)$ acts as 
$t(\al)v_\la:= \al(p_J(\la))v_\la$, $v_\la\in L(\La)_\la$, $\La\in P^+$. We have similar descriptions for the torus $T_{rest}$.
\vspace*{1ex}\\ 
To determine the closures of $T_J$, $T$, and $T_{rest}$ we need the coordinate rings 
$\FK{T_J}$, $\FK{T}$, and $\FK{T_{rest}}$. Except the coordinate ring of $T_{rest}$, these are in general only 
subalgebras of the classical coordinate rings:
\begin{Prop}\mb{}\label{Gq1}
Let $J\subseteq I\,$. We have: \vspace*{0.5ex}\\
1) $\qquad \FK{T_J}\:\,=\,\:\FK{\,p_J\left(X\cap P\right)}\,\;=\;\,
\FK{\,p_{J^\infty}(X\cap P)+P_{J^0}\,}\;\;$.\vspace*{0.5ex}\\
2) $\qquad \FK{T}\,\:=\,\:\FK{X\cap P}\;\;$.\vspace*{0.5ex}\\
3) $\qquad \FK{T_{rest}}\,\:=\,\:\FK{P_{rest}}\;\;$.
\end{Prop}
\Proof
We only show 1) in the nontrivial case $J\neq\emptyset$. The proofs of 2) and 3) are similar.\\
a) Let $v_\la\in L(\La)_\la$, $w_\mu\in L(\La)_\mu$, $\la,\mu\in P(\La)$, and $\La\in P^+$. By checking 
on the elements of $T_J$, we find 
\begin{eqnarray}\label{TResMat}
    f_{v_\la w_\mu}\res{T_J} &=& \kB{v_\la}{w_\mu} \, e_{p_J (\mu)} \;\,=\;\,f_{v_\la w_\mu}(1)\, e_{p_J (\mu)}\;\;.
\end{eqnarray}
Due to the Peter and Weyl theorem for $\FK{G}$, and 
$\bigcup_{\La\in P^+} P(\La) = X\cap P$, we get $\FK{T_J}\subseteq \FK{p_J(X\cap P)}$.
Due to the nondegeneracy of the forms $\kBl$ on the weight spaces, we have even equality.\vspace*{0.5ex}\\
b) Let $J^0\neq\emptyset\,$. We have 
$p_J(X\cap P) = (p_{J^\infty}+p_{J^0})(X\cap P)\subseteq p_{J^\infty}(X\cap P)+P_{J^0}$.
To show the reverse inclusion, note that the comorphism of the multiplication map 
$ G_{J^\infty}\times G_{J^0} \to G_J$ gives an isomorphism of algebras $\FK{G_J}\to
\FK{G_{J^\infty}}\otimes\FK{G_{J^0}}$. It is easy to conclude, that the comorphism of the multiplication map 
$m:\, T_{J^\infty}\times T_{J^0} \to T_J$ exists, and gives an isomorphism of algebras
\begin{eqnarray*}
  m^*:\;\FK{T_J} &\to & \FK{T_{J^\infty}}\otimes\FK{T_{J^0}}\;.
\end{eqnarray*} 
Due to a) we have $\FK{T_J}=\FK{p_J(X\cap P)} $, and 
$\FK{T_{J^\infty}}=\FK{p_{J^\infty}(X\cap P)}$. Because $G_{J^0}$ is a Kac-Moody group of 
finite type with maximal torus $T_{J^0}$, we have $\FK{T_{J^0}}=\FK{P_{J^0}}$. By checking on the 
elements of $T_J$, we find
\begin{eqnarray*}
  (m^{-1})^*\left(e_{\la_\infty}\otimes e_{\la_0}\right) 
  \;=\;e_{\la_\infty+\la_0} &,& \la_\infty\in p_{J^\infty}(X\cap P)\,,\;\la_0\in P_{J^0}\;\;.
\end{eqnarray*}
Therefore $p_{J^\infty}(X\cap P)+P_{J^0} \subseteq p_J(X\cap P)$.\\
\End
\begin{Theorem}\mb{}\label{Gq2}
Let $J\subseteq I\,$. We have: \vspace*{0.5ex}\\
1) $\quad \overline{T_J}\:=\:\TD_J\;\;$.\vspace*{0.5ex}\\
2) $\quad \Tq\:=\:\TD\;\;$.\vspace*{0.5ex}\\
3) $\quad T_{rest}\,$ is closed.
\end{Theorem}
\Proof
We only show 1) in the nontrivial case $J\neq\emptyset$. The proofs of 2) and 3) are similar to part c) of the following proof.\vspace*{0.5ex}\\
a) To prepare the proof, we show that the submonoid $p_J(X\cap P)$ of $P_J$ is saturated, and its faces are given by
\begin{eqnarray*}
   \Fa{\,p_J(X\cap P)\,} &=& \Mklz{p_J(R\cap P)\,}{\,R\in\RkX_J }\;\;.
\end{eqnarray*}
$\bullet$ We first treat the case $J=J^\infty$. Due to Theorem \ref{GD2} b), the kernel of $p_{J^\infty}$ is the hull subgroup of $R(J^\infty)\cap P$:
\begin{eqnarray}\label{kernelpinfty}
 \Kern\,p_{J^\infty} \;\,= \;\, \bigl(R(J^\infty)-R(J^\infty)\bigr)\cap P\;\;.
\end{eqnarray}
Recall that $(X-R(J^\infty))\cap P$ is a saturated submonoid of $P$, whose faces are given by 
\begin{eqnarray*}
    Fa(\,(X-R(J^\infty))\cap P\,) &=& \Mklz{(R-R(J^\infty))\cap P}{R\in \RkX\,,\,R\supseteq R(J^\infty)}\\
                              &=& \Mklz{(R-R(J^\infty))\cap P}{\,R\in\RkX_J}\;\;.
\end{eqnarray*} 
The surjective linear map $p_{J^\infty}:P\to P_{J^\infty}$ maps the saturated submonoid $(X-R(J^\infty))\cap P$ surjectively onto 
$p_{J^\infty}(X\cap P)$. Its kernel (\ref{kernelpinfty}) is smallest the face of $(X-R(J^\infty))\cap P$.\\
Using these facts, it is easy to check, that $p_{J^\infty}(X\cap P)$ is saturated. It is also not difficult to check, that the map  
\begin{eqnarray*}
    Fa(\,(X-R(J^\infty))\cap P\,) &\to & Fa(\,p_{J^\infty}(X\cap P)\,)\\
                   F\qquad &\mapsto& \quad p_{J^\infty}(F)
\end{eqnarray*}
is bijective. Using the description of the faces of $(X-R(J^\infty))\cap P$, and once more (\ref{kernelpinfty}), we get the required 
result.\\   
$\bullet$ For an arbitrary $J\subseteq I$, we have due to the last proposition:
\begin{eqnarray}\label{ZerpJXcapP} 
  p_J(X\cap P) &=& p_{J^\infty}(X\cap P)+P_{J^0}\;\,\subseteq \;\, P_{J^\infty}\oplus P_{J^0}\;\;.
\end{eqnarray} 
Using the preceeding results, it is easy to see, that $p_J(X\cap P)$ is saturated in $P_J$, and its faces are given by
\begin{eqnarray*}
  \Fa{\,p_J(X\cap P)\,} &=& \{ \, p_{J^\infty}(R\cap P)+P_{J^0}\,\mid \,R\in \RkX_{J^\infty}\,=\,\RkX_J \,\}\;.
\end{eqnarray*}
Because of $p_J = p_{J^\infty} + p_{J^0}$ and $p_{J^\infty} = p_J - p_{J^0}$, we find 
\begin{eqnarray*}
  p_J(R\cap P)\;\,\subseteq\;\, p_{J^\infty}(R\cap P)+P_{J^0} \;\,\subseteq\;\, p_J(R\cap P)+P_{J^0} \;\;.   
\end{eqnarray*}
To show that the first inclusion is an equality, it is sufficient to show $P_{J^0}\subseteq p_J(R\cap P)$. Let $\mu\in P_{J^0}$. Due to 
(\ref{ZerpJXcapP}) there exists an element $\la\in X\cap P$, such that $\mu=p_J(\la)$. 
From this follows $0=\mu(h_i)=p_J(\la)(h_i)=\la(h_i)$ for all $i\in J^\infty$. Therefore $\la\in X\cap ( R(J^\infty)+R(J^\infty) )\cap P=
R(J^\infty)\cap P$, and $\mu\in p_J(R(J^\infty)\cap P)\subseteq p_J(R\cap P)$.\vspace*{1ex}\\
b) Fix $\mu\in X\cap P$ and $R\in \RkX_J$. Next we show
\begin{eqnarray*}
 \mu \,\in\,  R\cap P &\iff& p_J(\mu)\,\in \,p_J(R\cap P)\;\;.
\end{eqnarray*}
The direction '$\Rightarrow$' is trivial. Now let $p_J(\mu)\in p_J(R\cap P)$. By using (\ref{kernelpinfty}) we find
\begin{eqnarray*}
  \mu &\in& \left((R\cap P)+\Kern\, p_J\right)\cap X \;\,\subseteq\; \,\left((R\cap P)+\Kern\,p_{J^\infty}\right)\cap X \\
  &=& \left(R-R(J^\infty)\right)\cap X\cap P\;\;.
\end{eqnarray*}
Suppose there exist elements $r_1\in R$, $r_2\in R(J^\infty)$ with $r_1-r_2\in X\setminus R$. Because $R$
is a face of $X$, we get $r_1\in X\setminus R+r_2 \subseteq X\setminus R$, which is a contradiction. Therefore
we conclude $(R-R(J^\infty))\cap X=R$.\vspace*{1ex}\\
c) Due to Proposition \ref{V4} the Zariski closure of $T_J$ is given by 
\begin{eqnarray*}
  \lefteqn{ \overline{T_J} \;=\; \biggl\{\:\phi\in \grAdj\:\biggl|\biggr.\:\:\exists\; \mb{hom. of alg.}\;\beta_\phi:\FK{T_J}\to \F } 
  \qquad\qquad\;\; \\
  & &\;\;\forall\;v,w\in L(\La),\,\La\in P^+\;:\quad\beta_\phi(f_{vw}\res{T_J}) =\kB{v}{\phi w} \:\biggr\}\;\,,
\end{eqnarray*}
and $\phi\in \overline{T_J}$ is uniquely determined by $\beta_\phi$. Due to the last proposition we have $\FK{T_J}=\FK{p_J(X\cap P)}$.\\ 
For $t\in T_J$ let $\al_t\in\mb{Hom}(P_J,\F^\times)$ be the corresponding homomorphism of monoids. Note that $p_J(X\cap P)$ spans $P_J$, 
because $X\cap P$ spans $P$. Due to a), and the remark following Proposition \ref{Fa4}, we get for every $t\in T_J$ and $R\in\RkX_J$ a 
homomorphisms of algebras $\beta_{t,\,R}:\FK{p_J(X\cap P)}\to \F$ by
\begin{eqnarray*}
  \beta_{t,\,R}\,(e_\la) &:=& \left\{\begin{array}{ccc}
  \al_t(\la) &\mb{if} & \la\,\in\, p_J(R\cap P) \\
      0      &\mb{if} &  \la\,\in\, p_J(X\cap P)\setminus p_J(R\cap P)
  \end{array}\right. \;\;.
\end{eqnarray*}
Using equation (\ref{TResMat}), and b), we find for all $v_\la\in L(\La)_\la$, $w_\mu\in L(\La)_\mu$, 
$\la,\mu\in P(\La)$ and $\La\in P^+$:
\begin{eqnarray*}
  \beta_{t,\,R}(\,f_{v_\la w_\mu}\res{T_J}\,) &=&\left\{\begin{array}{ccc}
   \al_t(p_J(\mu))\kB{v_\la}{w_\mu} &\mb{if} & p_J(\mu)\in p_J(R\cap P) \\
          0          &\mb{else} &     
  \end{array}\right\} \\  
  &=& \kB{v_\la}{te(R)w_\mu}\;\;.
\end{eqnarray*}
Therefore we have $\beta_{t,\,R}\,=\,\beta_{te(R)}$. Also by the remark following Proposition \ref{Fa4}, every homomorphism is of this form. 
Therefore $\overline{T_J}=\TD_J$.\\
\Ende
{\bf Remark:} The proof of the last theorem shows, that $(\,\Tq\,,\,\FK{\Tq}=\FK{X\cap P}\,)$ can be identified with the 
generalized affine toric variety, associated to the saturated submonoid $X\cap P$ of $P$. In particular all statements of Proposition 
\ref{Fa4} are valid.\\ 
Conjugation by elements of $N$ induces an action of the Weyl group $\We\cong N/T$ on $\Tq$, extending the action on $T$. It is given explicitely by 
\begin{eqnarray*}
    \sigma ( t e(R)) &=& \sigma(t) e(\sigma R) \quad,\quad t\in T\;,\; R\in\RkX\;\;.
\end{eqnarray*}  
%
%
%
\subsection{The orbits $\overline{G_J}\,(L(\La)_\La)$ ($J\subseteq I$) and $\Gq \,(L(\La)_\La)$}
%
%
%
Let $\La\in P^+$. The flag variety $G/P_{\{\,i\,\mid\,\La(h_i)=0\,\}}$ can be realized as a subvariety of $\mathbb{P}(L(\La))$. The corresponding 
affine cone, which we call the Kostant cone, is given by the orbit $G (L(\La)_\La)\subseteq L(\La)$. Although $\Gq$ does not act on the flag 
variety, we show that it acts on this cone, i.e., $\Gq (L(\La)_\La) = G (L(\La)_\La) $. This is advantageous, because of the rich combinatorial structure of the flag variety.\\
Furthermore we show $\overline{G_J}(L(\La)_\La)=G_J(L(\La)_\La)$. For this reason we first show the following Proposition: 
\begin{Prop}\mb{}\label{Gq3}
Let $\La\in P^+$ and $\emptyset\neq J\subseteq I$. With the notations of Proposition \ref{GD30-1} there exists a 
$\g(A_J)$-submodule of $L(\La)$, which is a orthogonal complement of $L_J(\La)$. 
\end{Prop}
\Proof\\
{\bf 1)} Let $V$ be a $\g$-module isomorphic to a sum of irreducible highest weight modules $L(N)$, $N\in\h^*$. It
decomposes in a direct sum
\begin{eqnarray*}                                              
    V &=& \bigoplus_{\:N\,\in\,{\bf h}^*}V(N)\;\,,
\end{eqnarray*} 
where $V(N)$ is the $L(N)$-isotypical component. Let $\kBl$ be a nondegenerate contravariant symmetric bilinear 
form on $V$. We show that different isotypical components are orthogonal:\vspace*{0.5ex}\\
Let $v_1\in V(N_1)\,$, $v_2\in V(N_2)\,$ and $N_1\neq N_2$. Because $V(N_1)$ is isomorphic to a sum of 
$L(N_1)$-mo\-du\-les, there exist elements $x\in{\cal U}(\n^-)$, 
$\ti{v}_1\in V(N_1)_{N_1}$ such that $v_1=x\ti{v}_1$. We get
\begin{eqnarray*}    
 \kB{v_1}{v_2} &=& \kB{x\ti{v}_1}{v_2} \;\,=\;\,    
 \langle\langle\,\underbrace{\ti{v}_1}_{\in\, V_{N_1}}\,\mid\,       
 \underbrace{ x^*\, v_2}_{\in\,\bigoplus_{q\,\in\, Q^+_0}V_{N_2-q}}\rangle\rangle \;\;.   
\end{eqnarray*}    
Different weight spaces are orthogonal. If $N_1 > N_2$, or if $N_1$, $N_2$ are incomparable, then $\kB{v_1}{v_2}=0$. 
Using the symmetry of $\kBl$, we conclude $\kB{v_1}{v_2}=0$ for $N_1<N_2$.\vspace*{0.5ex}\\     
{\bf 2)} Recall the notations of Proposition \ref{GD30-1}. Due to the orthogonality of different weight spaces, we have
\begin{eqnarray}\label{orthVC} 
 L(\La) &=& \obot_{\,C\,\in\, P(\La)/\sim}\; V_C \;\;.
\end{eqnarray}
$V_{[\La]}$ is a direct sum of admissible irreducible highest weight modules of $\g(A_J)$, and $L_J(\La)$ is an 
isotypical component. Now the proposition follows from 1).\\   
\End
Now we show:
\begin{Prop}\mb{}\label{Gq5}
Let $\La\in P^+$. Let $J\subseteq I$. We have: \vspace*{0.5ex}\\
a) $\quad\overline{G_J}\,(L(\La)_\La)\:=\:  G_J\, (L(\La)_\La)\;$.\vspace*{0.5ex}\\
b) $\quad\Gq\,(L(\La)_\La)    \:=\:  G\, (L(\La)_\La) \;$.
\end{Prop}
\Proof
The proof of b) is similar to the proof of a). We only show a) in the nontrivial case $J\neq\emptyset$.\\ 
We equip $L_J(\La)$ with the coordinate ring generated by $\bigoplus_{\la\in P(L_J(\La))} (L(\La)_\la)^*$. Then due to 
\cite{KP2} 3A, 3B, and Proposition \ref{GD30-1}, the orbit $G_J (L(\La)_\La)  \subseteq L_J(\La)$ is Zariski closed in $L_J(\La)$.\\
Due to  Proposition \ref{Gq3} there exists an orthogonal decomposition $L(\La)= L_J(\La)\,\obot\,R$. 
Using this decomposition, it is easy to see, that the coordinate ring of $L_J(\La)$ is the restriction of the coordinate ring of $L(\La)$. 
From this decomposition also follows, that $L_J(\La)$ is Zariski closed in $L(\La)$.\\
For every $v_\La\in L(\La)_\La$ the map 
\begin{eqnarray*}
  \Phi:\; gr\mb{-}Adj\left(\bigoplus_{N\in P^+} L(N)\right) &\to & 
   L(\La)\\
   \phi \qquad\quad &\mapsto & \;\phi v_\La
\end{eqnarray*}
is Zariski continuous. Therefore we get 
\begin{eqnarray*}
  \overline{G_J}\, v_\La \;=\; \Phi(\overline{G_J}) \;\subseteq\; \overline{\Phi(G_J)}
  \;\subseteq\;\overline{G_J\,(L(\La)_\La)}\;=\;G_J\,(L(\La)_\La)\;\;.
\end{eqnarray*}
\End
%
%
%
%
%
\subsection{The coordinate rings and closures of $U_J^\pm$, $(U^J)^\pm$ ($J\subseteq I$)}
%
%
%
Kac and Peterson showed in \cite{KP2}, Lemma 4.3, that the coordinate ring $\FK{U^+}$ is a symmetric algebra in
$\Mklz{\, f_{v_\La xv_\La}\res{U^+} }{ x\in\n^- \,}$, where $v_\La\in L(\La)_\La\setminus\{0\}$, $\La\in P^{++}$. From this result follows:
\begin{Cor}\mb{}\label{Gq6}
Let $J\subseteq I$. Let $\La\in F_{I\setminus J}\cap P^+$, and $v_\La\in L(\La)_\La\setminus\{0\}$. Then:\vspace*{0.5ex}\\
\hspace*{1em} $\FK{U_J^+}$ is a symmetric algebra in $\{\,f_{v_\La xv_\La}\res{U_J^+}\,|\,x\in \n_J^-\,\}\,$. 
\vspace*{0.5ex}\\
\hspace*{1em} $\FK{U_J^-}$ is a symmetric algebra in $\{\,f_{xv_\La v_\La}\res{U_J^-}\,|\,x\in \n_J^-\,\}\,$.
\end{Cor}
\Proof 
The case $J=\emptyset$ is trivial. If $J\neq\emptyset$, then
$G_J$, $\FK{G_J}$ can be identified with $G(A_J)'$, $\FK{G(A_J)'}$. We get the first statement by using Proposition 
\ref{GD30-1}, and the result of Kac and Peterson of above.\\ 
The second statement follows from the first by applying the comorphism of $*:U_J^-\to U_J^+$, which is an isomorphism of algebras 
$*:\FK{U_J^+}\to\FK{U_J^-}$.\\
\End\\ 
For $J\subseteq I$ set $(\W^J)^\pm := \W^\pm\setminus\sum_{j\in J}\Z\al_j$, and $(\n^J)^\pm:=
\bigoplus_{\al\in(\Delta^J)^\pm}\g_\al\,$. Write $\n^J:=(\n^J)^+$ for short.
\begin{Theorem}\mb{}\label{Gq7}
Let $J\subseteq I\,$. Let $\La\in F_J\cap P^+$, and $v_\La\in L(\La)_\La\setminus\{0\}$.\vspace*{0.5ex}\\
\hspace*{1em} The algebra $\FK{(U^J)^+}$ is generated by $\{\,f_{v_\La xv_\La}\res{(U^J)^+}\,|\,x\in (\n^J)^-\,\}\,$.\vspace*{0.5ex}\\
\hspace*{1em} The algebra $\FK{(U^J)^-}$ is generated by $\{\,f_{xv_\La v_\La}\res{(U^J)^-}\,|\,x\in (\n^J)^-\,\}\,$.
\end{Theorem}
\Proof
We only show the first statement. The second follows by applying the comorphism of $*:(U^J)^-\to (U^J)^+$, which is 
an isomorphism of algebras $*:\FK{(U^J)^+}\to\FK{(U^J)^-}$.\\
The case $J=I$ is trivial, let $J\neq I$. Main parts of the following proof are similar to the proof of the description of $\FK{U^+}$, 
given by Kac and Peterson. But we also use this result for the proof.\vspace*{0.5ex}\\    
$\bullet$ Let $h\in \h$. First we show
\begin{eqnarray*}
  uh\,-\,h\:\in\: \n^J  \;\mb{ for all }\; u\in U^J\;\;.
\end{eqnarray*}  
It is sufficient to show 
\begin{eqnarray*}
  u_1\exp(x_1)u_1^{-1}\,u_2\exp(x_2)u_2^{-1}\, \cdots\, u_p\exp(x_p)u_p^{-1}
  \,h\,-\,h \;\in\; \n^J
\end{eqnarray*}
for all $\beta_1,\ldots,\, \beta_p\in\W^J\cap\prW$, $x_1\in\g_{\beta_1},\,\ldots,\, x_p\in\g_{\beta_p}$, 
$u_1,\ldots,\, u_p\in U$, $p\in \N$, by induction over $p$.\\
Let $\al\in\pW$, $\beta\in(\W^J)^+$ such that $\al+\beta\in\pW$. The definition of
$(\W^J)^+$ implies $\al+\beta\in (\W^J)^+$. From this follows $U_\al\,\n^J\subseteq\n^J$ for all $\al\in\prW$. 
Because $U$ is generated by the root groups $U_\al$, $\al\in\prW$, we also get 
\begin{eqnarray}\label{UnJnJ}
   U \:\n^J \;\subseteq\;\n^J\;\;.
\end{eqnarray}
The begin of the induction $p=1$ follows by using this inclusion:
\begin{eqnarray*}
  u_1\exp(x_1)u_1^{-1}\, h\,-\,h 
     &=& \sum_{k\,=\,1}^\infty\,\frac{1}{k!}\,[u_1 x_1[\ldots [u_1 x_1,h]\ldots ]] \;\,\in\;\, \n^J\;\;.
\end{eqnarray*}
Now the step of the induction from $p$ to $p+1$:
\begin{eqnarray*}
  \lefteqn{u_1\exp(x_1)u_1^{-1}\,u_2\exp(x_2)u_2^{-1}\,\cdots\,
  u_{p+1}\exp(x_{p+1})u_{p+1}^{-1}\,h\,-\,h} \\
  &=& \underbrace{u_1\exp(x_1)u_1^{-1}\underbrace{\Bigl(u_2\exp(x_2)u_2^{-1}
  \,\cdots \,u_{p+1}\exp(x_{p+1})u_{p+1}^{-1}\,h\,-\,h
  \Bigr)}_{\in\,{\bf n}^J\;due\;to\; the\;induction\;assumption}}_{\in\,
  {\bf n}^J \;due \;to\;(\ref{UnJnJ})}\\ &&\,+\,
  \underbrace{u_1\exp(x_1)u_1^{-1}\,h\,-\,h}_{\in\,{\bf n}^J\;due\;to\;p\,=\,1}  
  \;\,\in\;\, \n^J\;\;.
\end{eqnarray*}
$\bullet$ Let $\FK{\n^J}$ be the algebra of functions generated by
\begin{eqnarray*}
    f_x:=\iB{x}{\;\cdot\;}\res{{\bf n}^J}\;,\quad x\in \n^-_J \;\;.
\end{eqnarray*}
Denote by $\iota: U^J\to U^J$ the inverse map. 
For $N\in P^+$ denote by $\psi_N:U^J\to\n^J$ the map $\psi_N(u):=u\,\nu^{-1}(N)-\nu^{-1}(N)$, $u\in U^J$. Choose $v_N\in L(N)_N\setminus\{0\}$. 
We show
\begin{eqnarray}\label{psfx}
  \psi_N^*(f_x) \;=\; \iota^*\Bigl(\,\frac{1}{\kB{v_N}{v_N}}\,f_{v_N\,xv_N}\res{U^J}\,\Bigr)\;\;, 
\end{eqnarray}
by checking this equality on the elements $u\in U^J$:
\begin{eqnarray*}
   \lefteqn{ f_x(\psi_N(u)) \,\;=\;\, \iB{x}{ u\,\nu^{-1}(N)-\nu^{-1}(N) }
   \,\;=\;\,  \iB{u^{-1}x}{\nu^{-1}(N)} }\\ 
   &=& N\left((u^{-1}x)_0\right)
   \;\,=\;\,  \frac{ \kB{v_N}{ (u^{-1}x)_0 v_N } } { \kB{v_N}{v_N} } 
   \;\,=\;\,  \frac{ \kB{v_N}{ (u^{-1}x) v_N   } } { \kB{v_N}{v_N} } 
   \nonumber\\
   &=&  \frac{ \kB{v_N}{ u^{-1}(x(uv_N)) } } { \kB{v_N}{v_N} } 
   \;\,=\;\,  \frac{ \kB{v_N}{ u^{-1}(x v_N)   } } { \kB{v_N}{v_N} } \;\;.
\end{eqnarray*}
Here $(u^{-1}x)_0$ denotes the $\g_0$-homogeneous part of $u^{-1}x$.  
We used several times the orthogonality of root spaces $\g_\al$ and $\g_\beta$, $\al+\beta\neq 0$, and the orthogonality 
of different weight spaces.\\
Due to \cite{KP2}, Lemma 4.2, the algebra $\FK{U}$ is a Hopf algebra. Therefore also $\FK{U^J}$ is a Hopf algebra. Because of (\ref{psfx}) 
we get $\psi_N^*\left(\FK{\n^J}\right)\subseteq\FK{U^J}$. Because of the description 
of $\FK{U}$ from Kac and Peterson of above, we have equality at least for $N\in P^{++}$.\\
Now fix an element $\ti{\La}\in P^{++}$. We have shown 
\begin{eqnarray*}
  \psi_\La^*\left(\FK{\n^J}\right)\;\subseteq\;\FK{U^J}\;=\;\psi_{\ti{\La}}^*\left(\FK{\n^J}\right)\;\;.
\end{eqnarray*}
The theorem is proved, if we show equality. 
For every $\al\in (\W^J)^+$  choose $\iB{\;}{\;}$-dual bases $(e_\al^{(i)})_{i=1,\,\ldots,\, m_\al}\subseteq \g_\al$,
$(e_{-\al}^{(i)})_{i=1,\,\ldots,\, m_\al}\subseteq\g_{-\al}$. It is sufficient to show
\begin{eqnarray*}
 \psi_{\ti{\La}}^* (f_{e_{-\al}^{(j)}})\;\in\;\psi_\La^*\left(\FK{\n^J}\right) 
  &\mb{ for all }& \al\in(\W^J)^+\,,\;\:j=1,\,\ldots,\, m_\al\;\;.
\end{eqnarray*}
For every $u\in U^J$ write the images $ \psi_{\ti{\La}}(u) \in \n^J$, $\psi_\La(u) \in \n^J$ as linear combinations in the basis 
$(e_\al^{(i)})_{\al\in(\Delta^J)^+\,,\,i=1,\,\ldots,\, m_\al}$ of $\n^J$:
\begin{eqnarray*}
  \psi_{\ti{\La}}(u) &=& \sum_{\al\in(\Delta^J)^+}\sum_{j} \,\ti{\phi}^{(j)}_\al (u)\,e^{(j)}_\al \;\;,\\
  \psi_\La(u) &=& \sum_{\al\in(\Delta^J)^+}\sum_{j} \,\phi^{(j)}_\al (u) \, e^{(j)}_\al \;\;,
\end{eqnarray*}
where
\begin{eqnarray*}
  \ti{\phi}^{(j)}_\al &:=& \iB{e_{-\al}^{(j)}}{\psi_{\ti{\La}}(\;\cdot\;)}\;\,=\;\, \psi_{\ti{\La}}^* (f_{e_{-\al}^{(j)}})   \;\;,\\
  \phi^{(j)}_\al &:=& \iB{e_{-\al}^{(j)}}{\psi_\La(\;\cdot\;)} \;\,=\;\, \psi_{\La}^* (f_{e_{-\al}^{(j)}})\;\in\;\psi_\La^*\left(\FK{\n^J}\right)\;\;.
\end{eqnarray*}
Because of $[\nu^{-1}(\La),\nu^{-1}(\ti{\La})]=0$, we get
\begin{eqnarray*}
  0 &=& u\, [\,\nu^{-1}(\La),\,\nu^{-1}(\ti{\La})] \;\,=\;\, [u\,\nu^{-1}(\La),u\,\nu^{-1}(\ti{\La})] \\ 
    &=& \left[\psi_\La(u),\psi_{\ti{\La}}(u)\right] +[\,\nu^{-1}(\La),\psi_{\ti{\La}}(u)]+
         [\psi_\La(u),\,\nu^{-1}(\ti{\La})]\\ 
    &=& \sum_{\al\,\in\,(\Delta^J)^+}\sum_j \left(-\iB{\al}{\ti{\La}}\phi_\al^{(j)}(u)\,+\,
        \iB{\al}{\La}\ti{\phi}_\al^{(j)}(u)\right) e_{\al}^{(j)} \\
    & &  \;+\; \sum_{\al,\beta\,\in\,(\Delta^J)^+}\sum_{j,k}\,\phi_\al^{(j)}(u)
               \ti{\phi}_\beta^{(k)}(u)\,[e_\al^{(j)},e_\beta^{(k)}]  \;\;.
\end{eqnarray*}
Taking the $\gamma$-homogeneous part, $\gamma\in(\W^J)^+$, we find
\begin{eqnarray}
 0 &=& \sum_j \left(-\iB{\gamma}{\ti{\La}}\phi_\gamma^{(j)}\,+\,
  \iB{\gamma}{\La}\ti{\phi}_\gamma^{(j)}\right) e_{\gamma}^{(j)} \nonumber\\
  && +\;\sum_{\al,\beta\,\in\,(\Delta^J)^+\atop \al+\beta=\gamma}\;\sum_{j,k}\:
  \phi_\al^{(j)}\ti{\phi}_\beta^{(k)}\,[e_\al^{(j)},e_\beta^{(k)}] \;\;.
  \label{Koeffptpe}
\end{eqnarray}
An element $\gamma\in (\W^J)^+$ is of the form $\gamma = \sum_{i\,\in\,I}\, p_i \al_i$ with 
$p_i\in\Nn$, and $p_i\neq 0$ for at least one $i\in I\setminus J$.
Because $\La$ is of the form $\La\,=\,\sum_{i\,\in\, I\setminus J} m_i \La_i\,+\,\La_{rest}$ with 
$m_i\in\N$, $\La_{rest}\in P_{rest}$, we find
\begin{eqnarray*}
  \iB{\gamma}{\La} &=& \sum_{i\in I,\,j\in I\setminus J} p_i\, m_j \iB{\al_i}
  {\La_j}\;+\;0 \;\,=\,\; 
  \sum_{i\,\in\, I\setminus J} p_i\, m_i \,\frac{\iB{\al_i}{\al_i}}{2} \;\,>\;\,0  \;\;.
\end{eqnarray*}
Now we show 
\begin{eqnarray}\label{tphipsiFn}
  \ti{\phi}_\gamma^{(j)}\;\,=\;\, \psi_{\ti{\La}}^* (f_{e_{-\al}^{(j)}}) &\in & \psi_\La^*\left(\FK{\n^J}\right)
\end{eqnarray}
by induction over the height of $\gamma$. To begin the induction, let $\gamma\in(\W^J)^+$ such that 
$ht(\gamma)=\min\Mklz{ht(\delta)}{\delta\in(\W^J)^+}$. Using (\ref{Koeffptpe}) we get
\begin{eqnarray*}
  0\;=\;\sum_j\left(-\iB{\gamma}{\ti{\La}}\phi_\gamma^{(j)}\,+\,\iB{\gamma}{\La}
  \ti{\phi}_\gamma^{(j)}\right) e_\gamma^{(j)} \,+\,0  \;\;.
\end{eqnarray*}
Therefore
\begin{eqnarray*}
  \ti{\phi}_\gamma^{(j)} \;=\; \frac{\iB{\gamma}{\ti{\La}}}{\iB{\gamma}{\La}}\,
  \phi_\gamma^{(j)}\;\in\;\psi_\La^*\left(\FK{\n^J}\right) &\qquad & 
  j=1,\,\ldots,\, m_\gamma  \;\;.
\end{eqnarray*}
Let (\ref{tphipsiFn}) be valid for all $\gamma\in(\W^J)^+$ with $ht(\gamma)<p$.
If there's no $\gamma\in(\W^J)^+$ with $ht(\gamma)\geq p$, there is nothing to show. Otherwise
let $\gamma\in(\W^J)^+$, such that $ht(\gamma)=\min\Mklz{ht(\delta)}{\delta\in(\W^J)^+
\,,\;ht(\delta)\geq p}\,$. Due to (\ref{Koeffptpe}) we have 
\begin{eqnarray*}
 \lefteqn{ \sum_j \iB{\gamma}{\La}\ti{\phi}_\gamma^{(j)} e_\gamma^{(j)}}\\ 
 &=& \sum_j\underbrace{ \iB{\gamma}{\ti{\La}}\phi_\gamma^{(j)}}_{\in\,\psi_\La^*
 ({\mathbb F}[{\bf n}^ J]) }e_\gamma^{(j)}\;-\;
  \sum_{\al,\beta\in(\Delta^ J)^+ \atop \al+\beta=\gamma}\;\sum_{j,k}
  \underbrace{ \phi_\al^{(j)} \ti{\phi}_\beta^{(k)} }_{\in\,
  \psi_\La^*({\mathbb F}[{\bf n}^ J])\;(ind.\;ass.)} \,[e_\al^{(j)},
  e_\beta^{(k)}]\;\;.
\end{eqnarray*}
Because of $\iB{\gamma}{\La}\neq 0$, we conclude $\,\ti{\phi}_\gamma^{(j)}\,\in\,
\psi_\La^*\left(\FK{\n^J}\right)$, $j=1,\,\ldots,\, m_\gamma$.\\
\End
\begin{Theorem}\mb{}\label{Gq8}
Let $J\subseteq I$. We have:\vspace*{0.5ex}\\
1) $\quad U_J^\pm$ is closed.\vspace*{0.5ex}\\
2) $\quad (U^J)^\pm$ is closed.
\end{Theorem}
\Proof
If 1) is valid, then due to
\begin{eqnarray*}
  (U^J)^\pm  &=&  \bigcap_{\sigma\,\in\, {\cal W}_J} \underbrace{\sigma 
  U^\pm\sigma^{-1}}_{closed} 
\end{eqnarray*}
also 2). It is sufficient to  show 1) only in the case of '$-$', because the ho\-meo\-mor\-phism 
$*:\,\mb{gr-Adj}(\ldots)\to\mb{gr-Adj}(\ldots)$ maps $U^-_J$ onto $U^+_J$. Due to Proposition 
\ref{V4}, the Zariski closure of $U_J^-$ is given by 
\begin{eqnarray*}
  \lefteqn{ \overline{U_J^-} \;=\; \biggl\{\;\phi\in\grAdj\:\:\biggl.\biggr|\:\:\exists\;
  \mb{hom. of alg.}\;\beta_\phi:\FK{U_J^-}\to \F } \qquad\qquad \\
  & &\;\;\forall\;v,w\in L(\La),\,\La\in P^+\;:\quad
  \beta_\phi(f_{vw}\res{U_J^-}) =\kB{v}{\phi w} \;\biggr\}\;,
\end{eqnarray*}
and $\phi\in\overline{U_J^-}$ is uniquely determined by $\beta_\phi$.\\
Let $\phi\in \overline{U_J^-}$. Choose an element $v_\La\in L(\La)_\La\setminus\{0\}$, 
$\La\in F_{I\setminus J}\cap P^+$. Due to Proposition \ref{Gq5}, and the Birkhoff decomposition of $G_J$, we have
\begin{eqnarray*} 
  \overline{U_J^-}\F v_\La\;\subseteq\; \overline{G_J}\,\F v_\La \;=\;G_J\,\F v_\La \;=\; U^-_J\,N_J\,\F v_\La \;\;.
\end{eqnarray*} 
Therefore there exist $u\in U_J^-$, $n_\sigma\in N_J$, $c\in\F$, such that
\begin{eqnarray*}
  \phi v_\La &=& c \,u\, n_\sigma\, v_\La\;\;.
\end{eqnarray*}
The coordinate ring $\FK{U_J^-}$ has been described in Corollary \ref{Gq6}. Applying $\beta_\phi$ to the unit 
of $\FK{U_J^-}$, we get
\begin{eqnarray*}
  1 &=& \beta_\phi \left( \frac{f_{v_\La v_\La}\res{U_J^-}}{\kB{v_\La}{v_\La}}\right) 
  \;\,=\;\, \frac{c\kB{u^* v_\La}{n_\sigma v_\La}}{\kB{v_\La}{v_\La}} \;\,=\;\, \frac{c\kB{v_\La}{n_\sigma v_\La}}{\kB{v_\La}{v_\La}} \;\;. 
\end{eqnarray*}
Therefore $c n_\sigma v_\La = v_\La $. Applying $\beta_\phi$ to the generators of 
$\FK{U_J^-}$, we find
\begin{eqnarray*}
   \beta_\phi \left(f_{xv_\La\,v_\La}\res{U_J^-}\right) &=& 
   \beta_{ u }\left(f_{xv_\La\,v_\La}\res{U_J^-}\right)\quad \mb{ for all }\quad x\in\n_J^- \;\;.
\end{eqnarray*}
We conclude $\beta_\phi = \beta_u$, and $\phi= u \in U_J^- $.\\
\End
%
%
%
%
%
\subsection{The closures of $G_J$ ($J\subseteq I$) and $G$}
%
%
%
%
Let $J\subseteq I\,$. The monoid $\overline{G_J}$ contains $G_J$ and $\Tq_J=\TD_J$. Because $G_J$ and $\TD_J$ generate the monoid $\GD_J$, 
we get
\begin{eqnarray}\label{AbschlussGJ}
   \overline{G_J} &\supseteq &  \GD_J\;\;.
\end{eqnarray}
Similarly we get
\begin{eqnarray}\label{AbschlussG}
 \Gq &\supseteq & \GD\;\;.
\end{eqnarray}
In this section we show that we have even $\overline{G_J}=\GD_J$ and $\Gq=\GD\,$.\\\\
%
%
%
%
Kac and Peterson introduced in \cite{KP2}, Section 4D, a multiplicative subset $\gt_\La$, $\La\in P^+_I$, of $\FK{G_I}$. Adopting in our context we 
have:
\begin{Prop}\mb{}\label{Gq9} For $\La\in P^+$ choose an element $v_\La\in L(\La)_\La\setminus\{0\}\,$. Then
\begin{eqnarray*}
  \gt_\La &:=& \frac{1}{\kB{v_\La}{v_\La}}\,\kB{v_\La}
  {\;\cdot\;v_\La}\res{\overline{G}}\;\;\,\in\;\:\FK{\Gq}\setminus\{0\}
\end{eqnarray*}
is independent of the chosen highest weight vector $v_\La$ and the chosen nondegenerate contravariant symmetric bilinear form. The map 
\begin{eqnarray*}
  (\,P^+\,,\,+\:) &\to & (\,\FK{\Gq}\setminus\{0\}\,,\,\cdot\:)\\
  \La \qquad &\mapsto & \quad \gt_\La
\end{eqnarray*}
is an injective homomorphisms of monoids.
\end{Prop}
The proof is essentially the same as in \cite{KP2}. This multiplicative set will be very important for the following considerations. For the 
convenience of the reader we sketch the {\bf proof:}\\
The map $P^+\to \FK{\Gq}\setminus\{0\}$ is injective, because already the restrictions $\gt_\La\,\res{T}\,=e_\La \in\FK{T}$, $\La\in P^+$, are different.\\
Let $\La,N\in P^+\,$. For $u_\pm\in U^\pm\,$, $t\in T\,$ we have
\begin{eqnarray*}
  (\gt_\La\gt_N)(u_-tu_+) &=& (e_\La e_N)(t) \;\,=\;\, e_{\La+N}(t)
  \,\;=\;\,\gt_{\La+N}(u_-tu_+) \;\;.
\end{eqnarray*}
The set $U^- T U^+$ is the principal open set of $G$ associated with $\gt_\La\res{G}$, $\La\in P^+\cap C$. $G$ is irreducible, because the coordinate 
ring $\FK{G}$ is an integral domain. Therefore $U^- T U^+$ is dense in $G$. It is also dense in $\Gq$, and we conclude $\gt_\La\gt_N=\gt_{\La+N}$.\\
\End\\
We will prove $\overline{G_J}=\GD_J$ by induction over $|J|$. To prepare the induction step, two theorems will be important. The first theorem relates 
certain pnc-varieties $D_{\hat{G}_J}(\gt_\La\res{\hat{G}_J})$ of principal open sets to the pnc-varieties $\GD_L$, $L\subsetneqq J$. The 
second theorem gives either a covering of $\overline{G_J}$, or a covering of $\overline{G_J}$ without a point, by certain principal open sets 
$D_{\overline{G_J}}(\gt_\La\res{\overline{G_J}})$. \vspace*{1ex}\\
The following easy proposition states, which of the these pnc-varieties of principal open sets are the same. For $J\subseteq I$ set 
$P_J^+:=P_J\cap P^+=\sum_{j\in J}\Nn\La_j$. Note that 
the faces of $P_J^+$ are given by $P_L^+$, $L\subseteq J\,$. The relative interior of a face 
$P_L$ is given by $ri\,P_L^+ := P_L^{++}:= \sum_{l\in L}\N\La_l \,$, $L\subseteq J$. These relative interiors give a partition of $P_J^+$.
\begin{Prop}
For $J\subseteq I$ let $T_J\subseteq A_J\subseteq \overline{G_J}$. The following map is an injective homomorphisms of monoids:
\begin{eqnarray*}
  (\,P_J^+\,,\,+\:) &\to & (\,\FK{A_J}\setminus\{0\}\,,\,\cdot\:)\\
  \La \qquad &\mapsto & \;\gt_\La\res{A_J}
\end{eqnarray*}
Let $L\subseteq J$ and $\La\in ri\,P_L^+$. The pnc-variety 
\begin{eqnarray*}
  (\,D_{A_J}(\gt_\La\res{A_J})\,,\,\FK{D_{A_J}(\gt_\La\res{A_J})}
\,,\, {\cal F}_{ D_{A_J} (\gt_\La\res{A_J}) } \,) 
\end{eqnarray*}
can be described in the following way: 
\begin{eqnarray*}
    D_{A_J}(\gt_\La\res{A_J}) &=& \Mklz{x\in A_J\,}{\,\gt_{\La'}(x)\neq 0
    \;\mb{ for all }\;\La'\in P^+_{L}} \\
    \FK{D_{A_J}(\gt_\La\res{A_J})} &=& \Mklz{\gt_{\La'}\res{A_J}}{\La'\in P_L^+}^{-1}\FK{A_J} \\
    {\cal F}_{ D_{A_J} (\gt_\La\res{A_J}) }     
    &=& \Mklz{ \Mklz{\gt_{\La'}\res{A_J}}{\La'\in P_L^+}^{-1}R\, }{\, R\in{\cal F}_{A_J} }  
\end{eqnarray*}
\end{Prop}
{\bf Remark:}  We often write only $\gt_\La$ instead of $\gt_\La\res{A_J}\,$ for short. We write 
$(\,D_{A_J}(P^+_L)\,,\,\FK{D_{A_J}(P^+_L)}\,,\,{\cal F}_{ D_{A_J}(P^+_L) }\,)$ 
for this pnc-variety.\vspace*{1ex}\\
\Proof
The homomorphism of monoids is injective, because the restrictions $\gt_\La\res{T_J}=e_\La\in\FK{T_J}$, $\La\in P^+_J$, are different.
\vspace*{1ex}\\
The inclusion $D_{A_J}(\gt_\La)\supseteq D_{A_J}(P^+_L)$ is obvious. To show the reverse inclusion let $x\in D_{A_J}(\gt_\La)$. It is easy 
to check, that 
\begin{eqnarray*}
  \Mklz{\La'\in P^+_J}{\gt_{\La'}(x)\neq 0}
\end{eqnarray*}
is a face of $P^+_J$. Because $P^+_L$ is the smallest face containing $\La$, we have $\gt_{\La'}(x)\neq 0$ for all $\La'\in P^+_L$.\vspace*{1ex}\\
Because of $D_{A_J}(\gt_\La)=D_{A_J}(P^+_L)$ and $\frac{f}{(\gt_\La)^n}=\frac{f}{\gt_{n\La}}$, $f\in\FK{A_J}$, $n\in\N$, we have an inclusion of 
algebras $\FK{D_{A_J}(\gt_\La)} \subseteq \FK{D_{A_J}(P^+_L)}$. To show equality, it is sufficient to show, that for every $\La'\in P^+_L$ there 
exist elements $f\in\FK{A_J}$, $n\in\N$, such that 
\begin{eqnarray*}
     \frac{f}{(\gt_\La)^n} &=& \frac{1}{\gt_{\La'}}\;\;.
\end{eqnarray*}
Because of $\La\in ri\,P^+_L$ and $\La'\in P^+_L$, there exists a positive integer $n\in\N$, such that $n\La-\La'\in P^+_L$. The function 
$f=\gt_{n\La-\La'}$ satisfies the last equation.\vspace*{1ex}\\ 
The natural maps $\FK{A_J}\to\FK{D_{A_J}(\gt_\La)}$ and $\FK{A_J}\to\FK{D_{A_J}(P^+_L)}$ coincide. Therefore the filters of ideals are the 
same.\\
\End
Next we want to describe the pnc-varieties $D_{\hat{G}_J}(P^+_{J\setminus L})$, where $L\subsetneqq J$. For this we first need the following 
proposition:
\begin{Prop}\mb{}\label{Gq11}
Let $K\subseteq I\,$. Then $\left(\,D_{T_K}(P^+_K)\,,\,\FK{D_{T_K}(P^+_K)}\,,\,
{\cal F}_{ D_{T_K}(P^+_K)}\,\right)\,$ is a variety, and we have
\begin{eqnarray*}
  D_{T_K}(P^+_K) &=& T_K\;\;,\\
  \FK{D_{T_K}(P^+_K)} &=& \FK{P_K}\;\;.
\end{eqnarray*}
\end{Prop}
\Proof
Clearly $D_{T_K}(P_K^+)\,=\,T_K\,$. Due to the definition of $\FK{D_{T_K}(P^+_K)}$ and Theorem \ref{Gq1} we get
\begin{eqnarray*}
   \FK{D_{T_K}(P^+_K)} &=& \FK{\,p_K(P\cap X)-P^+_K\,}\;\;.
\end{eqnarray*}
Because of $\,P^+_K\subseteq X\,$ we find
\begin{eqnarray*}
  P_K \;=\; P^+_K - P^+_K \;\subseteq\;
  p_K(P\cap X) - P^+_K \;\subseteq\;P_K\;\;.
\end{eqnarray*}
$\FK{P_K}$ is the classical coordinate ring of the torus $T_K\,$. The map $\,T_K\to \Spm\FK{P_K}\,$ is 
bijective. Therefore also the map $\,T_K\to \FSpm\FK{P_K}\,$ is bijective, and
$\,(\,D_{T_K}(P^+_K)\,,\,\FK{D_{T_K}(P^+_K)}\,,\,{\cal F}_{ D_{T_K}(P^+_K)}\,)\,$ is a variety.\\
\End 
For the next theorem equip the torus $T_{J\setminus L}$ with the variety structure of the last proposition.
\begin{Theorem}\mb{}\label{Gq12}
Let $L\subsetneqq J\subseteq I\,$. Set $\,(U^\pm_J)^L\,:=\,
\bigcap_{\sigma\,\in\,{\cal W}_L}\sigma U^\pm_J\sigma^{-1} \,$. 
We have:\vspace*{0.5ex}\\
a) $\,D_{\hat{G}_J}(P^+_{J\setminus L})\:=\:(U^-_J)^L\,\GD_L\, T_{J\setminus L}\,(U^+_J)^L\,$.\vspace*{1ex}\\
b) The multiplication map
\begin{eqnarray*}
 m:\;\; (U^-_J)^L\times \GD_L\times T_{J\setminus L}\times (U^+_J)^L &\to &
 D_{\hat{G}_J}(P^+_{J\setminus L}) 
\end{eqnarray*}
\hspace*{0.8em} is bijective. Its comorphism 
\begin{eqnarray*}
 m^*:\;\; \FK{D_{\hat{G}_J}(P^+_{J\setminus L})} &\to & 
  \FK{(U^-_J)^L}\otimes\FK{\GD_L}\otimes
  \FK{P_{J\setminus L}}\otimes\FK{(U^+_J)^L}
\end{eqnarray*}
\hspace*{0.8em} exists, and is an isomorphism of algebras.\vspace*{1ex}\\
c) $m^{-1}$ is a morphism of pnc-varieties.
\end{Theorem}
{\bf Remarks:}\\
1) A decomposition analogous to b), $J=I$, $L=\emptyset$, the monoid $\GD_I$ replaced by the 
group $G_I$, has been given by Kac and Peterson in \cite{KP2}, Lemma 4.4.\\
A similar decomposition of the coordinate rings analogous to the second part of b), $J=I$, has been given by 
Kashiwara in \cite{Kas}, Lemma 5.3.4 and Lemma 5.3.5. There coordinate rings are realized as subalgebras 
of the duals of universal enveloping algebras belonging to certain Lie algebras, and Kashiwara works with 
these algebras and its spectra. In particular the Kac-Moody group itself isn't available.\\
2) We do not know if $m^{-1}$ is an isomorphism, or if it is no isomorphism of pnc-varieties.\vspace*{1ex}\\
\Proof\\
{\bf To a):} Let $\La \in ri\,P^+_{J\setminus L}=P_J^+\cap F_L\,$. Due to Theorem \ref{GD33} b) an element
$x\in \hat{G}_J$ can be written in the form $x=u_- n_\sigma e(R)u_+\,$, where $u_{\pm}\in 
U^\pm_J\,$, $n_\sigma\in N_J\,$ and $R\supseteq R(J^\infty)\,$. We have 
\begin{eqnarray*}
 0\;\neq\;\gt_\La(x)\;=\;\frac{\kB{v_\La}{u_- n_\sigma e(R)u_+ 
 v_\La}}{\kB{v_\La}{v_\La}}\;=\; \frac{\kB{v_\La}{n_\sigma e(R)v_\La}} 
 {\kB{v_\La}{v_\La}} \\ 
  \iff   \quad \La\,\in\, R\quad\mb{ and }\quad \sigma\La\,=\,\La  \;\;.\qquad\qquad\qquad
\end{eqnarray*}
Because of $\La\in F_L=F_{L^\infty\cup L^0}\subseteq ri\,R(L^\infty)$ this is equivalent to
\begin{eqnarray*}
   R\supseteq R(L^\infty) \quad\mb{ and }\quad \sigma\in \We_L\;\;.
\end{eqnarray*}
Because of $\We_L\subseteq \We_J$ and $R(L^\infty)\supseteq R(J^\infty)$ we get 
\begin{eqnarray*}
 D_{\hat{G}_J}(P^+_{J\setminus L})&=&D_{\hat{G}_J}(\gt_\La)\\ 
 &=&\Mklz{u_- n_\sigma e(R) u_+}{u_\pm\in U^\pm_J,\,n_\sigma\in \We_L T_J,
 \,R\supseteq R(L^\infty)}\;\;.
\end{eqnarray*}
Using $U^-_J=(U^-_J)^L U_L^-$, $\,U^+_J=U_L^+ (U^+_J)^L$,
and Theorem \ref{GD33} b) once more, now for $\GD_L$, a) follows.\vspace*{1ex}\\ 
{\bf To b) and c):} Both parts follow directly from Proposition \ref{V2} 1), we only have to check 
that the prerequisites are satisfied:\vspace*{0.5ex}\\
$\bullet $ 
Clearly $T_{J\setminus L} = D_{T_{J\setminus L}}(P^+_{J\setminus L})$ is a sub-pnc-variety of 
$D_{\hat{G}_J}(P^+_{J\setminus L})\,$.\\ 
For all $\La\in P^+_{J\setminus L}$ we have $\gt_\La\res{(U_J^\pm)^L}\, = 1_{(U_J^\pm)^L}$. Due to Proposition \ref{GD37} b) 
we also have $\gt_\La\res{\hat{G}_L} \,= 1_{\hat{G}_L}$ . From this follows, that $(U_J^\pm)^L$ and $\GD_L$ are sub-pnc-varieties of 
$D_{\hat{G}_J}(P^+_{J\setminus L})\,$.\vspace*{0.5ex}\\
$\bullet$ The multiplication map is surjective due to part a) of this proof. The remaining demand 
on $m$ is satisfied, because $(U_J^\pm)^L\,$, $\hat{G}_L$ and $T_{J\setminus L}$ contain the unit.\vspace*{0.5ex}\\
$\bullet$ To show that the comorphism $m^*$ exists, and is surjective, first note:
Due to (\ref{AbschlussGJ}) we have $G_K\subseteq \GD_K\subseteq\overline{G_K}$, $K\subseteq I$. Therefore
$\FK{\GD_K}$ is isomorphic to $\FK{G_K}$ by the restriction map. Due to the Peter and Weyl
theorem for $\FK{G_K}$, and Proposition \ref{GD30-1} we conclude
\begin{eqnarray}\label{PWfGDK}
  \FK{\GD_K} &=& \mb{span}\Mklz{ f_{vw}\res{ \hat{G}_K}  }
  { v,w\in L_K(\La),\;\La\in P^+_K }\;\;.
\end{eqnarray} 
Now let $\La\in P^+_{J\setminus L}$. For $N\in P^+ $ we choose $\kBl$-dual bases of $L(N)$, by choosing
$\kBl$-dual bases
\begin{eqnarray*}
(a_{\la\,k})_{k=1,\ldots\, m_{\la}}  &\quad,\quad&
(b_{\la\,k})_{k=1,\ldots\, m_{\la}}
\end{eqnarray*}
of $L(N)_\la$ for every $\la\in P(N)$. Let $v\in L(N)_\la\,$, $w\in L(N)_\mu\,$, 
$\la,\mu\in P(N)\,$. 
By checking on the elements of $D_{\hat{G}_J}(P^+_{J\setminus L})\,$, using Proposition \ref{GD37} b), we find 
\begin{eqnarray*}
  \lefteqn{ m^*\left( \frac{f_{vw}}{\gt_\La}
  \res{D_{\hat{G}_J}(P^+_{J\setminus L})}\right)}\\ 
  &=& \sum_{\mb{$\la'\geq\la\;,\;i 
            \atop \mu'\geq\mu \;,\;j$}}
  f_{v\,a_{\la'i}}\res{(U^-_J)^L}\,\otimes\, 
  f_{b_{\la'i}\,a_{\mu'j}}\res{\hat{G}_L}\,\otimes\, 
  e_{p_{J\setminus L}(\mu')-\La} \,\otimes\,
  f_{b_{\mu'j}\,w}\res{(U^+_J)^L}\;\;.
\end{eqnarray*}
Since $P(N)\subset N - Q^+_0$ this sum is really finite. Using (\ref{PWfGDK}) for $K=J$ we find that 
$ m^* :\FK{D_{\hat{G}_J}(P^+_{J\setminus L})}\to \FK{(U^-_J)^L}\otimes\FK{\GD_L}\otimes
\FK{P_{J\setminus L}} \otimes \FK{(U^+_J)^L}$ exists. To show the surjectivity of $m^*$, it is sufficient to find 
elements of $\FK{D_{\hat{G}_J}(P^+_{J\setminus L})}$, which are mapped by $m^*$ onto a system of generators of
the algebra $\FK{(U^-_J)^L}\otimes\FK{\GD_L}\otimes\FK{P_{J\setminus L}}\otimes\FK{(U^+_J)^L}\,$:\\
By using Proposition \ref{GD37} b) it is easy to see, that for all $\La,N\in P^+_{J\setminus L}$ we have
\begin{eqnarray*}
 m^*\left(\frac{\gt_N}{\gt_\La}\res{D_{\hat{G}_J}(P^+_{J\setminus L})}\right)\;\,=\;\, 
  1\,\otimes\,1\,\otimes\,e_{N-\La}\,\otimes\,1\;\;.
\end{eqnarray*}
Using the same proposition we find for $\La \in ri\,P^+_{J\setminus L}= P_{J\setminus L}^{++} =
P^+_J\cap F_L\,$, $v_\La\in L(\La)_\La\setminus\{0\}\,$, $x\in(\n^-_J)^L\,$:
\begin{eqnarray*}
  m^*\left(\frac{f_{v_\La\,x v_\La}}{\gt_\La}\res{D_{\hat{G}_J}(P^+_{J\setminus L})} \right) &=& 
  1\,\otimes\,1\,\otimes\,1\,\otimes\,f_{v_\La\,xv_\La}\res{(U_J^+)^L}    \;\;,  \\ 
  m^*\left(\frac{f_{x v_\La\, v_\La}}{\gt_\La}\res{D_{\hat{G}_J}(P^+_{J\setminus L})} \right)&=& 
  f_{xv_\La\,v_\La}\res{(U_J^-)^L}\,\otimes\, 1\,\otimes\,1\,\otimes\,1   \;\;.
\end{eqnarray*}
The elements $f_{v_\La\,x v_\La}\res{(U^+_J)^L}$, $x\in(\n^-_J)^L$ generate $\FK{(U^+_J)^L}$, and the elements 
$f_{xv_\La\, v_\La}\res{(U^-_J)^L}$, $x\in(\n^-_J)^L$ generate $\FK{(U^-_J)^L}$.  By identifying $G_J$, $\FK{G_J}$ with 
$G(A_J)'$, $\FK{G(A_J)'}$, this follows from Theorem \ref{Gq7} and Proposition \ref{GD30-1}.\\
Let $N\in P_L^+\,$ and $v_\la\in L_L(N)_\la$, $w_\mu\in L_L(N)_\mu$ where $\la,\,\mu\in P(L_L(N))$.
Because the weight $\mu$ is of the form $\mu=N-\sum_{l\in L}k_l\al_l$ with $k_l\in\Nn$, we have
\begin{eqnarray*} 
   \mu(h_i) \;\,=\;\, \underbrace{N(h_i)}_{\geq\,0}\,-\,\sum_{l\in L}k_l
   \underbrace{\al_l(h_i)}_{\leq\,0}\;\,\geq\;\,0   &\mb{ for all }& i\in J\setminus L\;\;.
\end{eqnarray*}
Set $\La := \sum_{i\in J\setminus L}\mu(h_i)\La_i \in P_{J\setminus L}^+$. By using Proposition \ref{Gq4} we find
\begin{eqnarray*}
 m^*\left(\frac{f_{v_\la w_\mu}}{\gt_\La}\res{D_{\hat{G}_J}(P^+_{J\setminus L})}\right)&=&
    1\,\otimes\,f_{v_\la w_\mu}\res{\hat{G}_L}\,\otimes\,1\,\otimes\,1 \;\;. 
\end{eqnarray*}
Due to (\ref{PWfGDK}) the elements $f_{v_\la w_\mu}\res{\hat{G}_L}$ generate $\FK{\GD_L}\,$.\\
\End
We use the principal open sets $D_{\overline{G_J}}(P^+_{J\setminus L})$, $L\subsetneqq J$ to give a 
covering either of $\overline{G_J}$, or of $\overline{G_J}$ without one point. We first give a theorem, which is used to 
distinguish these two cases, and to describe the vanishing ideal of this point.\vspace*{1ex}\\ 
Let $\emptyset\neq J\subseteq I$. $\FK{\overline{G_J}}$ is isomorphic to $\FK{G_J}$ by the restriction map. Using the 
Peter and Weyl theorem for $\FK{G_J}$ and Proposition \ref{GD30-1}, we get an $\overline{G_J}\times 
\overline{G_J}$-equivariant linear bijective map:
\begin{eqnarray*}
  \overline{\Phi}_J:\;\;\bigoplus_{\La\in P^+_J} L_J(\La)\otimes L_J(\La) 
  &\to & \FK{\overline{G_J}} \\
  v\otimes w\qquad &\mapsto & \;f_{vw}\res{\overline{G_J}}
\end{eqnarray*}
\begin{Theorem} Set $\,{\cal N}_J\::=\:\overline{\Phi}_J\left(\bigoplus_{\La\,\in \,P^+_J\setminus
\{0\}} L_J(\La)\otimes L_J(\La) \right)\,$. Then we have:\vspace*{1ex}\\
\hspace*{1.25em} $\;{\cal N}_J\:=\:I\Bigl(e(R(J))\Bigr)$ if $J\,=\,J^\infty\,$.\vspace*{1ex}\\
\hspace*{1.25em} $\;{\cal N}_J$ is no ideal if $J\,\neq\,J^\infty\,$.   
\end{Theorem}
\Proof\\
a) If $J=J^\infty$ we have to show
\begin{eqnarray*}
  e(R(J))v &=& \left\{ \begin{array}{ccc}
    v  & \mb{ if } & v\,\in\, L_J(0) \\
    0  & \mb{ if } & v\,\in\, L_J(\La)\,,\;\La\in P_J^+\setminus\{0\} 
  \end{array}\right.\;\;.
\end{eqnarray*}
We have $P(L_J(0))= \{0\}\subseteq R(J)$. It remains to show $P(L_J(\La))\cap R(J)=\emptyset$ 
for all $\La\in P_J^+\setminus\{0\}\,$.\\
Suppose there exist elements $\La\in P_J^+\setminus\{0\}\,$, $q=\sum_{i\in J}m_i\al_i\in (Q_J)^+_0\,$, $c\in R(J)$ such that
$\La - q=c\,$. Then we get 
\begin{eqnarray*}
  0 \;\leq\; \La(h_j) \;=\; \sum_{i\in J} m_i\,\al_i(h_j)\,+\,0
  \;=\; \sum_{i\in J} a_{ji}\, m_i &\mb{ for all }& j\in J\;\;.
\end{eqnarray*}
The generalized Cartan matrix $A_J$ has no component of finite type. Due to \cite{K2}, Theorem 4.3, we get
\begin{eqnarray*}
  \La(h_j) \;=\; \sum_{i\in J} a_{ji}\, m_i \;=\; 0 &\mb{ for all }& j\in J\;\;,
\end{eqnarray*}
which contradicts $\La\neq 0\,$.\vspace*{1ex}\\
b) Let $J=J^0$ and suppose ${\cal N}_J$ is an ideal of $\FK{\overline{G_J}}$.
The algebra $\FK{\overline{G_J}}$ is isomorphic to $\FK{G_J}$ by the restriction map. The map $G_J\to\Spm\FK{G_J}$ is bijective, 
because $A_J$ has only components of finite type. Therefore there exists a $g\in G_J$, such that ${\cal N}_J=\Vi{g}\,$.
Choose an element $v\in L_J(\La)\setminus \{0\}\,$, $\La\in P^+_J\setminus \{0\}$. Then
\begin{eqnarray*} 
v\;=\;g^{-1}\underbrace{(gv)}_{=\,0}\;=\;0
\end{eqnarray*}
is a contradiction.\vspace*{1ex}\\
c) Let $J=J^0\cup J^\infty$ with $J^0\neq \emptyset$ and $J^\infty\neq\emptyset$. The generalized Cartan matrix $A_J$ decomposes into the sum 
$A_J=A_{J^0}\oplus A_{J^\infty}$. Recall that we have $G_J\,\cong\,G_{J^0}\times G_{J^\infty}$, and $\FK{G_J}\,\cong\,\FK{G_{J^0}}\otimes 
\FK{G_{J^\infty}}$. The subalgebra corresponding to $\FK{G_{J^0}}\otimes\,1\,$ is given by the image of the embedding  
\begin{eqnarray*}
  i:\,\FK{G_{J^0}}\,\to\,\FK{G_J}\;\;,
\end{eqnarray*}
defined by
\begin{eqnarray*}
  i(f_{vw}\res{G_{J^0}})\;:=\;f_{vw}\res{G_{J}}\quad,\quad v,w\in  L_{J^0}(\La)=L_{J}(\La)\,,\;\La\in P^+_{J^0} \;\;.
\end{eqnarray*}
Suppose ${\cal N}_J$ is an ideal and let $\phi:\FK{G_J}\to \F$ be the corresponding 
homomorphism of algebras with kernel ${\cal N}_J\res{G_J}\,$. Then $\phi\circ i:\FK{G_{J^0}}\to\F$ 
is a homomorphism of algebras with kernel ${\cal N}_{J^0}\res{G_{J^0}}\,$, which contradicts b).\\
\End\\
Let $\sigma,\,\tau\in\We_J$, and choose elements $n_\sigma,n_\tau\in N_J$ belonging to $\sigma, \tau$. Let $K\subseteq I$. Since 
$T_J \,D_{\overline{G_J}}(P^+_K)\,T_J=D_{\overline{G_J}}(P^+_K)$, the set
\begin{eqnarray*}
  \sigma D_{\overline{G_J}}(P^+_{K})\tau  &:=&  n_\sigma D_{\overline{G_J}}(P^+_{K})n_\tau
\end{eqnarray*}
is independent of the chosen elements $n_\sigma, n_\tau$.
\begin{Theorem}\mb{}\label{Gq10}
Let $\emptyset\neq J\subseteq I$. We have:
\begin{eqnarray*}
  \bigcup_{\sigma,\tau\,\in\, {\cal W}_J}\, \bigcup_{\emptyset\,\neq \,K\,\subseteq\, J} 
  \sigma D_{\overline{G_J}}(P^+_{K})\tau &=&
   \left\{\begin{array}{ccc}
           \overline{G_J}             & \mb{ if } & J\neq J^\infty \\
  \overline{G_J}\setminus\{ e(R(J))\} & \mb{ if } & J=J^\infty
  \end{array}\right.
\end{eqnarray*}
\end{Theorem}
{\bf Remarks:}\\
1) For $J=J^\infty$ the element $e(R(J))$ is the zero of $\overline{G_J}$. This follows from Proposition \ref{GD34}, and $\overline{G_J}=\GD_J$, 
which we will prove in the next theorem.\\
2) It is easy to see that it is sufficient to take the second union only over the sets $K=\{j\}$, $j\in J$.\\
In general it is not possible to omit the union over the elements of the Weyl group on the left or on the right. This can be checked by using the 
generalized Cartan matrix of Example 3) of Subsection \ref{FacesoftheTitscone}. Choose $J=I$. Not all idempotents of $E$ would be contained 
in these unions.\\
3) A similar decomposition for the full maximal spectrum of $\FK{G}$ has been given by Kashiwara in 
\cite{Kas}, Proposition 6.3.1.\vspace*{1ex}\\ 
\Proof
For every $\La\in P^+_J\setminus\{0\}$ choose an element $v_\La\in L(\La)_\La\setminus\{0\}$. \\
We have $\,D_{\overline{G_J}}(P^+_{K})\,=\,D_{\overline{G_J}}(\gt_\La)$, $\La\in ri\, P^+_{K}$. Because of 
$\bigcup_{\emptyset\neq K\subseteq J} ri\,P_{K}^+ = P_J^{+}\setminus\{0\}$ we find:
\begin{eqnarray*}
  \lefteqn{ \overline{G_J}\setminus \left(\bigcup_{\emptyset\,\neq\,K\,\subseteq\, J}\, \bigcup_{\sigma,\tau\,
  \in\,{\cal W}_J} \sigma D_{\overline{G_J}}(P^+_{K})\tau \right)}\\ 
  &=&\Mklz{x\in \overline{G_J}}{\forall\,\sigma,\tau\in\We_J\;\;\forall \La\in P_J^+
 \setminus\{0\}\;:\quad \kB{n_\sigma v_\La}{xn_\tau v_\La}\,=\,0 }
\end{eqnarray*}
If $J=J^\infty$, then due to the last theorem $e(R(J))$ belongs to this set. Now let $\emptyset\neq J\subseteq I$ be arbitrary, and 
$x\in \overline{G_J}$ an element of this set, i.e., 
\begin{eqnarray}\label{stWJLPpJkBnull}
   \kB{n_\sigma v_\La}{xn_\tau v_\La}\;\,=\;\,0 &\mb{ for all }& \sigma,\tau\in\We_J,\;\La\in P_J^+ \setminus\{0\}\;\;.
\end{eqnarray}
Fix $\La\in P_J^+\setminus\{0\}$. For an element $v\in L(\La)$ denote by $v_\la$ its component in the weight space $L(\La)_\la$, and set 
$supp(v):=\Mklz{\la}{v_\la\neq 0}$. Denote by $co(supp(v))$ the convex hull of $supp(v)$ in $\h_\R^*$.\\
Because $\overline{G_J}$ is a monoid containing $G_J$, we have $x n_\tau\in\overline{G_J}$. Due to Proposition \ref{Gq5} we find 
\begin{eqnarray*}
  x n_\tau v_\La \;\in\; \overline{G_J} v_\La \;\subseteq \;G_J L(\La)_\La \;\subseteq \;L_J(\La) &\mb{ for all }& \tau\in\We_J\;\;.
\end{eqnarray*}
Using \cite{KP1}, Lemma 4, the vertices of the convex hull $co(supp(xn_\tau v_\La))$ are elements of $\We_J\La$. The weight spaces 
$L(\La)_{\sigma\La}$, $\sigma\in\We_J$, are one-dimensional, and different weight spaces are orthogonal. Because of 
(\ref{stWJLPpJkBnull}) there are no vertices of $co(supp(xn_\tau v_\La))$. Therefore we have 
\begin{eqnarray*}
     0\;=\;xn_\tau v_\La &\mb{ for all }& \tau\in \We_J\;\;.
\end{eqnarray*}
From this follows
\begin{eqnarray*}
  0\;=\;\kB{gv_\La}{xn_\tau v_\La}\;=\;\kB{x^*gv_\La}{n_\tau v_\La} &\mb{ for all }& g\,\in\, G_J,\;\tau\in\We_J\;\;.
\end{eqnarray*}
The involution $*$ is Zariski continuous, leaving $G_J$ invariant. Therefore also $\overline{G_J}$ is invariant, and we get 
$x^*\in \overline{G_J}$. Now we can use the same argument as above to show $x^* g v_\La=0$, $g\in G_J$. Since $L_J(\La)$ is 
spanned by $G_J v_\La$, we find
\begin{eqnarray*}
  x^*L_J(\La) &=&\{0\}\;\;.
\end{eqnarray*}
Because this equation is valid for all $\La\in P_J^+\setminus\{0\}$, we get ${\cal N}_J\subseteq \Vi{x^*}$. We have even 
equality because both are 1-codimensional subspaces of $\FK{\overline{G_J}}$.\\ 
Due to the last theorem we conclude $J=J^\infty$, and $\Vi{\,x^*\,}=\Vi{\,e(R(J))\,}$. Because $\FK{\overline{G_J}}$ separates the 
points of $\overline{G_J}$, we get $x^*=e(R(J))$. Therefore $x=e(R(J))^*=e(R(J))$.\\
\End
%
%
%
%
%
%
%
%
%
%
%
The following theorem is one of the main results of this paper:
\begin{Theorem}\mb{}\label{Gq13}
Let $J\subseteq I$. We have:\vspace*{0.5ex}\\ 
1) $\quad \overline{G_J}\,=\,\hat{G}_J\;$.\vspace*{0.5ex}\\
2) $\quad \Gq=\GD\;$.
\end{Theorem}
\Proof
At the beginning of this subsection we have already proved the inclusions '$\supseteq$'. We only have to show the reverse 
inclusions.\vspace*{0.5ex}\\
{\bf To 1):} We show $\overline{G_J}\subseteq \hat{G}_J$ by induction over $|J|$:\vspace*{0.5ex}\\
The case $J =\emptyset$ is trivial. Let $|J|=1 $. Since $J^\infty = \emptyset$ we have 
$\GD_J = G_J\,$. The map $G_J\to\Spm\FK{G_J}$ is bijective, because ($G_J$, $\FK{G_J}$) can be identified with 
($SL(2,\F)$, $\FK{SL(2,\F)}$). Therefore also the map $G_J\to\FSpm\FK{G_J}$ is bijective, and $G_J$ is Zariski 
closed.\vspace*{0.5ex}\\
Now the step of the induction from $|J|\leq m$ to $|J|=m+1\,$, $(1<m<|I|)\,$:\\
Let $L\subsetneqq J$. Due to Proposition \ref{Gq11} $T_{J\setminus L}$ is a variety. 
Using the induction assumption we find that $\GD_L$ is Zariski closed. Therefore $\GD_L$ is a variety. Left and right translations 
with elements of $G$ are Zariski homeomorphisms. Using Theorem \ref{Gq8} we see that the groups 
\begin{eqnarray*}
   (U_J^\pm)^L &=& \bigcap_{\sigma\,\in\,{\cal W}_L}
   \underbrace{\sigma U^\pm_J \sigma^{-1}}_{closed}\;.
\end{eqnarray*}
are Zariski closed. Therefore also $(U_J^\pm)^L$ are varieties.\\
Due to Theorem \ref{Gq12} there is a bijective morphism
\begin{eqnarray*}
 m^{-1}:\;\;  D_{\hat{G}_J}(P^+_{J\setminus L}) &\to &   (U^-_J)^L\times \GD_L\times T_{J\setminus L}\times (U^+_J)^L \;\;,
\end{eqnarray*}
whose comorphism is also bijective. Applying Proposition \ref{V2} 2), we find that $D_{\hat{G}_J}(P^+_{J\setminus L})$ is a variety. 
We have $\GD_J\subseteq \overline{G_J}$. Applying Proposition \ref{V3} we get
$D_{\overline{G_J}}(P^+_{J\setminus L})= D_{\hat{G}_J}(P^+_{J\setminus L})$. Therefore the principal open set 
$D_{\overline{G_J}}(P^+_{J\setminus L})$ is contained in $\GD_J$. Using Theorem \ref{Gq10} we find
\begin{eqnarray*}
 \overline{G_J} \;\,=\;\, \bigcup_{\sigma,\tau\,\in\,{\cal W}_J}\,
   \bigcup_{L\,\subsetneqq\, J} \underbrace{ \sigma D_{\overline{G_J}}
   (P^+_{J\setminus L}) \tau }_{\subseteq\,\hat{G}_J}
   \,\;\cup\;\,\left\{\begin{array}{ccc}
   \emptyset &\mb{ if } & J\neq J^\infty \\
     e(R(J)) &\mb{ if } & J=J^\infty \end{array}\right\} \;\, \subseteq \;\,\hat{G}_J\;\;.
\end{eqnarray*}
{\bf To 2):} It is easy to check that Proposition \ref{V2}, 1) can be applied to the multiplication map 
$m:\hat{G}_I\times T_{rest}\to\hat{G}$. Therefore $m^{-1}$ is a bijective morphism, whose comorphism is also bijective. Due to part 1) 
$\hat{G}_I$ is Zariski closed. Due to Theorem \ref{Gq2} $T_{rest}$ is Zariski closed. Therefore $\hat{G}_I$ and $T_{rest}$ are varieties. 
Applying Proposition \ref{V2}, 2) we find that $\hat{G}$ is a variety. Therefore $\hat{G}$ is Zariski closed. We conclude $\Gq\subseteq \hat{G}$.\\
\End
%
%
%
%
\subsection{The closures of $N_J$ ($J\subseteq I$) and $N$}
%
%
%
%
The Weyl monoid $\WeD\cong \ND/T$ plays the same role for the Bruhat and Birkhoff decompositions, as the Renner monoid does for the 
Bruhat and Birkhoff decompositions of a reductive algebraic monoid. To make this analogy even closer, we show $\Nq=\ND$.\vspace*{1ex}\\
The proof of $\Nq=\ND$ can be looked at as a simplified version of the proof of $\Gq=\GD$. We only state the main steps, some more
details can be found in \cite{M}.\\\\
Using the same methods as in the proof of Theorem \ref{Gq12} we can show:
\begin{Theorem}\mb{}\label{Gq15}
Let $L\subsetneqq J\subseteq I\,$. We have:\vspace*{0.5ex}\\
a) $\,D_{\hat{N}_J}(P^+_{J\setminus L})\:=\:\ND_L\, T_{J\setminus L}$\vspace*{1ex}\\
b) The multiplication map 
\begin{eqnarray*}
  m:\;\; \ND_L\times T_{J\setminus L} &\to & D_{\hat{N}_J}(P^+_{J\setminus L})  
\end{eqnarray*}
\hspace*{0.8em} is bijective, and its comorphism 
\begin{eqnarray*}
 m^*:\;\; \FK{D_{\hat{N}_J}(P^+_{J\setminus L})} &\to & \FK{\ND_L}\otimes\FK{P_{J\setminus L}}
\end{eqnarray*}
\hspace*{0.8em} exists and is an isomorphism of algebras.\vspace*{1ex}\\
c) $m^{-1}$ is a morphism of pnc-varieties.
\end{Theorem}
We get by intersecting the equation of Theorem \ref{Gq10} with $\overline{N_J}\,$:
\begin{Cor}\mb{}\label{Gq14}
Let $\emptyset\neq J\subseteq I\,$. Then
\begin{eqnarray*}
  \bigcup_{\sigma,\,\tau\,\in\, {\cal W}_J}\, \bigcup_{\emptyset\,\neq \,K\,\subseteq\, J} 
  \,\sigma D_{\overline{N_J}}(P^+_{K})\tau &=&
   \left\{\begin{array}{ccc}
           \overline{N_J}             & \mb{ if } & J\neq J^\infty \\
  \overline{N_J}\setminus\{ e(R(J))\} & \mb{ if } & J=J^\infty
  \end{array}\right.
\end{eqnarray*}
\end{Cor}
Using this theorem and this corollary, we can show in the same way as Theorem \ref{Gq13}:
\begin{Theorem}\mb{}\label{Gq16}
Let $J\subseteq I$. We have:\vspace*{0.5ex}\\ 
1) $\quad \overline{N_J}\,=\,\ND_J$\vspace*{0.5ex}\\
2) $\quad \Nq=\ND$
\end{Theorem}
%
%
%
%
%
%
\subsection{The openness of the unit group}
%
%
%
The unit group of a reductive algebraic monoid is principal open. Here we have:
\begin{Prop}
The unit group $G$ is open in $\Gq$.
\end{Prop}
\Proof
We have $\Gq=\GD$. In the same way as in the proof of Theorem \ref{Gq12} a) it is possible to show
\begin{eqnarray*}
  D_{\overline{G}}(\gt_\La)\;=\;U^- T U^+ &\mb{for every}&\La\in P^{++} \;\;.
\end{eqnarray*}
Due to \cite{KP2}, Corollary 3.1, we have $G_I = \bigcup_{\sigma\,\in\,{\cal W}} \sigma U^- T_I U^+$. Note also that left 
translations with elements of $G$ are Zariski-homeomorphisms. Therefore
\begin{eqnarray*}
 G &=& G_I T_{rest} \;\,=\;\, \bigcup_{\sigma\,\in\,{\cal W}} 
 \underbrace{\sigma U^- T U^+}_{open} 
\end{eqnarray*}
is open in $\Gq\,$.\\
\End
%
%
%
%
%
%
%
\subsection{The Kac-Peterson-Slodowy part of $\Spm\FK{G}$}
%
%
%
Using the characterization $\GD=\Gq$, it is easy to prove the claim of Peterson about the Kac-Peterson-Slodowy part of the $\F$-valued points of 
the algebra of strongly regular functions.\vspace*{1ex}\\
A function $f\in\FK{\GD}$ induces a function on $\Spm\FK{\GD}$, assigning $x\in\Spm\FK{\GD}$ the value $\ti{x}(f)$, where $\ti{x}:\FK{\GD}\to\F$ is 
the homomorphism of algebras with kernel $x$.\\
In this way, $\Spm\FK{\GD}$ is equipped with a coordinate ring isomorphic to $\FK{\GD}$. Its Zariski topology coincides with the relative topology 
of the spectrum of $\FK{\GD}$.\\
For a set $M\subseteq \Spm\FK{\GD}$ we denote by $\overline{M}^{Specm}$ its Zariski closure.\vspace*{1ex}\\
We identify the elements of $\GD$ with the corresponding points of $\Spm\FK{\GD}$ given by the vanishing ideals. Note that for a set 
$M\subseteq \GD=\Gq$ we have 
\begin{eqnarray}\label{clSpecm}
  \overline{M} &\subseteq & \overline{M}^{Specm}\;\;\;.
\end{eqnarray}
\begin{Theorem} We have
\begin{eqnarray*}
  \GD &=& \bigcup_{m\in\N } \;\bigcup_{\beta_1,\,\ldots,\, \beta_m\in\Delta_{re}}  \overline{ U_{\beta_1}\,\cdots\,U_{\beta_m} T}^{Specm}\;\;\;.
\end{eqnarray*}
\end{Theorem}
\Proof
We first prove the inclusion '$\subseteq$': We have $\GD=G'\TD G'$, and $\TD=\Tq$. The derived Kac-Moody group $G'$ is generated 
by the root groups $U_\al$, $\al\in\rW$. Therefore, for every element $\hat{g}\in \GD$, there exist roots 
$\beta_1,\,\ldots ,\,\beta_p,\,\gamma_1,\,\ldots,\,\gamma_q\in\rW$, such that 
\begin{eqnarray*}
 \hat{g} &\in & U_{\beta_1}\,\cdots\,U_{\beta_p}\Tq U_{\gamma_1}\,\cdots\,U_{\gamma_q}\;\;.
\end{eqnarray*}
Because left and right multiplications with elements of $\GD$ are Zariski continuous on $\GD=\Gq$, we get
\begin{eqnarray*}
  \hat{g} &\in & \overline{ U_{\beta_1}\,\cdots\,U_{\beta_p} T U_{\gamma_1}\,\cdots\,U_{\gamma_q}}\;\;.
\end{eqnarray*}
Because of $U_{\beta_1}\,\cdots\,U_{\beta_p} T U_{\gamma_1}\,\cdots\,U_{\gamma_q}=
U_{\beta_1}\,\cdots\,U_{\beta_p}U_{\gamma_1}\,\cdots\,U_{\gamma_q}T$, and formula (\ref{clSpecm}), we find
\begin{eqnarray*}
  \hat{g} &\in & \overline{ U_{\beta_1}\,\cdots\,U_{\beta_p} U_{\gamma_1}\,\cdots\,U_{\gamma_q}T}^{Specm}\;\;.
\end{eqnarray*}
Next we prove the inclusion '$\supseteq$' of the theorem. Let $\La\in P^+$. Choose $\kBl$-dual bases 
\begin{eqnarray*}
    (a_{\la i})_{\la\in P(\La)\,,\,i=1,\ldots,m_\la} &,&  (b_{\la i})_{\la\in P(\La)\,,\,i=1,\ldots,m_\la}
\end{eqnarray*}
of $L(\La)$, by choosing
$\kBl$-dual bases
\begin{eqnarray*}
   (a_{\la i})_{i=1,\ldots,m_\la} &,&  (b_{\la i})_{i=1,\ldots,m_\la}
\end{eqnarray*}
of every weight space $L(\La)_\la$, $\la\in P(\La)$.\vspace*{1ex}\\
The root groups $U_\al$, $\al\in\rW$, and the torus $T$ act locally finite on $L(\La)$.  Therefore, for every 
$v\in L(\La)$ the linear space spanned by $U_{\beta_1}\cdots U_{\beta_m}Tv$ is finite dimensional. For a fixed pair $(\mu,j)$ we have 
\begin{eqnarray*}
   \kB{a_{\la i}}{U_{\beta_1}\cdots U_{\beta_m}T b_{\mu j}} \;\,=\;\, 0
\end{eqnarray*} 
for all pairs $(\la,i)$ except finitely many. Similarly, for a fixed pair $(\la,i)$ we have
\begin{eqnarray*}
   \kB{a_{\la i}}{U_{\beta_1}\cdots U_{\beta_m}T b_{\mu j}} \;\,=\;\, \kB{T U_{-\beta_m}\cdots U_{-\beta_1} a_{\la i}}{b_{\mu j}}\;\,=\;\,  0
\end{eqnarray*} 
for all pairs $(\mu,j)$ except finitely many.\\
This also applies to $\ti{x}(f_{a_{\la i} b_{\mu j}})$ for every $x\in \overline{U_{\beta_1}\cdots U_{\beta_m}T}^{Specm}$. Now fix such an 
element and define 
\begin{eqnarray*}
   \phi_\La\in End(L(\La)) & \mb{by} & \phi_\La b_{\mu j} \;:=\;\sum_{\la i} b_{\la i}\ti{x}(f_{a_{\la i} b_{\mu j}})\;,\; \mu\in P(\La), 
   \;j=1,\,\ldots,\,m_\mu\;, \\
   \psi_\La\in End(L(\La)) & \mb{by} & \phi_\La b_{\la i} \;:=\;\sum_{\mu j} a_{\mu j}\ti{x}(f_{a_{\la i} b_{\mu j}})\;,\; \la \in P(\La), 
   \;i=1,\,\ldots,\,m_\la\;. 
\end{eqnarray*}
It is easy to check, that $\phi_\La$ and $\psi_\La$ are adjoint maps.\\
The maps $\phi_\La$, $\La\in P^+$, define an element $\phi\in\grAdj$. Note that 
\begin{eqnarray*}
    \ti{x}(f_{a_{\la i} b_{\mu j}})\;=\; \kB{a_{\la i}}{\phi b_{\mu j}}\;\;\la,\mu\in P(\La),\;
    i=1,\,\ldots,\,m_\la,\;j=1,\,\ldots,\,m_\mu,\;\La\in P^+\;\;.
\end{eqnarray*}
From Proposition \ref{V4} a)  we get $\phi\in \overline{\GD}=\GD$. The evaluation homomorphism corresponding to $\phi$ coincides 
with $\ti{x}$. Therefore $\phi=x$.\\
\End
%
%
%
%
\newpage\section{The proof of $Lie(\overline{G})\protect\cong {\bf g}$}
%
%
%
%
In this section we show, that the Lie algebra of $\Gq$ is isomorphic to the Kac-Moody algebra $\g$. For the proof we use, that the 
tangent space of $\Gq$ at 1 is isomorphic to the tangent space of a principal open set $D_{\overline{G}}(\gt_\La)=U^- D_{\overline{T}}(\gt_\La) U^+$, 
$\La\in P^{++}$, at 1. To describe this tangent space, we first determine the Lie algebras of $\Tq$ and $U^\pm$.\vspace*{1ex}\\
To state the propositions and theorems in an easy way, we identify the Lie algebras of $\Tq$, $U^\pm$ and $\Gq$  with 
subalgebras of the Lie algebra of $\grAdj$ (via the tangential maps of the inclusion maps at 1).
This last Lie algebra is identified with $\grAdj$. We also identify the tangent space $T_1(D_{\overline{T}}(\gt_\La))$ with $T_1(\Tq)$ 
(via the tangential map of the inclusion map at 1).\\ 
Similarly for $v\in L(\La)$, we identify the tangent space $T_v({\cal V}_\La)$ of the Kostant cone  
${\cal V}_\La = G(L(\La)_\La)$ with the corresponding subspace of $T_v L(\La)$, which is identified with $L(\La)$.\\
In the proofs however we distinguish carefully between these different possibilities of giving tangent spaces.\vspace*{1ex}\\
The Kac-Moody algebra $\g$ acts faithfully on $\bigoplus_{\La\in P^+}L(\La)$. We also identify $\g$ with the 
corresponding subalgebra of the Lie algebra $\grAdj$.  
%
%
%
\subsection{The Lie algebra of $\Tq$}                              
%
%
%
\begin{Prop}\label{Lie1}                          
We have $\,Lie(\Tq)\,=\,\h\,$.              
\end{Prop}                            
\Proof   
Due to Proposition \ref{V4} the Lie algebra of $\Tq$ is given by    
\begin{eqnarray*}                                                            
   \lefteqn{ Lie(\Tq) \,\;=\;\, \biggl\{\;\,\phi\in\grAdj \;\,\biggl|\biggr.\;\,\exists\;    
  \eps_\phi\,\in\, Der_1(\FK{T})} \qquad\qquad \\    
  & &\;\;\forall\;v,w\in L(\La),\,\La\in P^+\;:\quad    
  \eps_\phi(f_{vw}\res{T}) =\kB{v}{\phi w} \;\,\biggr\}\;\;,    
\end{eqnarray*}                                                                
and $\phi\in Lie(\Tq)$ is uniquely determined by $\eps_\phi\,$. Due to Proposition \ref{Gq1} we have                     
$\FK{T}=\FK{X\cap P}\,$.\vspace*{1ex}\\                                                   
$\bullet$ First we show the inclusion $\h\subseteq Lie(\Tq)$. For $h\in\h$ we get an element 
$\eps_h \in  Der_1(\FK{T})$ by                      
\begin{eqnarray*}    
  \eps_h(e_\la) \;:=\; \la(h) &\quad,\quad & \la\in X\cap P   \;\;. 
\end{eqnarray*}
This derivation has the required properties, because for $v_\la\in L(\La)_\la$, $w_\mu\in L(\La)_\mu$, $\la,\,\mu\in P(\La)$, $\La\in P^+$, we have
\begin{eqnarray*}
  \eps_h(f_{v_\la w_\mu}\res{T}) &=& \eps_h\left(\,\kB{v_\la}{w_\mu}e_\mu\,\right)\;\,=\;\,\kB{v_\la}{w_\mu}\mu(h)\;\,=\;\,\kB{v_\la}{h w_\mu}\;\;.
\end{eqnarray*}                                                                 
$\bullet$ To show the reverse inclusion, we first prove for $\eps\in Der_1(\FK{T})$ the following statement by induction over 
$k\in\N$:\\
Let $\la_1,\ldots, \la_k\in X\cap P$ and $p_1,\ldots, p_k\in\Z$, such that $p_1\la_1+\ldots + p_k\la_k\in X\cap P$. Then we have                          \begin{eqnarray}\label{epsep1pk}                                                      
  \eps(e_{p_1\la_1+\ldots +p_k\la_k}) &=& p_1\eps(e_{\la_1})+\ldots +p_k \eps(e_{\la_k})\;\;.        
\end{eqnarray}                                                
The begin of the induction $k=1$:
If $p_1\geq 0$ this follows by using the derivation properties of $\eps$, and the monoid properties of $X\cap P$. 
If $p_1<0$ then $(-p_1)\la_1\in X\cap P$, and we have     
\begin{eqnarray*}                                                     
  0 \;=\; \eps(e_0) \;=\; \eps(e_{p_1\la_1}e_{-p_1\la_1}) \;=\; \eps(e_{p_1\la_1})+(-p_1)\eps(e_{\la_1})\;\;.                                     
\end{eqnarray*}         
The step of the induction from $k\to k+1$:       
If $p_1,\,\ldots,\,p_k,\,p_{k+1}$ are nonnegative, then $p_1\la_1+\ldots + p_k\la_k\,,\,p_{k+1}\la_{k+1}\in X\cap P$, and we have             
\begin{eqnarray*}                                                              
  \eps(e_{p_1\la_1+\ldots +  p_{k+1}\la_{k+1}}) \;=\;\eps(e_{p_1\la_1 +\ldots + p_k\la_k}e_{p_{k+1}\la_{k+1}})\;=\; 
  \eps(e_{p_1\la_1+\ldots + p_k\la_k})+\eps(e_{p_{k+1}\la_{k+1}})\;\;.                     
\end{eqnarray*}                                                                 
Let one of the $p_i$'s be negative, say $p_{k+1}<0$. Then $(-p_{k+1})\la_{k+1}\in X\cap P$. 
Therefore also $p_1\la_1+\ldots + p_k\la_k\,=\,(p_1\la_1+\ldots +p_{k+1}\la_{k+1})\,+\,(-p_{k+1}\la_{k+1})\,\in\,X\cap P$, and we have           
\begin{eqnarray*}                                               
  \eps(e_{p_1\la_1+\ldots + p_k\la_k}) \;=\;\eps(e_{p_1\la_1+\ldots + p_{k+1}\la_{k+1}}e_{-p_{k+1}\la_{k+1}}) \\         
  \;\,=\;\,\eps(e_{p_1\la_1+\ldots\, p_{k+1}\la_{k+1}})+\eps(e_{-p_{k+1}\la_{k+1}})\;\;.        
\end{eqnarray*}                    
In both cases, using the induction assumption and the begin of the induction, we get equation (\ref{epsep1pk}).\vspace*{1ex}\\                    
Now let $\phi \in Lie(\Tq)$ and set $h:=\sum_{i\,=\,1}^{2n-l}\eps_\phi(e_{\La_i})h_i$. For all $\la\in X\cap P$ we find        
\begin{eqnarray*}                                                             
  \eps_\phi(e_\la) \;\,=\;\, \eps_\phi\left(e_{\sum_i\la(h_i)\La_i}\right) \;\,=\;\, \sum_{i\,=\,1}^{2n-l}\la(h_i)\eps_\phi(e_{\La_i})         
  \;\,=\;\,\la(h) \;\,=\;\, \eps_h(e_\la) \;\;.     
\end{eqnarray*}                                    
Therefore $\eps_\phi=\eps_h$, and $\phi=h\in \h$.     
\End      
%
%
%
\subsection{The Lie algebras of $U^+$ and $U^-$}
%
%
%
To determine the Lie algebras of $U^\pm$ we need the tangent space of the Kostant cone ${\cal V}_\La:=G (L(\La)_\La)$ at a highest weight 
vector.
\begin{Theorem}\label{Lie2}
Let $\La\in P^+$ and $v_\La\in L(\La)_\La$. We have
\begin{eqnarray*}
  T_{v_\La}{\cal V}_\La &=& \left\{\begin{array}{ccc}
  \g L(\La)_\La &\mb{ if } & v_\La\neq 0\\
  L(\La) & \mb{ if } & v_\La=0
  \end{array}\right.\;\;.
\end{eqnarray*}
\end{Theorem}
\Proof
Denote by $L_{high}$ the $L(2\La)$-isotypical component of $L(\La)\otimes L(\La)$. For $\al\in\W\cup\{0\}$ choose ( $|$ )-dual 
bases $ { (e_\al^{(i)}) }_{i=1}^{m_\al}\subseteq \g_\al$, ${ (f_\al^{(i)}) }_{i=1}^{m_\al}\,\subseteq\g_{-\al}$, such that 
$f_\al^{(i)}=e_{-\al}^{(i)}$ for all $\al\in \W$, $i=1,\ldots,\,m_\al$.\vspace*{1ex}\\
$\bullet$ Let $v\in L(\La)$. First we show
\begin{eqnarray}
  T_v{\cal V}_\La &=& \Mklz{x\in L(\La)}{x\otimes v+v\otimes x \in L_{high} }\label{TvVLa1}\\ 
&=& \left\{ x\in L(\La) \:\left|\:  \iB{\La}{\La}\Bigl(x\otimes v + v\otimes x \Bigr)\,=
    \right.\right.\nonumber\\
   && \qquad\quad  \sum_{\al\,\in\,\Delta\cup\{0\}}\sum_i \left.\Bigl(
     e_\al^{(i)}x\otimes f_\al^{(i)}v \, + \,
     e_\al^{(i)}v\otimes f_\al^{(i)}x\Bigr) \:\right\}\;\;.\quad    \label{TvVLa2}
\end{eqnarray}
The equality of the two sets on the right follows, because of
\begin{eqnarray*}
  L_{high} &=& \Mklz{ v\in L(\La)\otimes L(\La)\; }{ \;\tilde{\Omega}\,v \,=\, 
                      \iB{\La}{\La} v\;}\;,
\end{eqnarray*} 
where $\tilde{\Omega}$ is given by the formal expression $\sum_{\al\in\Delta\cup\{0\}} e_\al^{(i)}\otimes f_\al^{(i)}$, compare \cite{K2}, 
Proposition 14.12.\\ 
Due to \cite{KP2}, Theorem 2 a), the vanishing ideal $\Vi{{\cal V}_\La}$ is generated by the functions
\begin{eqnarray*}
  \iB{\La}{\La}f_w f_{w'}\,-\,\sum_{\al\,\in\,\Delta\cup\{0\}}\,
  \sum_{i\,=\,1}^{m_\al}f_{(e_\al^{(i)})^* w}\,f_{(f_\al^{(i)})^* w'} 
  &\quad,\quad & w,w'\,\in\,L(\La)\;\;.
\end{eqnarray*}
Therefore $x\in T_v{\cal V}_\La$ if and only if the corresponding derivation $\de_x\in Der_v\left(\,\FK{L(\La)}\,\right)$ an\-ni\-hi\-lates these 
functions, which means that for all $w,\,w' \in L(\La) $ we have
\begin{eqnarray*}
 0  &=& \iB{\La}{\La}\biggl(\kkB{w}{x}\kkB{v}{w'}+\kkB{w}{v}\kkB{x}{w'}\biggr)\\
    && -\sum_{\al\,\in\,\Delta\cup\{0\}}\sum_i \biggl(
     \kkB{w}{e_\al^{(i)}x}\kkB{f_\al^{(i)}v}{w'} \, + \,
     \kkB{w}{e_\al^{(i)}v}\kkB{f_\al^{(i)}x}{w'} \biggr)\;\;.
\end{eqnarray*}
Because $\kBl$ is nondegenerate, this is equivalent for $x$ to be an element of the set given in (\ref{TvVLa2}).\vspace*{1ex}\\
$\bullet$ Using equation (\ref{TvVLa1}) we find $T_0{\cal V}_\La = L(\La)$. Now let $v_\La\neq 0$.
Since $v_\La\otimes v_\La \in L_{high}$, we get 
\begin{eqnarray*}
 (cy v_\La)\otimes v_\La + v_\La\otimes (cy v_\La) \;=\; 
  cy(v_\La \otimes v_\La)\:\in\: L_{high} &\mb{ for all }& y\in\g,\; c\in\F    \;\;.
\end{eqnarray*}   
Due to equation (\ref{TvVLa1}) we have $T_{v_\La}{\cal V}_\La\supseteq \g L(\La)_\La$. To show the reverse inclusion, let 
$x\in T_{v_\La}{\cal V}_\La $. It is sufficient to show that all the homogenous parts $x_\la\in L(\La)_\la$ of $x$ are contained in 
$\g L(\La)_\La$. Due to equation (\ref{TvVLa2}) we have
\begin{eqnarray*}
  \lefteqn{ \iB{\La}{\La}  \sum_{\la\,\in\, P(\La)} 
  \underbrace{\Bigl( x_\la\otimes v_\La+v_\La\otimes x_\la\Bigr)}_{\in\;
  \Bigl( L(\La)\otimes L(\La)\Bigr)_{\La+\la}} }\\ 
  &=& \sum_{\la\,\in\, P(\La)} \sum_{\al\in\Delta\cup\{0\}}
  \sum_i \underbrace{\left(e_\al^{(i)}x_\la\otimes f_\al^{(i)}v_\La
  +e_\al^{(i)}v_\La\otimes f_\al^{(i)}x_\la\right)}_{\in\;\Bigl(L(\La)\otimes 
  L(\La)\Bigr)_{\La+\la}}\;\;.
\end{eqnarray*}
We compare the $(\La+\la)$-homogenous parts of this equation. Since $f_\al^{(i)}=e_{-\al}^{(i)}$ for 
$\al\in \W$, and $e_\al^{(i)}v_\La=0$ for $\al\in\pW$, we find
\begin{eqnarray*}
  \lefteqn{ \Bigl(\iB{\La}{\La}-\iB{\La}{\la}\Bigr)\Bigl(x_\la\otimes v_\La+ 
  v_\La\otimes x_\la\Bigr)\qquad\qquad } \nonumber\\ 
  & &=\; \sum_{\al\in\Delta^+} \sum_i \left( e_\al^{(i)}x_\la\otimes 
  f_\al^{(i)}v_\La +   f_\al^{(i)}v_\La\otimes e_\al^{(i)}x_\la \right)\;\;.\qquad\qquad
\end{eqnarray*}
If $\iB{\La}{\la}=\iB{\La}{\La}$, then due to \cite{MoPi}, section 6.2, Proposition 9, we have 
$\la=\La$. Therefore $x_\la\in L(\La)_\La$.
If $\iB{\La}{\la}\neq\iB{\La}{\La}$ apply the map $\kB{v_\La}{\,\cdot\:}\otimes id$ to the last equation:
\begin{eqnarray*}
  \lefteqn{ \Bigl(\iB{\La}{\La}-\iB{\La}{\la}\Bigr)\Bigl(\kB{v_\La}{x_\la}v_\La 
   \,+\,\kB{v_\La}{v_\La}x_\la\Bigr) } \nonumber\\ 
  & &=\; \sum_{\al\in\Delta^+} \sum_i \Bigl(
      \kkB{v_\La}{e_\al^{(i)}x_\la}f_\al^{(i)}v_\La \,+\,
      \underbrace{\kkB{v_\La}{f_\al^{(i)}v_\La}}_{=\,0} e_\al^{(i)}x_\la  
      \Bigr)\;\;.\qquad\qquad
\end{eqnarray*}
Solving for the free $x_\la$ gives $x_\la\in\n^-v_\La\oplus\F\,v_\La =\g L(\La)_\La$.\\
\End
\begin{Theorem}\label{Lie3}
We have $\,Lie(U^\pm)\,=\,\n^\pm\,$.
\end{Theorem}
\Proof
Using the tangential map at 1 of the isomorphism of varieties $*:U^-\to U^+$, we find $Lie(U^+)=Lie(U^-)^*$. Therefore
it is sufficient to show the theorem for $U^-$.\vspace*{0.5ex}\\
$\bullet$ First we show $Lie(U^-)\supseteq \n^-$. The Lie algebra $\n^-$ is generated by $\g_\al$, $\al\in\nrW$. Therefore it is 
sufficient to show, that for an element $x\in \g_\al$, $\al\in \nrW$, the corresponding derivation $\de_x\in Der_1\left(\,\FK{\grAdj}\,\right)$ 
annihilates the vanishing ideal $\Vi{U^-}$. If $f\in \Vi{U^-}$, then we have
\begin{eqnarray*}
  0\;=\;f(\exp(tx)) \;=\; f(1)+t\de_x(f) + \mb{O}(t^2) &\mb{ for all }& t\in\F\;\;.
\end{eqnarray*}
Because of $|\F\,|=\infty$ we get $\de_x(f)=0$.\vspace*{0.5ex}\\
$\bullet$ To show the reverse inclusion let $\La\in P^{++}$, $v_\La\in L(\La)_\La\setminus\{0\}$. The map  
\begin{eqnarray*}
  \Phi:\;\;gr\mb{-}Adj\left(\bigoplus_{N\,\in\,P^+}L(N)\right) &\to & L(\La) \\
  \psi\qquad\qquad  &\mapsto & \:\psi v_\La
\end{eqnarray*}
is a morphism of varieties. Its tangential map at 1 is given by:
\begin{eqnarray*}
  T_1\Phi:\;\;T_1\biggl(gr\mb{-}Adj\biggl(\bigoplus_{N\,\in\,P^+}L(N)\biggr)\,\biggr) &\to & T_{v_\La}L(\La) \\
  \phi \qquad\qquad &\mapsto & \;\;\phi v_\La
\end{eqnarray*}
Because of $\Phi(U^-)\subseteq {\cal V}_\La$ we have $(T_1\Phi)(T_1 U^-)\subseteq T_{v_\La}{\cal V}_\La$. Let $\phi\in Lie(U^-)$ and 
$\de_\phi$ the corresponding derivation. Due to the last theorem there exist $y\in\n^-$, $c\in\F$, such that
\begin{eqnarray*}
  \phi v_\La \;=\; yv_\La + cv_\La\;\;.
\end{eqnarray*}
For $\ti{y}\in \n^-$ we find
\begin{eqnarray*}
  \de_{\phi}(f_{\ti{y}v_\La\,v_\La}\res{U^-})  
     \;=\; \kB{\ti{y}v_\La}{yv_\La}+ c\underbrace{\kB{\ti{y}v_\La}{v_\La}}_{=\,0} 
     \;=\; \de_y(f_{\ti{y}v_\La\,v_\La}\res{U^-})\;\;.
\end{eqnarray*}
Because the functions $f_{\ti{y}v_\La\,v_\La}\res{U^-}$, $\ti{y}\in \n^-$ generate the algebra $\FK{U^-}$, we have $\de_\phi = \de_y $. 
Therefore $\phi = y \in \n^-$.\\
\End
%
%
%
%
%
\subsection{The Lie algebra of $\overline{G}$}
%
%
%
To determine the Lie algebra of $\Gq$, we use a description of the variety of the principal open set $D_{\overline{G}}(\gt_\La)=U^- T U^+$, 
$\La\in P^{++}$. We state this description first. 
\begin{Prop}\label{Lie4} Let $\left(\,D_{\overline{T}}(\gt_\La)\,,\,\FK{D_{\overline{T}}(\gt_\La)}\,,\,
{\cal F}_{ D_{\overline{T}}(\gt_\La)}\,\right)$ be the variety of a principal open set $D_{\overline{T}}(\gt_\La)$, $\La\in P^{++}$.
\vspace*{1ex}\\
1) We have $D_{\overline{T}}(\gt_\La) = T$ and $\FK{D_{\overline{T}}(\gt_\La)} = \FK{P}$.\vspace*{1ex}\\
2) We have $T_1 ( D_{\overline{T}}(\gt_\La) ) =\h$. The derivation $\de_h$ corresponding to $h\in\h$ is given by $\de_h(e_\la)=\la(h)$, 
$\la\in P$.
\end{Prop}
\Proof Part 1) can be proved in a similar way as Proposition \ref{Gq11}. Due to our identifications, and Proposition \ref{Lie1}, we have 
$T_1 ( D_{\overline{T}}(\gt_\La) ) =T_1(\Tq)=\h$. The derivation $\eps_h\in Der_1(\FK{\Tq})$ given in the proof of this proposition extends to 
the derivation $\de_h\in Der_1(\FK{D_{\overline{T}}(\gt_\La)})$.\\
\End
\begin{Theorem} Let $\La\in P^{++}$. Let $v_\La\in L(\La)_\La\setminus\{0\}$. The multiplication map
\begin{eqnarray*}
  m:\; U^-\times  D_{\overline{T}}(\gt_\La)\times U^+ &\to & D_{\overline{G}}(\gt_\La) 
\end{eqnarray*}
is bijective. Its comorphism 
\begin{eqnarray*}
  m^*:\; \FK{D_{\overline G}(\gt_\La)} &\to & \FK{U^-}\otimes\FK{P}\otimes\FK{U^+}
\end{eqnarray*}
exists, and is an isomorphism of algebras. Furthermore we have
\begin{eqnarray}
  m^*\left(\,\frac{f_{yv_\La\,v_\La}}{\gt_\La}\res{D_{\overline G}(\gt_\La)}\,\right) 
    &=& f_{yv_\La\,v_\La}\res{U^-}\otimes\,1\,\otimes\,1 
    \quad,\quad y\in \n^- \,, \label{LetzteGleichung1}\\
  m^*\left(\,\frac{\gt_N}{(\gt_\La)^n}\res{D_{\overline G}(\gt_\La)} \,\right) 
    &=& 1\,\otimes e_{N-n\La}\otimes\,1 
    \quad,\quad  N\in P^+,\,n\in\Nn \,,\label{LetzteGleichung2}\\
  m^*\left(\,\frac{f_{v_\La\,yv_\La}}{\gt_\La}\res{D_{\overline G}(\gt_\La)} \,\right) 
    &=& 1\,\otimes\,1\,\otimes f_{v_\La\,yv_\La}\res{U^+}
    \quad,\quad y\in \n^- \,. \label{LetzteGleichung3}
\end{eqnarray}
The map $m^{-1}: D_{\overline G}(\gt_\La) \to  U^-\times  D_{\overline{T}}(\gt_\La) \times U^+$ is a morphism of varieties. Its tangent map
$T_1(m^{-1}):\,T_1(D_{\overline G}(\gt_\La)) \to  T_{(1,1,1)}\left( U^-\times D_{\overline{T}}(\gt_\La)) \times U^+ \right)$ is injective.
\end{Theorem}
\Proof The proof is similar to the proof of Theorem \ref{Gq12}. The injectivity of $T_1(m^{-1})$ follows from the surjectivity of $(m^{-1})^*$.\\
\End
The following theorem is one of the main results of this paper:
\begin{Theorem}\mb{}
We have $\,Lie(\Gq)\,=\,\g\,$.
\end{Theorem}
\Proof\\
$\bullet$ First we show $Lie(\Gq)\supseteq\g$. It is sufficient to show $\g_\beta\subseteq Lie(\Gq)$ for all $\beta\in\rW$, 
and $H\subseteq Lie(\Gq)$. Due to Proposition \ref{V4} the Lie algebra of $\Gq$ is given by    
\begin{eqnarray*}                                                            
   \lefteqn{ Lie(\Gq) \,\;=\;\, \biggl\{\;\,\phi\in\grAdj\;\,\biggl|\biggr.\;\,\exists\;    
  \de_\phi\,\in\, Der_1(\FK{G})} \qquad\qquad \\    
  & &\;\;\forall\;v,w\in L(\La),\,\La\in P^+\;:\quad    
  \de_\phi(f_{vw}\res{G}) =\kB{v}{\phi w} \;\,\biggr\}   \;\;. 
\end{eqnarray*}      
Let $x_\beta\in\g_\beta$, $\beta\in\rW$, and $h\in H$. Due to the Peter and Weyl theorem for $\FK{G}$, there exist linear maps 
\begin{eqnarray*}
\begin{array}{c}
  \de_{x_\beta}:\FK{G}\to\F \\
  \de_{h}:\FK{G}\to\F  
\end{array}\mb{ with }\begin{array}{c}
  \de_{x_\beta}(f_{vw}\res{G}):=\kB{v}{x_\beta\,w} \\
  \de_{h}(f_{vw}\res{G}):=\kB{v}{h\,w}
\end{array} ,\,\;v,w\in L(\La),\;\La\in P^+\:.
\end{eqnarray*}
Using the 1-parameter subgroups $exp(tx_\beta)$, $t\in \F$, and $t_h(s)$, $s\in \F^\times$, we can write 
these maps in the following forms, which show that $\de_{x_\beta}$ and $\de_h$ are derivations in $1$:
\begin{eqnarray*}
\begin{array}{c}
  \de_{x_\beta}(f)\;=\;\frac{d}{dt}\res{0} f(exp(tx_\beta)) \\
  \de_{h}(f)\;=\;\frac{d}{ds}\res{1}f(t_h(s))
\end{array} \;\;\:,\:\;\;f\in\FK{G}\;\;.
\end{eqnarray*}
$\bullet$ To prove the inclusion $Lie(\Gq)\subseteq \g$, we use the notations of the last theorem. We identify $T_{(1,1,1)}\left( U^-\times D_{\overline{T}}(\gt_\La)) \times U^+ \right)$ with $T_1 (U^-)\times T_1 (D_{\overline T}(\gt_\La))\times T_1 (U^+)$.\\
The tangent map $T_1\, j:T_1(D_{\overline G}(\gt_\La))\to T_1\Gq $ of the inclusion map $ j:D_{\overline G}(\gt_\La)\to\Gq$
is bijective. Therefore the following concatenation of maps
\begin{eqnarray*}
     T_1(\Gq)\;\,\stackrel{(T_1 j)^{-1}}{\longrightarrow} \;\,T_1\left(D_{\overline G}(\gt_\La)\right) \;\,\stackrel{T_1 (m^{-1})}{\longrightarrow} \;\, 
  T_1 (U^-)\times T_1 (D_{\overline T}(\gt_\La))\times T_1 (U^+)
\end{eqnarray*}
is injective. Due to Theorem \ref{Lie3} we have $T_1(U^\pm)=\n^\pm$. Due to Proposition \ref{Lie4} we have $T_1 (D_{\overline T}(\gt_\La))=\h$. 
We also know $\g\subseteq T_1(\Gq)$. To show equality, it is sufficient to show 
\begin{eqnarray*}
\left(T_1(m^{-1})\circ (T_1 j)^{-1}\right)(y+h+x)\;\,=\;\,(\,y\,,\,h\,,\,x\,) &\mb{ for all }& y\in\n^-,\;  h\in\h,\;  x\in \n^+  \;\;.
\end{eqnarray*}
To do this, we have to work with the corresponding derivations. Set $z=y+h+x$, and let $\de_z\in Der_1 (\,\FK{\Gq}\,)$ be the corresponding 
derivation. Let $\de_y\in Der_1 (\,\FK{U^-}\,)$, $\de_h\in Der_1 (\,\FK{D_{\overline T}(\gt_\La)}\,)$, and 
$\de_x\in Der_1 (\,\FK{U^+}\,)$ be the derivations corresponding to $y$, $h$, and $x$.\\
Using equation (\ref{LetzteGleichung3}) of the last theorem, we find for $\ti{y}\in \n^-$:
\begin{eqnarray*}
  \lefteqn{ \Bigl(\left(T_1 (m^{-1})\right)(T_1 j)^{-1} (\de_{z})\Bigr) (\,1\,
  \otimes\,1\,\otimes f_{v_\La\,\ti{y}v_\La}\res{U^+})}\\ 
  &=& \Bigl((T_1 j)^{-1}(\de_{z})\Bigr)\left(
  \frac{f_{v_\La\,\ti{y}v_\La}}{\gt_\La}\res{D_{\overline G}(\gt_\La)}\right) \\
  &=& \frac{1}{\gt_\La(1)^2}\,\Bigl(\,\de_{z}(f_{v_\La\,\ti{y}v_\La}
  \res{\overline G})\underbrace{\gt_\La(1)}_{=\,1}\,-\,
  \underbrace{f_{v_\La\,\ti{y}v_\La}(1)}_{=\,0}\de_{z}(\gt_\La)\,\Bigr) \\
  &=& 
  \de_x(f_{v_\La\,\ti{y}v_\La}\res{U^+})\;\;.
\end{eqnarray*} 
Because the algebra $\FK{U^+}$ is generated by $f_{v_\La\,\ti{y}v_\La}\res{U^+}$, $\ti{y}\in\n^-$, we get
\begin{eqnarray*}
 \left(T_1(m^{-1})\right)(T_1 j)^{-1}(\de_{z}) &=& (\,\ldots\,,\,\ldots\,,\,\de_x\,)\;\;.
\end{eqnarray*}
Similar we have $\,\left(T_1(m^{-1})\right)(T_1 j)^{-1}(\de_{z}) \,=\, (\,\de_{y}\,,\,\ldots\,,\,\ldots\,)\,$.\vspace*{0.5ex}\\
Let $N \in P^+$ and $n\in\Nn$. Using equation (\ref{LetzteGleichung2}) of the last theorem we find
\begin{eqnarray*}
  \lefteqn{  \Bigl(\left(T_1 (m^{-1})\right)(T_1 j)^{-1})(\de_{z})\Bigr) (\,1\,\otimes e_{N-n\La}\otimes \,1\,)  
 \;\,=\;\,\Bigl((T_1 j)^{-1}(\de_{z})\Bigr)\left(\frac{\gt_N}{(\gt_\La)^n}\res{D_{\overline G}(\gt_\La)}\right) }\\
  &=& \frac{1}{\gt_{n\La}(1)^2}\,\Bigl(\,\de_{z}(\gt_N)\gt_{n\La}(1)-
  \de_{z}(\gt_{n\La})\gt_N(1)\,\Bigr) 
  \;\,=\;\,N(h)- n\La(h)\quad\qquad\qquad\qquad\qquad\\
  &=&\de_h(e_{N-n\La})\;\;.
\end{eqnarray*}
Because $\FK{P}$ is spanned by $e_{N-n\La}$, $N\in P^+$, $n\in\Nn$, we conclude
\begin{eqnarray*}
 \left(T_1(m^{-1})\right)(T_1 j)^{-1}(\de_{z}) &=& (\,\ldots\,,\,\de_h\,,\,
 \ldots\,)\;\;.
\end{eqnarray*}
\End
%
%
%

%
\end{document}